\documentclass[12pt,twoside, reqno]{amsart}

\usepackage{amssymb}
\usepackage{amsmath}
\usepackage{mathrsfs}
\usepackage{amsthm}
\usepackage{amsfonts}
\usepackage{txfonts}
\usepackage{enumerate}
\usepackage{hyperref}
\hypersetup{hidelinks}
\usepackage{graphicx}
\usepackage{subfigure}
\usepackage{float}
\usepackage{latexsym,amssymb}
\usepackage{color}
\usepackage{graphicx}
\usepackage{pgf,tikz}
\usepackage{mathrsfs}
\usepackage{fancyhdr}
\usetikzlibrary{arrows}

\allowdisplaybreaks
\usepackage{indentfirst}
\usepackage{setspace}
\usepackage{geometry}
\newtheorem{theorem}{Theorem}[section]
\newtheorem{proposition}[theorem]{Proposition}
\newtheorem{corollary}[theorem]{Corollary}
\newtheorem{lemma}[theorem]{Lemma}
\newtheorem{remark}[theorem]{Remark}

\numberwithin{equation}{section} 

\numberwithin{equation}{section}

\newcommand{\R}{\mathbb{R}}

\newcommand{\ben}{\begin{eqnarray*}}
\newcommand{\enn}{\end{eqnarray*}}
\newcommand{\pa}{\partial}

\newcommand{\g}{\gamma}

\newcommand{\ve}{\varepsilon}

\newcommand{\ol}{\overline}

\newcommand{\half}{\frac{1}{2}}

\newcommand{\na}{\nabla}
\newcommand{\be}{\begin{equation}}
\newcommand{\ee}{\end{equation}}
\newcommand{\ba}{\begin{aligned}}
\newcommand{\ea}{\end{aligned}}
\newcommand{\lf}{\left}
\newcommand{\rt}{\right}

\newcommand{\wt}{\tilde{w}}

\def\9{{\infty}}
\def\a{{\alpha}}
\def\b{{\beta}}
\def\g{{\gamma}}

\def\lbb{{\lambda}}
\def\t{{\theta}}
\def\call{{\mathcal{L}}}

\def\calo{{\mathcal{O}}}
\def\calp{{\mathcal{P}}}

\def\calu{{\mathcal{U}}}

\def\bbr{{\mathbb{R}}}
\def\ve{{\varepsilon}}
\def\vf{{\varphi}}

\def\wt{\widetilde}
\def\wh{\widehat}
\def\ol{\overline}
\def\({\left(}
\def\){\right)}
\def\<{\langle}
\def\>{\rangle}

\begin{document}

\title
[uniqueness of multi-bubble blow-up solutions and multi-solitons to NLS]{On uniqueness of multi-bubble blow-up solutions and multi-solitons
to $L^2$-critical nonlinear Schr\"odinger equations}

\author{Daomin Cao}
\address{Institute of Applied Mathematics, Chinese Academy of Sciences, Beijing 100190, and University of Chinese Academy of Sciences, Beijing 100049, P.R. China}
\email{dmcao@amt.ac.cn}
\thanks{}
\author{Yiming Su}
\address{Department of mathematics,
Zhejiang University of Technology, 310014 Zhejiang, China}
\email{yimingsu@zjut.edu.cn}
\thanks{}

\author{Deng Zhang}
\address{School of mathematical sciences,
Shanghai Jiao Tong University, 200240 Shanghai, China }
\email{dzhang@sjtu.edu.cn}
\thanks{}



\date{}


\begin{abstract}
We are concerned with the focusing $L^2$-critical nonlinear Schr\"odinger equations
in $\bbr^d$ for $d=1,2$.
The uniqueness is proved for a large energy class of multi-bubble blow-up solutions,
which converge to a sum of $K$ pseudo-conformal blow-up solutions
particularly with low rate $(T-t)^{0+}$,
as $t\to T$, $1\leq K<\9$.
Moreover,
we also prove the uniqueness in the energy class of multi-solitons
which converge to a sum of $K$ solitary waves with convergence rate $(1/t)^{2+}$,
as $t\to \9$.
The uniqueness class is further enlarged to contain the multi-solitons
with even lower convergence rate $(1/t)^{\frac 12+}$
in the pseudo-conformal space.
The proof is mainly based on the pseudo-conformal invariance and the monotonicity properties
of several functionals adapted to the multi-bubble case, the latter is crucial towards the
upgradation of the convergence
to the fast exponential decay rate.
\end{abstract}

\maketitle
\begin{spacing}{1.02}
\tableofcontents
\end{spacing}

\section{Introduction and formulation of main results} \label{Sec-Intro}

\subsection{Introduction} \label{Subsec-Intro}

We are concerned with the focusing $L^2$-critical nonlinear Schr\"odinger equations
in $\bbr^d$ for $d=1,2$,
\begin{align}  \label{equa-NLS}
 &i\partial_tu+\Delta u+|u|^\frac{4}{d}u=0,  \tag{NLS}
\end{align}
where $u:I\times \bbr^d \to \mathbb{C}$, $I\subseteq \bbr$ is a time interval.

Equation \eqref{equa-NLS} has a variety of applications
in nonlinear optics, Bose-Einstein condensation and plasma physics.
It can be a model for the propagation of intense laser beams in bulk media with Kerr nonlinearity,
and has also relationship with the weak turbulence theory.
See, e.g., \cite{DNPZ92,GS11,SS99}.

Mathematically, it is well-known (see, e.g., \cite{C03,W83}) that
equation \eqref{equa-NLS} is locally well-posed in the space $H^1$.
An important role is played by the {\it ground state}
$Q$, which is the unique positive radial solution to the elliptic equation
\begin{align} \label{equa-Q}
    \Delta Q - Q + Q^{1+\frac{4}{d}} =0.
\end{align}
The mass of ground state is exactly the threshold of
global well-posedness and blow-up of solutions to \eqref{equa-NLS}:
the solutions exist globally in the subcritical mass case where $\|u\|_{L^2} < \|Q\|_{L^2}$,
while singularities can be formed in the (super-)critical mass case
where $\|u\|_{L^2} \geq \|Q\|_{L^2}$.

Moreover, the solutions  to \eqref{equa-NLS} satisfy
the following three conservation laws:
\begin{enumerate}
  \item[$\bullet$] Mass conservation:
  \begin{align} \label{mass-conserv}
     M(u)(t):=\int_{\R^d}|u(t)|^2dx=M(u_0).
  \end{align}
  \item[$\bullet$] Energy conservation:
  \begin{align}   \label{energy-conserv}
  E(u)(t) := \half\int_{\R^d}|\nabla u(t)|^2dx-\frac{d}{2d+4}\int_{\R^d}|u(t)|^{2+\frac{4}{d}}dx=E(u_0).
  \end{align}
  \item[$\bullet$] Momentum conservation:
  \begin{align}  \label{momentum-conserv}
  Mom(u)(t):={\rm Im} \int_{\R^d}\nabla u(t) \ol{u}(t)dx =  Mom(u_0).
  \end{align}
\end{enumerate}
Equation \eqref{equa-NLS} also admits the invariance under the translation, scaling, phase rotation and Galilean transform,
i.e.,
if $u$ solves (\ref{equa-NLS}),
then so does
\begin{align}\label{symmerty}
 \mathscr{T}u(t,x)=\lambda_0^{-\frac{d}{2}}
u \(\frac{t-t_0}{\lambda_0^2},\frac{x-x_0 }{\lambda_0}- \frac{\beta_0(t-t_0)}{\lbb_0}\)
e^{i\frac{\beta_0}{2}\cdot(x-x_0)-i\frac{|\beta_0|^2}{4}(t-t_0)+i\t_0},
\end{align}
where
$(\lambda_0, \beta_0, \theta_0) \in \bbr^+ \times \bbr^d \times \bbr$,
$x_0\in\R^d$, $t_0\in\R$.
In particular,
the $L^2$-norm of solutions is  preserved under the scaling,
and thus
\eqref{equa-NLS} is  called the {\it $L^2$-critical} equation.

Another invariance,
particularly important in the $L^2$-critical case,
is the {\it pseudo-conformal invariance}
related to the
{\it pseudo-conformal transformation}, defined by
\begin{align} \label{pseu-conf-transf}
   P_T(u)(t,x):= \frac{1}{(T-t)^{\frac d2}} u \(\frac{1}{T-t}, \frac{x}{T-t}\) e^{-i\frac{|x|^2}{4(T-t)}}, \ \ t\not =T,\ u\in \Sigma,
\end{align}
where $\Sigma$ denotes the pseudo-conformal space
$\Sigma:=\{u\in H^1: \|u\|_{\Sigma}:= \|u\|_{H^1} + \|xu\|_{L^2}<\9\}$.

Furthermore,
the pseudo-conformal invariance relates both multi-solitons
and multi-bubble blow-up solutions with
pseudo-conformal blow-up rate.
These two special families of solutions
are indeed of significant importance to
describe the dynamics of solutions to \eqref{equa-NLS}:
one in the large time behavior regime,
and the other in the singularity regime.
Hence,
the pseudo-conformal invariance
provides an alternative way to study the uniqueness of multi-solitons
from that of multi-bubble blow-up solutions,
to which more sharper singularity analysis can be performed.
This is actually the main motivation of present work.

To be precise,
on one hand,
given $K \in \mathbb{N}\setminus\{0\}$,
the {\it solitary waves} $W_k$, $1\leq k\leq K$, are defined by
\begin{align}  \label{Wj-soliton}
W_k(t,x):=\omega_k^{-\frac d2}Q \(\frac{x-v_k t}{\omega_k} \)e^{i(\half v_k\cdot x-\frac{1}{4}|v_k|^2t+\omega_k^{-2}t+\vartheta_{k})},
\end{align}
where the parameters $\omega_k \in \bbr^+$, $v_k\in \bbr^d$ and $\vartheta_{k}\in \bbr$,
corresponding to the frequency,
propagation speed
and phase, respectively, $1\leq k\leq K$.
A {\it multi-soliton} (or, {\it multi-solitary wave solution})
is a solution to \eqref{equa-NLS} defined on $[T_0, \9)$ for some $T_0 \in \bbr$
and such that
\begin{align} \label{u-W-H1-o1}
 \|u(t) - \sum\limits_{k=1}^K W_k(t)\|_{H^1} =o(1),\ \ as\ t\ \to \9,
\end{align}
where and hereafter $o(1)$ means small quantities that converge to zero.
This means that the multi-solitons behave
exactly as a sum of  solitary waves
without loss of mass by dispersion.
We would like to mention that,
multi-solitions with space translations are also studied in literature,
see, e.g., \cite{CF20,CL11,CMM11,MM06,M90}.
For the simplicity of exposition,
we focus on the multi-solitons of form \eqref{Wj-soliton}.

Multi-solitons have been studied extensively  in the integrable case,
see, e.g.,  \cite{M76,ZS72}.
For the nonintegrable equations,
the construction of multi-solitons to \eqref{equa-NLS}
was initiated by Merle \cite{M90}.
The proof in  \cite{M90} is based on the pseudo-conformal invariance and
the construction of multi-bubble blow-up solutions.
Later,
multi-solitons are constructed in the subcritical case
by Martel and Merle \cite{MM06},
and in the supercritical case by C\^ote, Martel and Merle \cite{CMM11}. Moreover, Martel, Merle and Tsai \cite{MMT06} proved that  multi-solitons  are stable by the perturbation of  initial data in the energy space.
We also refer to \cite{LeLP15,LeT14}
for the construction of the infinite trains of solitons.

Multi-solitons have been also constructed in various other settings.
For the generalized Korteweg-de Vries (gKdV) equations,
we refer to \cite{Ma05} for the subcritical and critical cases,
and \cite{Co11} for the super-critical case.
See also \cite{CM14} for the Klein-Gordon equation,
\cite{K-M-R} for the Hartree equation
and \cite{MRT15} for the water-waves system.

Despite the extensive study of constructions,
it remains open for
the uniqueness or classification of multi-solitons to
nonlinear Schr\"odinger equations.
This was first pointed out
by Martel and Merle \cite{MM06}
in the $L^2$-subcritical case,
and later was raised as an open problem by Martel \cite{M18}
in both the $L^2$-subcritical and critical cases.
It is also expected that no uniqueness holds
in the $L^2$-supercritical case,
see C\^ote and Le Coz \cite{CL11}.

To our knowledge,
the only complete study of the uniqueness problem of multi-solitons
was done for the gKdV equations
in the pioneering work of Martel \cite{Ma05} for the $L^2$-subcritical
and critical cases.
Later, multi-solitons were classified by C\^ombet \cite{Co11}
in the $L^2$-supercritical case.
An important ingredient in the proof of \cite{Ma05,Co11}
is the almost monotonicity of local mass and energy,
which allows to gain the fast exponential decay rate.

However, this property fails for NLS,
which makes it quite difficult to study the uniqueness  of multi-solitons to NLS.
The recent progress has been made by C\^ote and Friederich \cite{CF20}
in the $L^2$-subcritical and critical cases.
In \cite{CF20},
the uniqueness class is achieved for the multi-solitons
$u$ to NLS that converge to the sum of solitary waves with a high power of $1/t$, i.e.,
\begin{align} \label{u-W-H1-N}
    \|u(t) -\sum\limits_{k=1}^K W_k(t) \|_{H^1} = \calo\(\frac{1}{t^{N}}\), \ \ as\ t\to \9,
\end{align}
for some $N$ large enough.
Note that,
this result allows to break the class of exponential convergence.
Moreover,
the smoothness of multi-solitons and
general nonlinearities are also studied in  \cite{CF20}.
The challenge is hence to prove the uniqueness of multi-solitons
in the low convergence regime.

Furthermore,
the above uniqueness problem in the single bubble case ($K=1$)
is closely related to the {\it solitary wave conjecture},
which states that non-scattering $L^2$ global solutions to equation \eqref{equa-NLS}
with critical mass $\|Q\|^2_{L^2}$
shall coincide with the solitary wave up to the symmetries \eqref{symmerty}.
This conjecture is affirmative in the pseudo-conformal space
by the rigidity result of Merle \cite{M93}
and the pseudo-conformal invariance.
In the Sobolev space,
when $d =2,3$ it is proved in \cite{LZ12} for $H^1$ radial solutions,
and when $d\geq 4$
it is proved in \cite{KLVZ09} and \cite{LZ09},
respectively, for $H^1$ and $H^s$ ($s>0$) radial solutions.

On the other hand,
given $T \in \bbr$,
the {\it  pseudo-conformal blow-up solutions}
$S_k$, $1\leq k\leq K$, are defined by
\begin{align}  \label{Sj-blowup}
S_k(t,x):=(\omega_k(T-t))^{-\frac d2} Q \(\frac{x-x_k}{\omega_k(T-t)}\) e^{-\frac i4\frac{|x-x_k|^{2}}{T-t}+\frac{i}{\omega_k^2(T-t)}+i\vartheta_{k}}.
\end{align}
Note that,
$S_k$ has the critical mass $\|Q\|^2_{L^2}$
and blows up at time $T$ with the blow-up rate $\|\nabla S_k(t)\|_{L^2}\sim(T-t)^{-1}$.
More importantly,
by the seminal work of Merle \cite{M93},
the pseudo-conformal blow-up solutions
are exactly the unique minimal mass blow-up solutions to \eqref{equa-NLS},
up to the symmetries \eqref{symmerty}.
Furthermore,
via the
the pseudo-conformal transformation,
the pseudo-conformal blow-up solutions
are  closely related to the solitary waves
\begin{align}
    S_k = P_T (W_k), \ \ with\ x_k = v_k,\ 1\leq k\leq K.
\end{align}

Thus, a natural question, as in the case of multi-solitons with asymptotic behavior \eqref{u-W-H1-o1},
is whether there exists a unique {\it multi-bubble blow-up solution} $v$ to \eqref{equa-NLS}
such that
\begin{align} \label{v-S-H1-o1}
    \|v(t) -\sum\limits_{k=1}^K S_k(t) \|_{H^1} =o(1),\ \ as\ t\to T.
\end{align}

It should be mentioned that,
the multi-bubble blow-up solutions to \eqref{equa-NLS}
with pseudo-conformal blow-up rate were first constructed by Merle \cite{M90}.
Bubbling phenomena have been also
exhibited in various other settings.
We would like to refer to
\cite{F17} for the multi-bubble solutions with log-log blow-up rate,
\cite{J17} for  the energy-critical NLS,
and \cite{MP18} for the blow-up solutions with multiple
bubbles concentrating at the same point.
We also refer to \cite{CM18,L17} for the gKdV equations,
\cite{JL18,KST18} for the wave maps,
and \cite{SG19} for the nonlinear Schr\"{o}dinger system.

The understanding of above uniqueness problem
enables to enlarge the uniqueness class of multi-solitons to \eqref{equa-NLS}
particularly in the low convergence regime.

Furthermore,
this uniqueness problem provides an illustration of the rigidity of equation \eqref{equa-NLS}
around the pseudo-conformal blow-up solutions.
It seems also helpful to understand the {\it mass quantization conjecture}
in \cite{MR05}.
It is conjectured by Merle and Rapha\"el \cite{MR05} that,
the blow-up solutions to \eqref{equa-NLS} shall concentrate
the mass $m_k(\geq \|Q\|_{L^2}^2)$ at the singularities
and converge strongly in $L^2$ to a residue $u^*$ away from the singularities.
Then, intuitively,
the blow-up solutions with mass $K\|Q\|_{L^2}^2$ and $K$ distinct singularities
shall distribute the mass $\|Q\|_{L^2}^2$ to each blow-up bubble,
and thus the residue vanishes (i.e., $u^* = 0$).
Inspired by the rigidity result in \cite{M93},
each blow-up bubble with critical mass $\|Q\|_{L^2}^2$
is expected to behave as a pseudo-conformal blow-up solution near the corresponding singularity.
Thus, the above intuition leads naturally to the multi-bubble blow-up solutions
with asymptotic behavior \eqref{v-S-H1-o1}
(and an extra energy information),
and the uniqueness problem states that
this class of blow-up solutions shall be unique.

Let us also mention that,
this kind of uniqueness problem shares interesting similarities with the local uniqueness problem of peak or bubbling solutions
to nonlinear elliptic equations, which has attracted a lot of attentions during the last decades.
For instance,
Cao and Heinz \cite{CH03} studied the local uniqueness of multi-lump solutions
$u_\varepsilon$ to stationary nonlinear Schr\"{o}dinger equations  such that
\begin{align*}
\|u_\varepsilon-\sum_{k=1}^{K}Q_k(\frac{x-x_{k,\varepsilon}}{\varepsilon})\|_{H^1}
=\calo(\varepsilon^{\frac{d}{2}+2}), \ \ and  \ x_{k,\varepsilon}\rightarrow x_k,\ \ as \ \ \varepsilon\rightarrow0,
\end{align*}
where $\{x_1,\dots,x_K\}$ are the distinct nondegenerate critical points
of potential $V$,
$Q_k(x)=Q(\sqrt{V(x_k)}x)$.
See also \cite{CLP20} for
the local uniqueness of multi-peak solutions $u_\varepsilon$ to Brezis-Nirenberg problem  concentrating at different points $\{x_1,\dots,x_K\}$ satisfying
\begin{align*}
\|u_\varepsilon-\sum_{k=1}^{K}PU_{x_{k,\varepsilon},\lambda_{k,\varepsilon}}\|_{H^1}=o(1),\  \ x_{k,\varepsilon}\rightarrow x_k, \ and \ \lambda_{k,\varepsilon}\rightarrow+\infty,\ \ as \ \ \varepsilon\rightarrow0,
\end{align*}
where $P$ is the projection from $H^1(\Omega)$ onto $H^1_0(\Omega)$,
$\Omega$ is a smooth and bounded domain in $\bbr^d$,
and
$U_{x,\lambda}$
solves the elliptic equation $\Delta u + u^{\frac{d+2}{d-2}}=0$ in $\bbr^d$, $d\geq 5$.
We refer  to \cite{CGPY19,CPY20} and the references therein for more recent progresses on local uniqueness problem
for other elliptic equations.

As mentioned above,
the unconditional uniqueness  of critical mass blow-up solutions to \eqref{equa-NLS}
has been completely understood by the seminal work of Merle \cite{M93}.
Such strong rigidity result has been also proved for
the minimal mass blow-up solutions to
the inhomogeneous NLS by Rapha\"el and Szeftel \cite{RS11},
under the sharp non-degenerate condition of nonlinearity.
In particular,
a robust modulation method and upgradation procedures
to upgrade the estimates of remainder have been developed in \cite{RS11},
particularly
in the absence of the pseudo-conformal invariance.
The uniqueness of minimal mass blow-up solutions
has been also completely studied for the $L^2$-critical gKdV equation by
Martel, Merle and Rapha\"el \cite{MMR15}.
For the conditional uniqueness results in the single bubble case,
we would like to refer to \cite{MRS13} for the Bourgain-Wang blow-up solutions to NLS,
and \cite{KK20} for the Chern-Simons-Schr\"odinger equations.

In the recent work \cite{SZ20},
we study the multi-bubble blow-up solutions to stochastic nonlinear Schr\"odinger equations driven by Wiener processes,
in the absence of the pseudo-conformal invariance and the conservation law of energy
due to the presence of noise.
The multi-bubble blow-up solutions are constructed
to converge exponentially fast in $\Sigma$
to a sum of $K$ pseudo-conformal blow-up solutions, $1\leq K<\9$.
Moreover, the conditional uniqueness  is proved in \cite{SZ20} for the class of multi-bubble blow-up solutions
$v$ such that
\begin{align} \label{v-S-H1-3+}
    \|v(t) - \sum\limits_{k=1}^K S_k(t)\|_{H^1}
    = \calo((T-t)^{3+}), \ \ as\ t\to T,
\end{align}
provided that the frequencies $\{\omega_k\}$ are close
or the relative distances $|x_j - x_k|$ are large, $j\not= k$
(see {\rm Case (I)} and {\rm Case (II)} below).
Here $(T-t)^{3+}$ means that $(T-t)^{3+\zeta}$ for any given $\zeta >0$.
In particular, these results are applicable to equation \eqref{equa-NLS},
without the driven noise.
By virtue of the pseudo-conformal invariance,
the  above conditional uniqueness result also gives the uniqueness of multi-solitons to \eqref{equa-NLS} either
with convergence rate $(1/t)^{5+}$ in the energy space
or with rate $(1/t)^{4+}$ in the pseudo-conformal space.
See Remark \ref{Rem-uniq-blow-H13+} $(iii)$ below.

In this paper,
we obtain the uniqueness in a large energy class of multi-bubble blow-up solutions $v$ to \eqref{equa-NLS}
particularly in the low asymptotic regime where
\begin{align} \label{v-S-H1-0+}
    \|v(t) -\sum\limits_{k=1}^K S_k(t)  \|_{H^1} =\calo((T-t)^{0+}), \ \ as\ t\to T,
\end{align}
in both  {\rm Case (I)} and {\rm Case (II)}.
Note that, the convergence rate in \eqref{v-S-H1-0+} almost reaches
the one $o(1)$ in \eqref{v-S-H1-o1}.
The  condition \eqref{v-S-H1-0+} can be even weakened
by the asymptotic behavior
\begin{align} \label{v-S-L2-H1-o1}
    \|v(t) -\sum\limits_{k=1}^K S_k(t)  \|_{L^2}
    + (T-t)  \|\na v(t) -\sum\limits_{k=1}^K \na S_k(t)  \|_{L^2}  = o(1), \ \ as\ t\to T,
\end{align}
plus additionally a double average condition (see \eqref{nav-naS-iint-0+-Thm} below).
In both cases, we show that the convergence rate indeed can be upgraded to the
much faster exponential decay rate
$e^{- \frac{\delta}{T-t}}$, $\delta>0$,
in the more regular pseudo-conformal space $\Sigma$.
In particular,
these results apply to the case where
the blow-up points $\{x_k\}$ can be arbitrarily distinct when
the frequencies are the same.

Furthermore,
the above results also allow to enlarge the uniqueness class of multi-solitons to \eqref{equa-NLS}.
By virtue of the pseudo-conformal invariance,
we are able to return back to multi-solitons
and prove that the uniqueness holds in the energy class of multi-solitons $u$ to \eqref{equa-NLS}
such that
\begin{align} \label{u-W-H1-2+}
   \|u(t) - \sum\limits_{k=1}^K W_k(t) \|_{H^1} = \calo\(\frac{1}{t^{2+}}\),\ \ for\ t\ large\ enough,
\end{align}
in both  {\rm Case (I)} and {\rm Case (II)},
where the blow-up points $\{x_k\}$ are replaced by the speeds $\{v_k\}$.
Again, the convergence rate in \eqref{u-W-H1-2+} can be upgraded to the exponential decay rate
in $\Sigma$.
One interesting application
is that
the speeds $\{v_k\}$ can be arbitrarily distinct
when the frequencies $\{\omega_k\}$ are the same.

In the single soliton case (i.e., $K=1$),
this result in particular shows that
the solitary wave conjecture is affirmative
for general $H^1$ solutions
with asymptotic behavior \eqref{u-W-H1-2+}.
Additionally, if the propagation speed is zero,
the uniqueness class can be further enlarged to contain the solutions
with lower convergence rate $(1/t)^{1+}$.

Moreover,
we also prove the uniqueness of multi-solitons $u$ to \eqref{equa-NLS}
with even lower convergence rate $(1/t)^{\frac 12+}$
in the pseudo-conformal space, i.e.,
\begin{align} \label{u-W-Sigma-12+}
   \|u(t) - \sum\limits_{k=1}^K W_k(t) \|_{\Sigma} = \calo\(\frac{1}{t^{\frac 12+}}\),\ \ for\ t\ large\ enough,
\end{align}
provided additionally that the speeds $\{v_k\}$ are non-zero.

The main idea of proof is to reduce,
via the pseudo-conformal invariance,
the proof of the uniqueness of multi-solitons to that of
multi-bubble blow-up solutions.
The main effort is hence dedicated to the latter issue,
to which, inspired by \cite{RS11},
the sharp singularity analysis is applied.
We show that the convergence indeed can be upgraded to
the much faster exponential decay rate,
which in particular is beyond the third order  $(T-t)^{3+}$ in \eqref{v-S-H1-3+}
and thus yields the desirable uniqueness of multi-bubble blow-up solutions.

The crucial ingredients in the upgradation procedure
are the monotonicity properties of different functionals
adapted to the multi-bubble case.
More delicately,
in each upgradation step
the monotonicity of functionals relies on suitable estimates of the remainder
and geometrical parameters in the previous step.
We show that the initial low convergent rate $(T-t)^{0+}$ in \eqref{v-S-H1-0+}
is effective to run the whole upgradation procedure.

One major difficulty, particularly in the multi-bubble case,
is that
in all the controls of functionals arises the localized mass
which may restrict the upgradation strength.
Such a problem is absence in the single bubble case,
because the localized mass vanishes due to the conservation law of mass.
One keypoint is that
two more orders (i.e., $(T-t)^{2+}$)
can be explored
for the localized mass by the local virial identities in \cite{M93}.
This convergence rate is effective to serve as the basis of estimates of localized mass.
More importantly,
we relate the localized mass and the remainder together
and upgrade their estimates simultaneously
by iteration arguments through several Gronwall type inequalities.

Let us also mention that,
another challenge in the low convergence regime is to identify
the exact value of the energy of solutions.
The key fact here is that
the remainder exhibits dispersion in the energy space along a sequence,
which enables us to obtain the energy quantization,
that is, the energy of multi-bubble blow-up solutions
admits the quantization into the sum of the energies of pseudo-conformal blow-up solutions.
This is the key towards the derivation of refined energy estimate
in the upgradation procedure. \\

{\bf Notations.}
We use the standard Sobolev spaces
$H^{s,p}(\bbr^d)$, $s\in \bbr, 1\leq p\leq\9$.
In particular,
$L^p := H^{0,p}(\bbr^d)$ is
the space of $p$-integrable (complex-valued) functions,
$L^2$ denotes the Hilbert space endowed with the scalar product
$\<v,w\> =\int_{\bbr^d} v(x) \ol w(x)dx$,
and $H^s:= H^{s,2}$.
As usual,
if $B$ is a Banach space,
$L^q(0,T;B)$ means the space of all integrable $B$-valued functions $f:(0,T)\to B$ with the norm
$\|\cdot\|_{L^q(0,T;B)}$,
and $C([0,T];B)$ denotes the space of all $B$-valued continuous functions on $[0,T]$ with the sup norm over $t$.

We also use the notation $\dot{g} = \frac{d}{dt}g$ for any $C^1$
function $g$ defined on an open time interval.
As $t\to T$ or $t\to \9$,
$f(t)=\calo(g(t))$
means that $|f(t)/g(t)|$ stays bounded,
and $f(t)=o(g(t))$ means that $|f(t)/g(t)|$ converges to zero.

Throughout this paper,
the positive constants $C$ and $\delta$ may change from line to line.

\subsection{Formulation of main results} \label{Subsec-Main}

Throughout this paper
we mainly consider two cases below:

{\rm \bf Case (I).}
$\{x_k\}_{k=1}^K$ are arbitrarily distinct points in $\bbr^d$,
and $\{\omega_k\}_{k=1}^K (\subseteq \bbr^+)$ satisfy that
for some $\omega>0$, $|\omega_k - \omega| \leq \ve$ for every $1\leq k \leq K$,
where $\ve >0$;

{\rm \bf Case (II).}
$\{\omega_k\}_{k=1}^K$ are arbitrary points in $\bbr^+$,
and $\{x_k\}_{k=1}^K (\subseteq \bbr^d)$ satisfy that
$|x_j-x_k| \geq \ve^{-1}$ for every $1\leq j \neq k \leq K$,
where $\ve>0$.

One interesting application is the case where $\{x_k\}$
are arbitrarily distinct points in $\bbr^d$ and
the frequencies $\{\omega_k\}$ are the same.
In the single bubble case where $K=1$,
the blow-up point $x$ and the frequency $\omega$ can be
arbitrary points in $\bbr^d$ and $\bbr^+$, respectively.

Let us first consider the multi-bubble blow-up solutions.
As mentioned above,
the multi-bubble blow-up solutions with pseudo-conformal blow-up rate
were first constructed
in the pioneering work of Merle \cite{M90}.
Theorem \ref{Thm-Uniq-Blowup-3+} below gives the uniqueness class of multi-bubble blow-up solutions with
convergence rate $(T-t)^{3+}$.

\begin{theorem} \label{Thm-Uniq-Blowup-3+}
(\cite[Theorems 2.7 and 2.15]{SZ20})
Consider equation \eqref{equa-NLS} in $\bbr^d$ for $d=1,2$.
Let $T\in \bbr$, $K \in \mathbb{N}\setminus\{0\}$.
Let $\{\vartheta_k\}_{j=k}^K \subseteq \bbr$,
$\{x_k\}_{k=1}^K$
and $\{\omega_k\}_{k=1}^K$
satisfy  either {\rm Case (I)} or {\rm Case (II)}.
Then, for any $\zeta\in (0,1)$,
there exists $\ve^*>0$, such that for any $0<\ve<\ve^*$,
there exists a unique multi-bubble blow-up solution $v$ to \eqref{equa-NLS}
such that
\begin{align} \label{v-S-H1-3+-Thm}
\|v(t)-\sum_{k=1}^KS_k(t)\|_{H^1} = \calo((T-t)^{3+\zeta}),\ \ for\ t\ close\ to\ T,
\end{align}
where $S_k$, $1\leq k\leq K$, are the pseudo-conformal blow-up solutions given by \eqref{Sj-blowup}.
Moreover, the unique multi-bubble blow-up solution $v$ satisfies that
for some  $\delta>0$.
\begin{align} \label{v-S-sig-exp+-Thm}
\|v(t)-\sum_{k=1}^KS_k(t)\|_{\Sigma} = \calo(e^{-\frac{\delta}{T-t}}),\ \ for\ t\ close\ to\ T.
\end{align}
\end{theorem}

\begin{remark} \label{Rem-uniq-blow-H13+}
$(i)$ In the recent work \cite{SZ20},
the construction and conditional uniqueness of
multi-bubble blow-up solutions are also obtained for stochastic nonlinear Schr\"odinger equations
(SNLS)
driven by Wiener processes,
particularly in the absence of the pseudo-conformal invariance and the
conservation law of energy.
For the interested readers,
we would like to refer to
\cite{FSZ20,SZ19} for the construction of
stochastic blow-up solutions with loglog blow-up rate or with critical mass,
and \cite{BD02,BD05} for the noise effect on blow-up.
For other related topics of SNLS,
see \cite{BRZ16.1,BRZ18,BM14,FX19,HRZ18,Z19}
and the references therein.

$(ii)$ For equation \eqref{equa-NLS},
the smallness of $T$ assumed in \cite{SZ20} can be removed.
It is assumed there simply because the Wiener processes start moving at time zero.
One may replace the smallness of $T$ by taking $t$ close enough to $T$ in the proof
of \cite[Theorem 2.15]{SZ20}.

$(iii)$
By virtue of the pseudo-conformal invariance,
Theorem \ref{Thm-Uniq-Blowup-3+} also yields the uniqueness of multi-solitons $u$ to \eqref{equa-NLS}
satisfying either
\begin{align} \label{u-W-H1-6+}
    \| u(t) - \sum\limits_{k=1}^K W_k(t)\|_{H^1} = \calo\( \frac{1}{t^{5+}}\),\ \ for\ t\ large\ enough,
\end{align}
or
\begin{align} \label{u-W-Sigma-4+}
    \| u(t) - \sum\limits_{k=1}^K W_k(t)\|_{\Sigma} = \calo\( \frac{1}{t^{4+}}\),\ \ for\ t\ large\ enough.
\end{align}

Actually, by Lemma \ref{Lem-uW-H1-Sigma} and the inequalities \eqref{v-u-L2-pct} and \eqref{v-u-H1-pct} below,
the condition \eqref{u-W-H1-6+} or  \eqref{u-W-Sigma-4+}
suffices to verify \eqref{v-S-H1-3+-Thm}
for $v=P_T(u)$ given by \eqref{pseu-conf-transf}.
Hence, Theorem \ref{Thm-Uniq-Blowup-3+} yields the uniqueness of multi-bubble blow-up solutions,
and thus of the corresponding multi-solitons.
\end{remark}

The first main result of this paper is concerned with the uniqueness class of multi-bubble blow-up solutions
to \eqref{equa-NLS}.
The precise statements are formulated below.

\begin{theorem}\label{Thm-Uniq-Blowup}
Consider equation \eqref{equa-NLS} in dimensions $d=1,2$.
Let $T\in \bbr$, $K \in \mathbb{N}\setminus\{0\}$.
Let $\{\vartheta_k\} \subseteq \bbr$,
$\{\omega_k\}$ and $\{x_k\}$ satisfy either {\rm Case (I)} or {\rm Case (II)}.
Then, for any $\zeta\in (0,1)$,
there exists $\ve^*>0$ such that the following holds.
For any $0<\ve<\ve^*$, there exists a unique multi-bubble blow-up solution $v$ to \eqref{equa-NLS}
such that
\begin{align} \label{v-S-H1-o1-Thm}
    \| v(t) - \sum\limits_{k=1}^K S_k(t) \|_{L^2}
    + (T-t) \| \na v(t) - \sum\limits_{k=1}^K \na S_k(t) \|_{L^2} = o(1),\ \ as\ t\ close\ to\ T,
\end{align}
and additionally
\begin{align} \label{nav-naS-iint-0+-Thm}
     \frac{1}{T-t} \int_t^T \frac{1}{T-s} \int_s^T \|v (r) - \sum\limits_{k=1}^K S_k(r)  \|_{H^1}^2 dr ds
     = \calo((T-t)^{\zeta}),
\end{align}
where $S_k$, $1\leq k\leq K$, are the
pseudo-conformal blow-up solutions
given by \eqref{Sj-blowup}.
Moreover, the unique solution $v$ converges exponentially fast to
$\sum_{k=1}^K S_k$   in
the pseudo-conformal space, i.e.,
there exists $\delta>0$ such that
\begin{align} \label{v-S-Sigma-exp-Thm}
   \|v(t) - \sum\limits_{k=1}^K S_k(t) \|_{\Sigma} = \calo(e^{-\frac{\delta}{T-t}}), \ \ for\ t\ close\ to\ T.
\end{align}

In particular,
the above results hold for
the multi-bubble blow-up solutions $v$ to \eqref{equa-NLS} such that
\begin{align} \label{v-S-H1-0+-Thm}
   \| v(t) - \sum\limits_{k=1}^K S_k(t) \|_{H^1}  = \calo( (T-t)^{\zeta}),\ \ for\ t\ close\ to\ T.
\end{align}
\end{theorem}

\begin{remark}
$(i)$ The conditions \eqref{v-S-H1-o1-Thm} and \eqref{nav-naS-iint-0+-Thm} (or, \eqref{v-S-H1-0+-Thm})
are much weaker than the previous one \eqref{v-S-H1-3+-Thm}.
Actually, the proof of Theorem 2.15 in \cite{SZ20}
is based on the control of modified generalized energy,
and the convergence rate $(T-t)^{3+}$
is enough to provide a Lyapunov type rigidity.
However,
the proof of Theorem \ref{Thm-Uniq-Blowup} is much more delicate in the low convergence regime.
It requires several upgradation procedures and relies crucially on the monotonicity of
several new functionals, including the (modified) localized virial functionals.

$(ii)$
It is conjectured by Merle and Rapha\"el \cite{MR05} that,
the general blow-up solutions to \eqref{equa-NLS} can be decomposed into
multi blow-up bubbles and a regular residue.
Theorem \ref{Thm-Uniq-Blowup} gives one uniqueness class of such blow-up solutions
in the case where the mass is $K\|Q\|_{L^2}^2$
and the asymptotic behaviors
\eqref{v-S-H1-o1-Thm} and \eqref{nav-naS-iint-0+-Thm} (or, \eqref{v-S-H1-0+-Thm})
are satisfied.
\end{remark}

Our next goal is to enlarge the uniqueness class of multi-solitons to \eqref{equa-NLS}.
The precise statement for the uniqueness in energy class is formulated below.

\begin{theorem}\label{Thm-Uniq-Solitons-H1}
Consider equation \eqref{equa-NLS} in dimensions $d=1,2$.
Let  $K \in \mathbb{N}\setminus\{0\}$.
Let $\{\vartheta_k\}\subseteq \bbr$, $\{\omega_k\}$ and $\{v_k\}$
satisfy either {\rm Case (I)} or {\rm Case (II)}
with $v_k$ replacing $x_k$, $1\leq k\leq K$.
Then, for any $\zeta\in (0,1)$,
there exists $\ve^*>0$ such that the following holds.
For any $0<\ve<\ve^*$,  there exists a unique multi-soliton $u$ to \eqref{equa-NLS}
such that
\begin{align} \label{u-W-H1-2+-Thm}
\|u(t)- \sum_{k=1}^K W_k(t)\|_{H^1} = \calo\(\frac{1}{t^{2+\zeta}}\), \ \ for\ t\ large\ enough,
\end{align}
where $\{W_k\}$ are the solitons given by \eqref{Wj-soliton}.
Moreover,
the unique multi-soliton $u$ converges exponentially fast to $\sum_{k=1}^K W_k$
in the more regular pseudo-conformal space,
i.e., for some $\delta>0$,
\begin{align} \label{u-W-Sigma-exp-Thm}
\|u(t)- \sum_{k=1}^K W_k(t)\|_{\Sigma} = \calo(e^{- \delta t}),  \ \ for\ t\ large\ enough.
\end{align}
\end{theorem}

\begin{remark}
$(i)$ The uniqueness of multi-solitons to NLS
was first obtained by C\^ote and Friederich \cite{CF20}
in the $L^2$-subcritical and critical cases,
provided that the convergence rate is $(1/t)^{N}$ for $N$ large enough.
Theorem \ref{Thm-Uniq-Solitons-H1} shows that
the uniqueness class of multi-solitons to $L^2$-critical \eqref{equa-NLS}
can be enlarged in the low convergence regime with rate $(1/t)^{2+}$.
Differently from \cite{CF20},
the proof of Theorem \ref{Thm-Uniq-Solitons-H1} is based on
the pseudo-conformal invariance and
the uniqueness of multi-bubble blow-up solutions
in Theorem \ref{Thm-Uniq-Blowup} .

$(ii)$
Thanks to the smoothness result in \cite{CF20},
the unique multi-solitons in Theorem \ref{Thm-Uniq-Solitons-H1}
belong to the space $C([T, \9); H^\9)$ for some $T\in \bbr$
and satisfy that for all integer $s\geq 0$,
\begin{align}
\|u(t)- \sum_{k=1}^K W_k(t)\|_{H^s}  \leq Ce^{- \delta t},  \ \ for\ t\ large\ enough,
\end{align}
where $C,\delta>0$.
\end{remark}

In the  pseudo-conformal space,
the uniqueness class of multi-solitons to \eqref{equa-NLS}
can be further enlarged with even lower convergence rate $(1/t)^{\frac 12+}$.
This is the content of Theorem \ref{Thm-Uniq-Solitons-Sigma} below.

\begin{theorem} \label{Thm-Uniq-Solitons-Sigma}
Assume the conditions in Theorem \ref{Thm-Uniq-Solitons-H1} to hold.
Assume additionally that $v_k \not=0$, $1\leq k\leq K$.
Then, for any $\zeta\in (0,1)$,
there exists $\ve^*>0$ such that for any $0<\ve<\ve^*$,
there exists a unique multi-soliton $u$ to \eqref{equa-NLS}
such that
\begin{align} \label{u-W-Sigma-12+-Thm}
\|u(t)- \sum_{k=1}^K W_k(t)\|_{\Sigma} = \calo\(\frac{1}{t^{\frac 12+\zeta}}\), \ \ for\ t\ large\ enough.
\end{align}
Moreover, the unique multi-soliton $u$ converges exponentially fast to $\sum_{k=1}^K W_k$ in $\Sigma$.
\end{theorem}

\begin{remark}
By the pseudo-conformal invariance,
it is easy to see from  \eqref{v-S-H1-0+-Thm} and
the inequalities \eqref{v-u-L2-pct}-\eqref{v-u-H1-pct} below that,
the uniqueness holds for the multi-solitons $u$ to \eqref{equa-NLS} such that
\begin{align}
    \|u(t) - \sum\limits_{k=1}^K W_k(t) \|_{\Sigma} = \calo\(\frac{1}{t^{1+}}\), \ \ for\ t\ large\ enough.
\end{align}
However,
the regime of low convergence rates
from $\frac 12 +$ to $1$ requires more delicate proof.
The key ingredients are the
monotonicity of virial functional
and the refined space time estimate of the gradient of multi-solitons,
which allow to verify the weaker double average condition \eqref{nav-naS-iint-0+-Thm}
for the corresponding multi-bubble blow-up solutions,
rather than the pointwise decay condition \eqref{v-S-H1-0+-Thm}.
\end{remark}

We conclude this subsection with an interesting application to the single soliton case $K=1$,
which is
related to the solitary wave conjecture.

\begin{corollary}\label{Cor-Uniq-Soliton-Single}
Consider equation \eqref{equa-NLS} in dimensions $d=1,2$. For any $\zeta\in(0,1)$,
assume that $u$ is a solution to \eqref{equa-NLS}
satisfying
\begin{align} \label{u-W-H1-2+-Cor}
\|u(t)-  \omega^{-\frac d2}Q \(\frac{\cdot-v t}{\omega} \)e^{i(\half v\cdot x-\frac{1}{4}|v|^2t+\omega^{-2}t+\vartheta)}\|_{H^1} = \calo\(\frac{1}{t^{2+\zeta}}\), \ \ for\ t\ large\ enough,
\end{align}
where $v, \vartheta\in \bbr, \omega>0$.
Then, $u$ coincides with the solitary wave solution, i.e.,
\begin{align} \label{u-W-H1-2+-single-Thm}
   u(t,x)=\omega^{-\frac d2}Q \(\frac{x-v t}{\omega} \)e^{i(\half v\cdot x-\frac{1}{4}|v|^2t+\omega^{-2}t+\vartheta)}.
\end{align}
Moreover, if $v=0$,
then the same uniqueness holds for solutions  with lower convergence rate
\begin{align} \label{u-W-H1-1+-single-Thm}
\|u(t)-  \omega^{-\frac d2}Q \(\frac{\cdot}{\omega} \)e^{i(\omega^{-2}t+\vartheta)}\|_{H^1} = \calo\(\frac{1}{t^{1+\zeta}}\), \ \ for\ t\ large\ enough.
\end{align}
\end{corollary}

\begin{remark}
The solitary wave conjecture is affirmative in the
pseudo-conformal space by the rigidity result in \cite{M93}
and the pseudo-conformal invariance.
In the Sobolev space,
when $d=2,3$ it is proved in \cite{LZ12}  for $H^1$ radial
solutions,
and when $d\geq 4$ it is proved in \cite{KLVZ09} and \cite{LZ09}, respectively,
for $H^1$ and $H^s$ ($s>0$) radial solutions.
Corollary \ref{Cor-Uniq-Soliton-Single} shows that,
in dimensions $d=1,2$,
this conjecture is affirmative
for general $H^1$ solutions
satisfying either the asymptotic behavior \eqref{u-W-H1-2+-Cor} or \eqref{u-W-H1-1+-single-Thm}
when $v=0$.
\end{remark}

\subsection{Strategy of the proof} \label{Subsec-strategy-proof}

We make use of the pseudo-conformal invariance
to reduce the uniqueness  of multi-solitons to that of multi-bubble blow-up solutions,
and then perform several upgradation procedures to upgrade the
convergence
to the exponential decay rate.
Let us give a sketch of the proof below.

{\bf Step 1.} Geometrical decomposition and preliminary controls.

We first obtain the geometrical decomposition of multi-bubble blow-up solution $v$ to \eqref{equa-NLS}, i.e.,
\begin{align}  \label{v-dec-intro}
    v (t,x) =\sum\limits_{k=1}^K \lbb_k^{-\frac d2}(t) Q_k\(t,\frac{x-\a_k(t)}{\lbb_k(t)}\) e^{i\theta_k(t)} + R(t,x)
    \ \(:=\sum_{k=1}^{K} U_k(t,x) +R(t,x)\),
\end{align}
where
$Q_k(t,y) = Q(y) e^{i(\beta_k(t)\cdot y - \frac 14 \g_k(t)|y|^2)}$
and the remainder $R$ satisfies the orthogonality conditions in \eqref{ortho-cond-Rn-wn} below.
The geometrical decomposition enables us to perform the robust modulation method
and to reduce the analysis of the blow-up dynamics of $v$ to that of
finite dimensional geometrical parameters $\calp=(\lbb, \a, \beta, \g, \theta)$
and the remainder $R$.

It should be mentioned that,
the geometrical decomposition \eqref{v-dec-intro}
should be valid for any time close to the blow-up time,
which, however, may give rise to the singularity of Jacobian matrix
with respect to $\calp$.
This is quite different from the construction of blow-up solutions
(see, e.g., \cite{SZ19,SZ20}),
where the geometrical decomposition is performed on the intervals away from the blow-up time.

To fix this problem,
inspired by \cite{MR05.2},
the idea is to introduce a new family of parameters
$\wt \calp = (\wt \lbb, \wt \a, \wt \beta, \wt \g, \wt \theta)$ to measure the small deformations
of pseudo-conformal blow-up solutions, implied by \eqref{v-S-H1-o1-Thm},
and to obtain a universal bound of the corresponding determinant of Jacobian matrix with respect to these new parameters.

As a byproduct,
the preliminary controls of remainder and geometrical parameters are obtained as well,
which serve as the basic estimates at the beginning of the upgradation procedure.

{\bf Step 2.} Controls of functionals adapted to the multi-bubble case.

The first upgradation step relies on the coercivity control of energy
around the blow-up profile:
\begin{align}  \label{energy-esti-intro}
     & \sum_{k=1}^{K}\frac{|\beta_{k}|^2 }{2\lbb_{k}^2} \|Q\|_{2}^2
        + \sum_{k=1}^{K}  \frac{\g_k^2}{8\lbb_k^2} \|yQ\|_{L^2}^2
        + C_1 \frac{D^2}{(T-t)^2}  \nonumber \\
\leq& E(v)
     + \sum_{k=1}^K \frac{|\beta_k|^2}{2\lbb_k^2} M_k
      + \sum_{k=1}^{K}\frac{1}{2\lbb_{k}^2} M_k
      + C_2  \(\sum\limits_{k=1}^K \frac{M_k^2}{(T-t)^2} + e^{-\frac{\delta}{T-t}}\) .
\end{align}
Here $C_1, C_2, \delta>0$,
$D(t):= \|R(t)\|_{L^2} + (T-t)\|\nabla R(t)\|_{L^2}$,
and $M_k$ denotes the localized mass
\begin{align*}
   M_k:= 2 {\rm Re} \< R_k, U_k\> + \int |R|^2 \Phi_k dx, \ \ with\ R_k = R\Phi_k,\ 1\leq k\leq K,
\end{align*}
where $\{\Phi_k\}$ are the localization functions (see \eqref{phi-local} below).

Note that,
the quantity $D$ measures the smallness of remainder and indeed
plays the key role in the proof of the uniqueness of multi-bubble blow-up solutions.
As a matter of fact,
the upgradation procedure is mainly dedicated to upgrading the convergence rate of $D$.

In view  of the conservation law of energy,
estimate \eqref{energy-esti-intro} suffices to upgrade the convergence rate to the first order,
i.e., $D(t)=\calo(T-t)$.

Two major problems, however, arise in the above energy estimate \eqref{energy-esti-intro}.
The first one is, that the exact value of energy $E(v)$ is a priori unclear due to the low convergence rate in \eqref{v-S-H1-o1-Thm}.

In order to identify the energy of solutions,
inspired by \cite{RS11},
we introduce the localized virial functional adapted to the multi-bubble case
\begin{align} \label{local-virial-def-intro}
   \mathcal{L}:=\sum_{k=1}^K \half {\rm Im}
            \int (\nabla\chi_A) \(\frac{x-\alpha_{k}}{\lambda_k}\)\cdot\nabla R \ol{R} \Phi_kdx
            - \frac{\g_k}{4\lbb_k}\|xQ\|_{L^2}^2,
\end{align}
where $\chi_A$ is a suitable localized function (see Subsection \ref{Subsec-Modif-Loc-Vir} below).

Let us mention that,
localized virial functionals were first introduced by Martel and Merle \cite{MM01}
to study the $L^2$-critical gKdV equation,
and later used by Rapha\"{e}l and  Szeftel \cite{RS11} in the inhomogeneous NLS setting.
The corresponding monotonicity property
in particular yields a space time estimate of remainder around the singularity.

In the multi-bubble case,
the localization functions $\{\Phi_K\}$ in \eqref{local-virial-def-intro}  are introduced
 in this particular way
mainly to ensure the monotonicity property (see Theorem \ref{Thm-Loc-virial} below).
Actually,
one difficulty in the multi-bubble case is to
beat the interactions between different localized remainders $R_j$ and $R_k$, $j\not=k$,
which are usually not easy to control,
due to the few knowledge of remainder.
The keypoint here is,
that cancellations up to high orders can be gained for these coupling terms,
by suitably choosing the localization functions and test functions.
See, e.g., the proof of $H^1$ dispersion in  Lemma \ref{Lem-R-H1-In},
see also the proof of the monotonicity of localized virial functional in Theorem \ref{Thm-Loc-virial}
and the related arguments of \cite[Lemma 5.12]{SZ20}.

The refined space time estimate of remainder yields
the precise asymptotic order of parameter $\g_k$
which, however, can not be obtained from the previous geometrical decomposition.
More importantly,
it reveals the key $H^1$ dispersion of remainder
along a sequence $\{t_n\}$ to $T$,
which  enables us to obtain the crucial energy quantization in Theorem \ref{Thm-Energy-Ident}.

As a consequence,  the refined energy estimate is derived: for some $C_1, C_2>0$,
\begin{align}  \label{energy-esti-refined-intro}
      \sum_{k=1}^{K}\frac{|\beta_{k}|^2 }{2\lbb_{k}^2} \|Q\|_{2}^2
        + C_1 \frac{D^2}{(T-t)^2}
\leq& \frac{ \|yQ\|_{L^2}^2}{8} \sum\limits_{k=1}^K\(\omega_k^2 - \frac{\g_k^2}{\lbb_k^2}\)
     + \sum_{k=1}^K \frac{|\beta_k|^2}{2\lbb_k^2} M_k
    + \sum_{k=1}^{K}\frac{1}{2\lbb_{k}^2}M_k   \nonumber \\
    &  + C_2 \(\sum\limits_{k=1}^K \frac{M_k^2}{(T-t)^2} + e^{-\frac{\delta}{T-t}}\).
\end{align}

The second problem of energy estimate
is that,
at this stage, the first term on the R.H.S. above is merely of order $\calo(1)$,
which indeed ceases the further upgradation.

This leads us to introduce another new modified localized virial functional
adapted to the multi-bubble case
\begin{align} \label{wtI-def}
     \mathscr{L}
           :=\sum_{k=1}^K  \frac{\g_k}{2\lbb_k} {\rm Im}
            \int (\nabla\chi_A) \(\frac{x-\alpha_{k}}{\lambda_k}\)\cdot\nabla R\ol{R}\Phi_kdx
            - \frac{\g_k^2}{8\lbb_k^2} \|xQ\|_2^2,
\end{align}
and to derive the corresponding monotonicity property,
which in particular yields the coercivity type control
\begin{align}  \label{Rj-lbbj3-refine-intro}
  \tilde c \sum_{k=1}^{K}\int_{t}^{T} \frac{\|R_{k}\|_{L^2}^2}{\lbb_{k}^3} ds
   \leq& \frac{\|yQ\|_{L^2}^2}{8} \sum_{k=1}^{K}\(\frac{\g_k^2}{\lbb_k^2}-\omega^2_k\)
         +C\|R\|_{L^2}\|\nabla R\|_{L^2}
        + \bigg|\int_{t}^{T}\sum_{k=1}^{K}\frac{\g_k}{\lbb_k^4}M_kds\bigg|
          + C \int_{t}^{T} Er ds,
\end{align}
where $ \tilde c, C>0$ and $Er$ is the error term given by \eqref{Error} below.
The key fact is, that
the trouble $\calo(1)$ terms in both estimates \eqref{energy-esti-refined-intro} and \eqref{Rj-lbbj3-refine-intro}
cancel each other out
and the higher order terms remain.
This makes it possible to further upgrade the convergence rate of remainder.

Thus, combining estimates \eqref{energy-esti-refined-intro} and \eqref{Rj-lbbj3-refine-intro} altogether
we lead to the refined estimate:
\begin{align} \label{D-iter-2-intro}
&\sum_{k=1}^{K}\frac{|\beta_{k}(t)|^2 }{\lbb_{k}^2} \|Q\|_{2}^2
     + \frac{D^2(t)}{(T-t)^2}
     + \sum_{k=1}^{K}\int_{t}^{T} \frac{\|R_{k}\|_{L^2}^2}{\lbb_{k}^3}  ds   \nonumber \\
\leq& C\( (T-t)^2
      +\bigg|\sum_{k=1}^{K}\frac{M_k}{\lbb_{k}^2}\bigg|
      +\bigg|\sum\limits_{k=1}^K \int_{t}^{T}\frac{\g_k}{\lbb_k^4}{M_k}ds\bigg|
       +\int_{t}^{T}   \frac{D^2 + \sum_{k=1}^K |M_k|}{(T-t)^2} ds \).
\end{align}

Note that,
the localized mass $M_k$ arises
in all the controls of functionals above.
The main estimates of localized mass are collected in Theorem \ref{Thm-Loc-Mass} below.

{\bf Step 3.} Upgradation to the higher convergence rate.

Quite differently from the first upgradation step,
the upgradation to the higher convergence rate requires
delicate iteration arguments through different Gronwall type inequalities.

More precisely,
two Gronwall type inequalities will be derived from the refined estimate \eqref{D-iter-2-intro},
which enable to upgrade the convergence rate
to the second order, i.e.,
$D(t) = \calo((T-t)^{2})$.

For the further upgradation,
we use the modified generalized energy introduced in \cite{SZ20}
\begin{align} \label{gen-energy-def}
 \mathscr{I}:= &\frac{1}{2}\int |\nabla R|^2dx+\frac{1}{2}\sum_{k=1}^K\int\frac{1}{\lambda_{k}^2} |R|^2 \Phi_kdx
           -{\rm Re}\int F(u)-F(U)-f(U)\ol{R}dx \nonumber \\
&+\sum_{k=1}^K\frac{\gamma_{k}}{2\lambda_{k}}{\rm Im} \int (\nabla\chi_A) \(\frac{x-\alpha_{k}}{\lambda_{k}}\)\cdot\nabla R \ol{R}\Phi_kdx,
\end{align}
where $F(z):= \frac{d}{2d+4} |z|^{2+ \frac 4d}$,
$f(z):= |z|^{\frac 4d}z$, $z\in \mathbb{C}$.
Again, the corresponding monotonicity property
allows to derive a new Gronwall type inequality,
which enables us to upgrade the convergence rate to the much faster exponential decay rate,
that is, $D(t) = \calo(e^{-\frac{\delta}{T-t}})$
for some $\delta>0$.
The exponential decay rate of geometrical parameters and modulation equations is then obtained as well.

{\bf Step 4.} Proof of the main results.

By virtue of the exponential decay results in the previous step,
the multi-bubble blow-up solutions indeed converge exponentially fast
to the sum of pseudo-conformal blow-up solutions,
which, via Theorem \ref{Thm-Uniq-Blowup-3+}, yields the uniqueness in Theorem \ref{Thm-Uniq-Blowup}.

The proof of Theorems \ref{Thm-Uniq-Solitons-H1} and \ref{Thm-Uniq-Solitons-Sigma}
is based on Theorem \ref{Thm-Uniq-Blowup} and the pseudo-conformal invariance.
In particular,
the key observation in the proof of Theorem \ref{Thm-Uniq-Solitons-Sigma} is that,
an improved convergence rate can be gained for the
space time estimate of gradient
by the monotonicity of virial functional
(see Lemma \ref{Lem-Mono-virial-soliton} and Corollary \ref{Cor-v-H1-Sig-2nu} below).
This is quite different from the proof of Theorem \ref{Thm-Uniq-Solitons-H1},
which relies on the pointwise decay rate in $\Sigma$.
The improved space time estimate is the key towards the verification
of the weaker double average condition \eqref{nav-naS-iint-0+-Thm},
which enables to enlarge the
uniqueness class of multi-solitons
with even lower convergence rate $(1/t)^{\frac 12+}$.

The remainder of this paper is structured as follows.
Section \ref{Sec-Geom-Decomp} contains the geometrical decomposition
and preliminary controls of the remainder and modulations equations.
Then,
Sections \ref{Sec-Energy-Locvirial} is dedicated to the controls of
localized mass and energy and the key monotonicity properties of (modified) localized virial functionals.
The crucial upgradation procedures
are then performed in Section \ref{Sec-Upgrad-Remainder}.
Eventually, the main results are proved in  Section \ref{Sec-Proof-Main}.
For simplicity, some tools and the proof of modulation equations are presented in the Appendix, i.e., Section \ref{Sec-App}.

\section{Geometrical decomposition and modulation equations} \label{Sec-Geom-Decomp}

The main result of this section is the geometrical decomposition of
multi-bubble blow-up solutions for all time close to the blow-up time $T$.

Let us use the notations $\lbb=(\lbb_k)_{1\leq k\leq K}, \g=(\g_k)_{1\leq k\leq K}, \theta=(\theta_k)_{1\leq k\leq K}$
for the vectors in $\bbr^K$,
and $\a=(\a_k)_{1\leq k\leq K}$, $\beta=(\beta_k)_{1\leq k\leq K}$
for the vectors in $\bbr^{dk}$, where $\a_k=(\a_{k,i})_{1\leq i\leq K}, \beta_k=(\beta_{k,i})_{1\leq i\leq d} \in \bbr^d$, $1\leq k\leq K$.
Hence, we have
$\calp=(\lbb, \a, \beta, \g, \theta)
\in \mathbb{X}:=\bbr^K \times \bbr^{dK} \times \bbr^{dK} \times \bbr^K \times \bbr^K$.

Let $|\lbb|:= \sum_{k=1}^K |\lbb_k|$.
Similar notations also apply to the vectors $\g, \theta \in \bbr^K$,
$\a, \beta \in \bbr^{dK}$ and $\calp\in \mathbb{X}$.

\subsection{Geometrical decomposition}   \label{Subsec-geom-dec}

\begin{theorem} (Geometrical decomposition)   \label{Thm-geometri-dec}
Let $v$ be the multi-bubble blow-up solution to \eqref{equa-NLS} satisfying \eqref{v-S-H1-o1-Thm}.
Then, for $T^*$ sufficiently close to $T$,
there exist  unique modulation parameters
$\mathcal{P} = (\lbb, \a, \beta, \g, \theta)
\in C^1([T^*,T); \mathbb{X})$,
such that
$v$ admits the geometrical decomposition
\begin{align} \label{v-dec}
    v(t,x)=\sum_{k=1}^{K} U_k(t,x)+R(t,x)\ (:= U(t,x)+R(t,x)),\ \ t\in [T^*, T),\ x\in \bbr^d,
\end{align}
where
\begin{align} \label{Uj-Qj-Q}
   U_k(t,x) = \lbb_k^{-\frac d2}(t) Q_k\(t,\frac{x-\a_k(t)}{\lbb_k(t)}\) e^{i\theta_k(t)} \ \ with\
   Q_k(t,y) = Q(y) e^{i(\beta_k(t)\cdot y - \frac 14 \g_k(t)|y|^2)},
\end{align}
and  the following orthogonality conditions hold on $[T_*,T)$:
for each $1\leq k\leq K$,
\be\ba\label{ortho-cond-Rn-wn}
&{\rm Re}\int (x-\a_{k}) U_{k}(t)\ol{R}(t)dx=0,\ \
{\rm Re} \int |x-\a_{k}|^2 U_{k}(t) \ol{R}(t)dx=0,\\
&{\rm Im}\int \nabla U_{k}(t) \ol{R}(t)dx=0,\ \
{\rm Im}\int \Lambda_k U_{k}(t) \ol{R}(t)dx=0,\ \
{\rm Im}\int \varrho_{k}(t) \ol{R}(t)dx=0.
\ea\ee
Here,
\be
\Lambda_k:=\frac{d}{2}I_d+(x-\a_k)\cdot \nabla,
\ee
\begin{align} \label{rhon}
   \varrho_{k}(t,x)= \lbb^{-\frac d2}_{k} \rho_{k}\(t,\frac{x-\a_{k}}{\lbb_n}\) e^{i\theta_{k}}
 \ \ with\  \  \rho_{k}(t,y) := \rho(y)^{i(\beta_{k}(t)\cdot y - \frac 14 \g_{k}(t) |y|^2)},
\end{align}
and $\rho$ is given by \eqref{rho-def} below.

Moreover, the following  estimates hold for the geometrical parameters and remainder:
\begin{align}
   &\sum_{k=1}^{K}(|\lbb_k(t)-\omega_k(T-t)|+|\alpha_k(t)-x_k|)= o(T-t), \label{lbb-a-t-0} \\
   &\sum_{k=1}^{K}(|\b_k(t)|+|\g_k(t)-\omega_k^2(T-t)|+|\t_k(t)-\omega_k^{-2}(T-t)^{-1} - \vartheta_k|)= o(1), \label{P-t-0} \\
   &  \|R(t)\|_{L^2}+(T-t)\|\nabla R(t)\|_{L^2}= o(1). \label{R-t-0}
\end{align}
\end{theorem}

\begin{remark}
$(i)$  The geometrical decomposition  \eqref{v-dec} is valid for any time close to the blow-up time,
and thus may give rise to the singularity of Jacobian matrix used in \cite{SZ20}.

$(ii)$ Estimate \eqref{lbb-a-t-0} gives the precise leading asymptotic order of $\lbb_k$,
i.e.,
\begin{align} \label{lbbk-t-approx}
     \lbb_k(t) \thickapprox \omega_k (T-t),\ \ for\ t\ close\ to\ T,
\end{align}
which particularly characterizes the blow-up rate.
However, estimate \eqref{P-t-0} is insufficient to yield
\begin{align} \label{gk-t2}
     \g_k(t) \thickapprox \omega_k^2(T-t),\ \ for\ t\ close\ to\ T.
\end{align}
As a matter of fact,
such precise estimate \eqref{gk-t2} will be derived after analyzing the localized virial functional
(see Theorem \ref{Thm-gamj-wj2-t} below).
It should be also mentioned that,
\eqref{gk-t2} is important in the derivation of refined energy estimate,
which allows to upgrade the convergence rate of remainder to the second order
(see Subsections \ref{Subsec-Energy-Quant} and \ref{Subsec-D-2} below).
\end{remark}

In order to avoid the singularity of Jacobian matrix,
inspired by \cite{MR05.2},
the idea here is
to take, via  \eqref{v-S-H1-o1-Thm},
the solution $v$ as a small perturbation of $\sum_{k=1}^K S_k$.
Hence, a new family of parameters
$\wt \calp_k=(\wt \lbb_k, \wt \a_k, \wt \beta_k, \wt \g_k, \wt \theta_k) \in \mathbb{X}$
is introduced to measure the small deformations,
and a universal bound will be obtained for the corresponding  determinant of  Jacobian matrix.

More precisely,
given any $L>0$,
$\omega_k>0$, $x_k\in \R^d$, $\vartheta_{k}\in\R$, $1\leq k\leq K$,
we set
\begin{align} \label{SL-def}
S_L(x):=\sum_{k=1}^{K}S_{k,L}(x)
=\sum_{k=1}^{K}(\omega_kL)^{-\frac d2} Q_{0,k,L}\(\frac{x-x_k}{\omega_kL}\)e^{i\t_{0,k}},
\end{align}
where
\begin{align}
   Q_{0,k,L}(y) = Q(y)e^{i(\beta_{0,k}\cdot y - \frac 14 \g_{0,k} |y|^2)}
\end{align}
with $\beta_{0,k} = 0$,
$\g_{0,k} = \omega_k^2 L$,
and $\theta_{0,k} = \omega^{-2}_{k} L^{-1} + \vartheta_k$,
$1\leq k\leq K$.
Note that,
in the case $L= T-t$,
$S_{k,L}$ is exactly the pseudo-conformal blow-up solution $S_k$
given by \eqref{Sj-blowup},
and thus $S_L = \sum_{k=1}^K S_k$.

The proof of Theorem \ref{Thm-geometri-dec} is based on Lemma \ref{Lem-Jacob} below.
\begin{lemma} \label{Lem-Jacob}
There exists a universal small constant $\delta_*>0$ such that
the following holds.
For any $0<\delta, L<\delta_*$ and for any $v$ satisfying $\|v-S_L\|_{L^2}\leq \delta$,
there exist unique $C^1$ parameters
$\wt \calp(v)=(\wt \lbb, \wt \a, \wt \b, \wt \g, \wt \t)
\in \mathbb{X}$ with respect to $v$,
such that
$v$ admits the decomposition
\begin{align} \label{v-dec-L}
     v =\sum_{k=1}^{K}(\wt \lbb_k \omega_k L)^{-\frac d2}
         {Q}_{k,L}\(\frac{x-x_k-\wt \a_k \omega_k L}{\wt \lbb_k \omega_k L}\)
         e^{i (\wt \t_{k}+\omega_k^{-2} L^{-1}+ \vartheta_k)}+ R_L
        \ \(=:\sum\limits_{k=1}^KU_{k,L} +R_L\)
\end{align}
and the following orthogonality conditions hold:
for each $1\leq k\leq K$,
\be\ba \label{Ortho-cond-L-Lemma}
&{\rm Re}\int (x- x_k - \wt {\a}_{k}\omega L) U_{k,L}\ol{R_L} dx=0,\ \
{\rm Re} \int |x- x_k - \wt {\a}_{k}\omega L|^2 U_{k,L}  \ol{R_L} dx=0,\\
&{\rm Im}\int \nabla U_{k,L}  \ol{R_L} dx=0,\ \
{\rm Im}\int \Lambda_k U_{k,L}  \ol{R_L} dx=0,\ \
{\rm Im}\int {\varrho}_{k,L}  \ol{R_L} dx=0.
\ea\ee
Here, for each $1\leq k\leq K$,
\begin{align}
   &{Q}_{k,L}(y)=Q (y) e^{i\b_k\cdot y-i\frac{\g_k}{4}|y|^2},    \label{Qj-Q}  \\
  {\varrho}_{k,L}(x)= {\lbb}^{-\frac d2}_{k} & \rho_{k}\(\frac{x- {\a}_{k}}{{\lbb}_k}\) e^{i{\theta}_{k}}\ with\
    {\rho}_k(y)=\rho(y) e^{i\b_k\cdot y-i\frac{\g_k}{4}|y|^2}, \label{rhok-rho-L}
\end{align}
where $\rho$ is given by \eqref{rho-def} in Appendix,
and the parameters $(\lbb_k, \a_k, \beta_k, \g_k, \theta_k)$ depend on
$(\wt \lbb_k, \wt \a_k, \wt \beta_k, \wt \g_k, \wt \theta_k)$ and $L$ as follows
\be\ba \label{calp-wtcalp}
  &{\lbb}_k=\wt \lbb_k \omega_k  L,\ \  {\a}_k=x_k+ \wt \a_k  \omega_k L, \ \ \beta_k = \wt \beta_k,    \\
  & \g_k = \wt \g_k + \omega_k^2 L, \ \ \theta_k = \wt \theta_k + \omega_k^{-2} L^{-1} + \vartheta_k, \ \ 1\leq k\leq K.
\ea\ee
Moreover, there exists a universal constant $C>0$ such that
\begin{align}
& \sum_{k=1}^{K}(|\wt \lbb_k-1|+|\wt \a_k|+|\wt \b_k|+|\wt \g_k|+|\wt \t_k|)
  \leq C \|v- S_L\|_{L^2}, \label{wtcalp-wtR-esti} \\
& \|R_L\|_{L^2}  \leq C \|v-S_L\|_{L^2}, \label{RL-v-SL} \\
&  L \|\na R_L\|_{L^2} \leq C (\|v-S_L\|_{L^2} + L\|\na v-\na S_L\|_{L^2}). \label{naRL-v-SL}
\end{align}
\end{lemma}

{\bf Proof.}
{\bf  Set-up of the notations.}
To simplify the notations, we write $U_L:= \sum_{k=1}^K U_{k,L}$,
$R:= R_L$,
and
\begin{align} \label{wtR-def}
   \wt R:= v-S_L.
\end{align}

Set  $\mathbb{Y}:= \bbr \times \bbr^{d} \times \bbr^{d} \times \bbr \times \bbr$,
$\wt \calp_{0,j}: =(1,0,0,0,0) \in \mathbb{Y}$
and $\wt \calp_0 = (\wt \calp_{0,1}, \cdots, \wt \calp_{0,K})\in \mathbb{Y}^K$.
Similarly,
let $\wt \calp_j:= (\wt \lbb_j, \wt \a_j, \wt \beta_j, \wt \g_j, \wt \theta_j) \in \mathbb{Y}$,
$\wt \calp:= (\wt \calp_1, \cdots, \wt \calp_K) \in \mathbb{Y}^K$.
Let $B_\delta(v_0, \wt \calp_0)$ denote
the closed ball centered at $(v_0, \wt\calp_0)$ of radius $\delta$,
i.e.,
\begin{align}
  B_\delta(v_0, \wt \calp_0) :=\{(v, \wt \calp):\ \|v-v_0\|_{L^2}\leq \delta,\ \ |\wt \calp- \wt \calp_0|\leq \delta\},
\end{align}
where $\delta$ is a small constant to be chosen later,
and
\begin{align}
   |\wt \calp-\wt \calp_0|:=
   \sum\limits_{j=1}^K |\wt \calp_j - \wt \calp_{0,j}|
   =\sum_{j=1}^{K}(|\wt \lbb_j-1|+|\wt \a_j|+|\wt \b_j|+|\wt \g_j|+|\wt \t_j|).
\end{align}

For every $1\leq k\leq K$, define the functionals by
\begin{align*}
& f^k_{1}(v, \wt \calp) := \lbb_k^{-2} {\rm Re} \int |x- {\a}_{k}|^2 U_{k,L} \ol{R}dx, \ \
  f^k_{2,i}(v,\wt \calp) :=\lbb_k^{-1}{\rm Re}\int (x_i- {\a}_{k,i}) U_{k,L}\ol{R}dx,  \nonumber \\
&f^k_{3,i}(v,\wt \calp) : =\lbb_k {\rm Im}\int \partial_i U_{k,L} \ol{R}dx,\ \
 f^k_{4}(v,\wt \calp) := {\rm Im}\int\Lambda_k U_{k,L} \ol{R}dx,\ \
 f^k_{5}(v,\wt \calp) : ={\rm Im}\int {\varrho}_{k,L} \ol{R}dx,
\end{align*}
where $1\leq i\leq d$,
and the parameters $\lbb_k, \a_k, \beta_k, \g_k, \theta_k$ are defined as in \eqref{calp-wtcalp}.

Let $F^k:= (f^k_{1},f^k_{2,1}, f^k_{2,2}, \cdots, f^k_{5})$ and
$\frac{\partial F^k}{\partial \wt \calp_j}$ denote the Jacobian matrix of $F^k$ with respect to $\wt \calp_j$
\begin{align}
\frac{\partial F^k}{\partial \wt \calp_j}:=\left(
                                           \begin{array}{ccc}
                                             \frac{\partial f^k_{1}}{\partial \wt \lbb_j} & \cdots & \frac{\partial f^k_{1}}{\partial \wt\t_j} \\
                                             \vdots &  & \vdots \\
                                             \frac{\partial f^k_{5}}{\partial \wt \lbb_j} &  \ldots &\frac{\partial f^k_5}{\partial \wt\t_j} \\
                                           \end{array}
                                         \right), \ \ 1\leq j, k\leq K.
\end{align}
Similarly,
let $F:=(F^1, \cdots, F^K)$
and $\frac{\partial F}{\partial \wt  \calp} := (\frac{\partial F^k}{\partial \wt \calp_j})_{1\leq j,k\leq K}$
be the corresponding Jacobian matrix.

Note that,
for every $(v, \wt \calp)\in B_\delta(S_L, \wt \calp_0)$,
since
\begin{align*}
   R = v- U_L = \wt R + S_L - U_L,
\end{align*}
we have
\begin{align} \label{R-wtR-SLUL}
\|R\|_{L^2} \leq \|\wt R\|_{L^2} + \|S_L - U_L\|_{L^2}.
\end{align}
Using the explicit expressions of $S_L$ and $U_L$
in \eqref{SL-def} and \eqref{v-dec-L}, respectively,
we compute that for some $C>0$
(see the proof of \cite[Theorem 2.15]{SZ19}
and \cite[Theorem 2.14]{SZ20}),
\begin{align} \label{SL-UL-diff}
\|S_L-U_L\|_{L^2}
\leq& C\sum_{j=1}^{K} \bigg(\bigg|\frac{{\lambda}_j^{\frac{d}{2}}-(\omega_jL)^{\frac{d}{2}}}{(\omega_jL)^{\frac{d}{2}}}\bigg|
+\bigg|\frac{{\alpha}_j-x_j}{\lbb_j}\bigg|
+\(\frac{{\lbb}_j}{\omega_jL}\)^{\frac d2}|\wt \beta_j|
+\(\frac{{\lbb}_j}{\omega_jL}\)^{\frac d2}|\wt \gamma_j|
   +|\wt \theta_j|
   + \bigg|\frac{\omega_j L}{\lbb_j} -1 \bigg|\bigg)  \nonumber  \\
   \leq& C |\wt \calp - \wt \calp_0|.
\end{align}
Thus,
we infer that there exists a universal constant $\wt C>0$
such that
\begin{align} \label{R-wtR-PP0}
\|R\|_{L^2}\leq \wt C (\|\wt R\|_{L^2} + |\wt \calp - \wt \calp_0|)
           \leq 2 \wt C \delta,\ \ \forall (v, \wt \calp)\in B_\delta(S_L, \wt \calp_0).
\end{align}

{\bf Step 1. Nondegeneracy of Jacobian. }
We claim that, there exist small constants $\delta_*, c_1, c_2>0$
such that for any $0<\delta, L \leq \delta_*$,
\begin{align} \label{Jacob-posit}
   0<c_1\leq \bigg|\det \frac{\partial F}{\partial \wt \calp}(v, \wt \calp) \bigg| \leq c_2<\9, \ \ \forall (v, \wt \calp)\in B_{\delta}(S_L, \wt \calp_0).
\end{align}

In order to prove \eqref{Jacob-posit},
using the results in the proof of \cite[Lemma 4.2]{SZ19}
and the chain rule we compute that for any $(v,\wt \calp)\in B_\delta(S_L, \wt \calp_0)$
and for any $1\leq k\leq K$, $1\leq i,h\leq d$,
\begin{align}
  & \pa_{\wt \lbb_k} f^k_{1}  = - \frac{1}{\wt \lbb_k} \|xQ\|_{L^2}^2  + \mathcal{O}(\|R\|_{L^2}),  \ \
    \pa_{\wt \a_{k,h}} f^k_{2,i}  = -\frac {\delta_{ih}}{2\wt \lbb_k}  \|Q\|_{L^2}^2 + \mathcal{O}(\|R\|_{L^2}),\   \label{Jacob-f1f2} \\
  & \pa_{\wt \lbb_k} f^k_{3,i}   =  \frac{ \beta_{k,i}}{2\wt \lbb_k}\|Q\|_{L^2}^2 + \mathcal{O}( \|R\|_{L^2}),\ \
    \pa_{\wt \beta_{k,h}} f^k_{3,i}  = -\frac{\delta_{ih}}{2} \|Q\|_{L^2}^2 + \mathcal{O}(\|R\|_{L^2}), \label{Jacob-f4} \\
   &  \pa_{\wt \a_{k,h}} f^k_{4}   = - \frac{\beta_{k,h}}{2\wt \lbb_k} \|Q\|_{L^2}^2  + \mathcal{O}(\|R\|_{L^2}),\ \
   \pa_{\wt \g_k} f^k_{4}   =  \frac{1}{4} \|xQ\|_{L^2}^2  + \mathcal{O}(\|R\|_{L^2}),  \label{Jacob-f3} \\
  & \pa_{\wt \lbb_k} f^k_{5}  =  \frac{\g_k}{2 \wt \lbb_k}   \<\rho, |x|^2 Q\> + \mathcal{O}(\|R\|_{L^2}) ,  \ \
    \pa_{\wt \a_{k,h}} f^k_{5} = \frac {\beta_{k,h}}{2 \wt \lbb_k}   \|xQ\|_{L^2}^2 +  \mathcal{O}(\|R\|_{L^2}),  \label{Jacob-f5.1} \\
  & \pa_{\wt \g_k} f^k_{5} = -\frac 14\<\rho, |x|^2Q\>  +  \mathcal{O}(\|R\|_{L^2}) ,   \ \
    \pa_{\wt \theta_k} f^k_{5} = -\frac 12 \|xQ \|_{L^2}^2 +  \mathcal{O}(\|R\|_{L^2}) . \label{Jacob-f5.2}
\end{align}
Moreover,
using the exponential decay of $Q$ and $\rho$
in \eqref{Q-decay} and \eqref{rho-decay}, respectively,
we infer that,
there exists $\mu>0$ such that
for any $1\leq j\not = k\leq K$,
\begin{align} \label{Jacob-Fk-Pj-exp}
     |\pa_{\wt \lbb_j} F^k|   +|\pa_{\wt \a_{j,h}} F^k|
    +|\pa_{\wt \b_{j,h}} F^k|+|\pa_{\wt\g_j} F^k| +|\pa_{\wt \t_j} F^k|
      =\calo (e^{-\frac{\mu |x_j-x_k|}{L}}).
\end{align}

Thus, we conclude that for any $(v, \wt \calp) \in B_\delta(S_L, \wt \calp_0)$,
\begin{align}
  \bigg|\det \frac{\partial F^k}{\partial \wt \calp_k} (v, \wt  \calp)\bigg|
=& 2^{-(3+2d)}  \wt \lbb_k^{-(1+d)}  \|Q\|_{L^2}^{4d}  \|xQ\|_{L^2}^6
       + \mathcal{O}(\|R\|_{L^2}),  \ \ 1\leq k\leq K,  \label{Jacob-FjPj} \\
  \bigg|\det \frac{\partial F^k}{\partial \wt \calp_j} (v, \wt  \calp)\bigg|
 = &  \calo(e^{-\frac{\mu |x_j-x_k|}{L}}), \ \ j\not =k.  \label{Jacob-FjPk}
\end{align}
Taking into account \eqref{R-wtR-PP0},
we obtain that
for $\delta, L$ small enough and
for any $(v, \wt \calp) \in B_\delta(S_L,  \wt  \calp_0)$,
\begin{align} \label{Jacob-FP}
    \bigg|\det \frac{\partial F}{\partial \wt \calp} (v,  \wt \calp)\bigg|
    =& 2^{-(3+2d)K}
      \prod\limits_{k=1}^K  \wt \lbb_k^{-(1+d)}
      \|Q\|_{L^2}^{4dK} \|xQ\|_{L^2}^{6K}
      + \calo(\|R\|_{L^2} + e^{-\frac{\mu \wt \sigma}{L}}) \nonumber \\
    =&  2^{-(3+2d)K}
      \prod\limits_{k=1}^K  \wt \lbb_k^{-(1+d)}
      \|Q\|_{L^2}^{4dK} \|xQ\|_{L^2}^{6K}
      + \calo(\delta + e^{-\frac{\mu \wt \sigma}{L}}),
\end{align}
where $\mu>0$ is a universal small constant,
$\wt \sigma:= \min_{j\not = k}\{|x_j-x_k|\}>0$.
This yields \eqref{Jacob-posit}, as claimed.

{\bf Step 2. Estimate of modulation parameters.}
Below we prove that
there exists a universal constant $C_*(\geq 1)$ such that, for any
$0<\delta, L\leq \delta_*$
and for any
$(v_1, \wt \calp(v_1)), (v_2, \wt \calp(v_2))\in B_\delta(S_L, \wt \calp_0 )$,
if $F(v_1, \wt \calp(v_1))= F(v_2, \wt \calp(v_2))=0$,
then
\begin{align} \label{Pv1-Pv2}
|\wt \calp(v_1)-\wt \calp(v_2)|\leq C_* \|v_1-v_2\|_{L^2}.
\end{align}

For this purpose,
we note that
\begin{align*}
 0=& F(v_1, \wt \calp(v_1)) -  F(v_2, \wt \calp(v_2))   \nonumber \\
  =& (F(v_1, \wt \calp(v_1)) -  F(v_1, \wt \calp(v_2)) )
    + (F(v_1,\wt  \calp(v_2)) - F(v_2, \wt \calp(v_2)) ),
\end{align*}
which yields immediately that
\begin{align} \label{Fv12-Fv21}
F(v_1, \wt \calp(v_1))-F(v_1, \wt \calp(v_2))= F(v_2, \wt \calp(v_2) ) - F(v_1, \wt \calp(v_2)).
\end{align}
By the differential mean value theorem,
this yields that
\begin{align} \label{Fv12-pFpP}
   \(\frac{\partial F}{\partial \wt \calp}(v_1, \wt \calp_{r})\) (\wt \calp(v_1)-\wt \calp(v_2))^t
   = (F(v_2, \wt \calp(v_2))-F(v_1,\wt  \calp(v_2)))^t,
\end{align}
where $\wt \calp_r=r\wt \calp(v_1)+(1-r)\wt \calp(v_2)$
for some $0<r<1$,
and the superscript $t$ means the transpose of matrices.

In view of \eqref{Jacob-posit},
the Jacobian matrix
$\frac{\partial F}{\partial \wt \calp}(v_1, \wt \calp_r)$ is invertible,
and thus
we may reformulate equation \eqref{Fv12-pFpP} as follows
\begin{align} \label{Fv12-Fv21-A}
   (\wt \calp(v_1)-\wt \calp(v_2))^t= \(\frac{\partial F}{\partial \wt \calp}(v_1,\wt \calp_r) \)^{-1} (F(v_2, \wt \calp(v_2)) - F(v_1, \wt \calp(v_2)))^t.
\end{align}
Then, using \eqref{Jacob-f1f2}-\eqref{Jacob-Fk-Pj-exp}
we get that for a universal constant $C>0$,
\begin{align} \label{Fp-Jacob-invert}
   \bigg\|\bigg(\frac{\partial F}{\partial \wt \calp}(v_1, \wt \calp_r)\bigg)^{-1} \bigg\|\leq C,
\end{align}
where $\|\cdot\|$ denotes the Hilbert-Schmidt norm of matrices.

Moreover,
by the exponential decay of $Q$ and $\rho$ in \eqref{Q-decay} and \eqref{rho-decay}, respectively,
\begin{align} \label{f-v2v1-diff}
&|F(v_2, \wt \calp(v_2)) - F(v_1, \wt \calp(v_2))|\leq C \|v_2-v_1\|_{L^2}.
\end{align}

Thus, combining \eqref{Fv12-Fv21-A}, \eqref{Fp-Jacob-invert} and \eqref{f-v2v1-diff} altogether
we prove \eqref{Pv1-Pv2}.

{\bf Step 3. Proof of main results. }
Let $\delta_*, C_*$ be the universal constants
as in Step 1 and Step 2, respectively.
Define the set $B$ in the ball $B_{\frac{\delta_*}{C_*}}(S_L)$
by
\begin{align}
B:=\{ v\in B_{\frac{\delta_*}{C_*}}(S_L): \exists \wt \calp \in B_{\delta_*}(\wt \calp_0),
     \ such\ that\ F(v,\wt \calp) =0\}.
\end{align}

In order to prove the geometrical decomposition \eqref{v-dec-L},
it suffices to prove that
\begin{align} \label{B-BSL}
   B = B_{\frac{\delta_*}{C_*}}(S_L).
\end{align}
For this purpose,
since $B_{\frac{\delta_*}{C_*}}(S_L)$ is connected
and $S_L \in B$ due to the fact $F(S_L, \wt \calp_0) =0$,
we only need to prove that
$B$ is both open and closed in $B_{\frac{\delta_*}{C_*}}(S_L)$.

To this end,  suppose that $v\in B$,
then by definition
there exists $\wt \calp(v) \in B_{\delta_*}(\wt \calp_0)$ such that $F(v,\wt \calp(v))=0$.
By \eqref{Jacob-posit}, the corresponding Jacobian matrix at $(v, \wt \calp(v))$ is non-degenerate.
Thus,  the implicit function theorem yields the existence of
a small open neighborhood $\calu (v)$ of $v$ in the ball $B_{\frac{\delta_*}{C_*}}(S_L)$
such that $\calu (v) \subseteq B$,
and so the set $B$ is open in $B_{\frac{\delta_*}{C_*}}(S_L)$.

It remains to show that $B$ is also closed in $B_{\frac{\delta_*}{C_*}}(S_L)$.

To this end,
take any sequence $\{v_n\} \subseteq B$
such that $v_n \to v_*$ in $L^2$ for some $v_*\in B_{\frac{\delta_*}{C_*}}(S_L)$.
Then, by definition
there exist modulation parameters $\wt \calp(v_n)\in  B_{\delta_*}(\wt \calp_0)$, $n\in \mathbb{N}$,
such that $F(v_n, \wt \calp(v_n))=0$.

In particular,
$\{\wt \calp(v_n)\}$ is uniformly bounded in the finite dimensional space $\mathbb{Y}^K$.
This yields that along a subsequence (still denoted by $\{n\}$),
$\wt \calp(v_{n}) \to \wt \calp_*\ (\in B_{\delta_*}(\wt \calp_0))$
for some $\wt \calp_* \in \mathbb{Y}^K$.

Let $U_{k,L, \wt \calp(v_n)}$ and $U_{k,L,\wt \calp_*}$
be the $k$-th blow-up profiles corresponding to $\wt \calp(v_n)$
and $\wt \calp_*$, respectively.
Similarly denote $\varrho_{k,L, \wt \calp(v_n)}$ and $\varrho_{k,L,\wt \calp_*}$.
Then,
using the explicit expressions \eqref{v-dec-L} and \eqref{rhok-rho-L}
we infer that
$\<x\>^2U_{k,L, \wt \calp(v_{n})} \to \<x\>^2 U_{k,L, \wt \calp_*}$ in $H^1$,
$\varrho_{k,L, \wt \calp(v_{n})} \to \varrho_{k,L, \wt \calp_*}$ in $L^2$,
and $v_{n} - \sum_{k=1}^KU_{k,L, \wt \calp(v_{n})}  \to v_* - \sum_{k=1}^K U_{k,L,\wt \calp_*}$ in $L^2$.
Since $F(v_{n}, \wt \calp(v_{{n}})) =0$,
letting $n\to \9$
we thus obtain
$F(v_*, \wt \calp_*) = 0$.
In particular, this yields that $v_*\in B$.
Thus,
$B$ is a closed set in $B_{\frac{\delta_*}{C_*}}(S_L)$
and \eqref{B-BSL} is proved.

Therefore, we obtain the geometrical decomposition \eqref{v-dec-L}
satisfying the orthogonality conditions in \eqref{Ortho-cond-L-Lemma}.
Moreover, estimate \eqref{wtcalp-wtR-esti} follows from \eqref{Pv1-Pv2}
by taking $v_1= v$ and $v_2 = S_L$.
Estimate   \eqref{RL-v-SL} then
follows from \eqref{wtcalp-wtR-esti}, \eqref{R-wtR-SLUL} and \eqref{SL-UL-diff}.

Regarding estimate \eqref{naRL-v-SL},
similarly to \eqref{SL-UL-diff}, we have
\begin{align} \label{naUL-naSL-L2}
    \|\na U_L - \na S_L\|_{L^2}
   \leq& C \sum\limits_{k=1}^K
        \bigg(  \frac{1}{\omega_k L} \bigg| \frac{\omega_k L}{\lbb_{k}} -1\bigg|
             + \bigg|\frac{\a_{k}-x_k}{\omega_k L\lbb_{k}}\bigg|
        + \bigg| \frac{\beta_{k}}{\omega_k L} \bigg|
         +  \bigg|\frac{\g_{k} - \omega_k^2 L}{\omega_k L} \bigg| \nonumber  \\
       & \qquad  \ \
         + \bigg|\frac{\lbb_{k}^{1+\frac d2} -(\omega_k L)^{1+\frac d2}}{\lbb_{k} (\omega_k L)^{1+\frac d2}} \bigg|
           + \bigg|\frac{\theta_{k} - \omega_k^{-2}L^{-1} - \vartheta_k}{\lbb_k} \bigg|   \bigg)  \nonumber \\
   \leq& C L^{-1} |\wt \calp - \wt \calp_0|,
\end{align}
which along with \eqref{wtcalp-wtR-esti} yields that
\begin{align}
    \|\na R_L\|_{L^2}
    \leq& \|\na v - \na S_L \|_{L^2} + \|\na S_L - \na U_L\|_{L^2}
    \leq   C (\|\na v - \na S_L \|_{L^2} + L^{-1} \| v- S_L\|_{L^2}),
\end{align}
thereby proving \eqref{naRL-v-SL}.
Therefore, the proof of Lemma \ref{Lem-Jacob} is complete.
\hfill $\square$.

{\bf Proof of Theorem \ref{Thm-geometri-dec}.}
We take $L=T-t$ in Lemma \ref{Lem-Jacob} and get
$S_L = \sum_{k=1}^K {S_k}$,
where $\{S_k\}$ are given by \eqref{Sj-blowup}.
By \eqref{v-S-H1-o1-Thm},
$\|v-S_L\|_{L^2} =o(1)$.
This yields that
$\|v-S_L\|_{L^2} \leq \delta_*$  for any $t$ close to $T$,
where $\delta_*$ is as in Lemma \ref{Lem-Jacob}.
Thus, Lemma \ref{Lem-Jacob} gives the existence of geometrical parameters
$\calp =(\lbb, \a, \beta, \g, \theta) \in \mathbb{X}$,
such that the geometrical decomposition \eqref{v-dec} and
the orthogonality conditions in \eqref{ortho-cond-Rn-wn} hold.
Moreover, the estimates \eqref{lbb-a-t-0}-\eqref{R-t-0} follow directly from \eqref{calp-wtcalp}-\eqref{naRL-v-SL}.

The $C^1$-regularity of $\calp$ can be proved by using the regularization arguments as in \cite{MM}.
More precisely, we may take a sequence $v_n(T^*)\in H^3$, $n\geq 1$,
such that $v_n(T^*) \to v(T^*)$ in $H^1$ as $n\rightarrow\infty$.
Let $v_n$ be the solutions to \eqref{equa-NLS} corresponding to $v_n(T^*)$, $n\geq 1$.
Then, by the local well-posedness theory,
for any $\tilde t \in [T^*, T)$,
we have that,
for $n$ large enough,
$v_n\in C([T^*,\tilde t];H^3)\cap C^1([T^*,\tilde t];H^1)$ and $v_n\rightarrow v$ in $C([T^*,\tilde t];H^1)$ as $n\rightarrow\infty$.
Let $\calp_n$ and $R_n$ be the geometrical parameter and remainder corresponding to $v_n$, $n\geq 1$.
Using Lemma \ref{Lem-Jacob} and \eqref{Pv1-Pv2},
we infer that for $n$ large enough, $\calp_n\in C^1([T^*,\tilde t];\mathbb{X})$
and  $\calp_n\rightarrow \calp$ in $C([T^*,\tilde t]; \mathbb{X})$ as $n\rightarrow\infty$.
Moreover, the computations of modulation equations in Subsection \ref{Subsec-mod-equa} below show that
$\dot{\calp}_n$ can be expressed in terms of $\calp_n$ and the remainder $R_n$,
and thus $\dot{\calp}_n$ converges in $C([T^*,\tilde t]; \mathbb{X})$.
This yields that $\calp\in C^1([T^*,\tilde t];\mathbb{X})$.
Thus, the $C^1$-regularity of $\calp$ on $[T^*, T)$ follows,
due to the arbitrariness of $\tilde t$.

Therefore, the proof of Theorem \ref{Thm-geometri-dec} is complete.
\hfill $\square$

Below we introduce the localized remainder $R_k$
and the localization function $\Phi_k$
as in \cite{ SZ20}, $1\leq k\leq K$.
We also refer the readers to \cite{CF20,Ma05, MMT06},
where different types of localization functions are introduced in the setting of multi-solitons.

Because \eqref{equa-NLS} is invariant under the orthogonal transform,
we may take an orthonormal basis $\{\textbf{v}_j\}_{j=1}^d$ of $\bbr^d$ as in \cite{MM06},
such that $(x_j-x_k)\cdot \textbf{v}_1 \not = 0$ for any $j\not = k$.
Without loss of generality,
we may assume that $x_1\cdot \textbf{v}_1<x_2\cdot \textbf{v}_1<\cdots<x_K\cdot \textbf{v}_1 $.
Then,
\begin{align} \label{sep-xj-0}
\sigma :=\frac{1}{12}\min_{1\leq k \leq K-1}\{(x_{k+1}-x_k)\cdot \textbf{v}_1\}> 0.
\end{align}

Let $\Phi(x)$ be a smooth function on $\R^d$ such that $0\leq \Phi(x)\leq 1$,
$\Phi(x)=1$ for $x\cdot \textbf{v}_1\leq 4\sigma$,
$\Phi(x)=0$ for $x\cdot \textbf{v}_1\geq 8\sigma$
and $\|\nabla \Phi\|_{L^\9}\leq C\sigma^{-1}$.
The localization functions $\Phi_k$ are defined by
\be\ba \label{phi-local}
&\Phi_1(x) :=\Phi(x-x_1), \ \ \Phi_K(x) :=1-\Phi(x-x_{K-1}),  \\
&\Phi_k(x) :=\Phi(x-x_{k})-\Phi(x-x_{k-1}),\ \ 2\leq k\leq K-1.
\ea\ee
In particular,
we have the  partition of unity $1= \sum_{k=1}^K \Phi_k$
and the decomposition
\begin{align} \label{R-Rk}
   R = \sum_{k=1}^K R_k,\ \ with\ \ R_k:= R\Phi_k.
\end{align}

As a consequence of Theorem \ref{Thm-geometri-dec}
and the decoupling Lemma \ref{Lem-inter-est} in Appendix,
the following almost orthogonality  holds for the localized remainders.

\begin{lemma} (Almost orthogonality)  \label{Lem-almost-orth}
There exist $C, \delta >0$ such that for
$1\leq k\leq K$ and
$t$ close to $T$,
\be\ba \label{Orth-almost}
&|{\rm Re}\int (x-\a_{k}) U_{k} \ol{R_{k}}dx|
  + |{\rm Re} \int |x-\a_{k}|^2 U_{k} \ol{R_{k}}dx|\leq Ce^{-\frac{\delta}{T- t}}\|R\|_{L^2},  \\
& |{\rm Im}\int \nabla U_{k} \ol{R_{k}}dx|
  + |{\rm Im}\int \Lambda_k U_{k} \ol{R_{k}}dx|
  + |{\rm Im}\int  \varrho_{k} \ol{R_{k}}dx|\leq Ce^{-\frac{\delta}{T- t}}\|R\|_{L^2}.
\ea\ee
\end{lemma}

\subsection{Modulation equations} \label{Subsec-mod-equa}

The leading order of the dynamics of geometrical parameters are characterized by
the \emph{vector of modulation equations}, defined by
\begin{align*}
   Mod_{k}:= |\lambda_{k}\dot{\lambda_{k}}+\gamma_{k}|+|\lambda_{k}^2\dot{\gamma}_{k}+\gamma_{k}^2|
   +|\lambda_{k}\dot{\alpha}_{k}-2\beta_{k}|
                +|\lambda_{k}^2\dot{\beta}_{k}+\gamma_{k}\beta_{k}|
              +|\lambda_{k}^2\dot{\theta}_{k}-1-|\beta_{k}|^2|,
\end{align*}
where $\dot{\lbb_k}:= \frac{d}{dt}\lbb_k$
and similar notations apply to the remaining four geometrical parameters.
Set
\begin{align*}
Mod:=\sum_{k=1}^{K}Mod_{k}.
\end{align*}
In particular,
for the pseudo-conformal blow-up solution $S_k$ given by \eqref{Sj-blowup},
we have $Mod =0$.

We also use the notation
\begin{align} \label{P-def}
   P := \sum_{k=1}^K (|\lbb_{k}|+ |\a_{k}-x_k|  + |\beta_{k}| + |\g_{k}|).
\end{align}

The main control of modulation equations is presented below.
\begin{theorem} (Control of modulation equations) \label{Thm-Mod}
There exists $C>0$ such that for $t$ close to $T$,
\begin{align} \label{Mod-bdd}
Mod(t)\leq
C\( \sum_{k=1}^{K}|M_k(t)| + P^2 D(t)
+ D^2(t) + e^{-\frac{\delta}{T-t}}\),
\end{align}
where the quantity
\begin{align}   \label{D-def}
D(t):= \|R(t)\|_{L^2} + (T-t)\|\nabla R(t)\|_{L^2},
\end{align}
and for every $1\leq k\leq K$,
\begin{align} \label{Mj-def-Mod}
   M_k(t):= 2 {\rm Re} \<R_k(t), U_k(t)\> + \int |R(t)|^2 \Phi_k dx,\ \ with\ \ R_k = R \Phi_k.
\end{align}
Moreover, the following identity holds
\begin{align}  \label{lbb2g+g2-f''Qve2}
 \frac{1}{4}\|yQ\|_2^2(\lbb_k^2\dot{\g_k}+\g_k^2)
=& M_k - \int |R|^2\Phi_kdx+{\rm Re}\int
   (1+\frac 2d)|Q_k|^{\frac {4}{d}}|\varepsilon_k|^2
   + \frac 2d |Q_k|^{\frac 4d-2}Q_k^2 \ol{\varepsilon_k}^2 dy  \nonumber \\
&+{\rm Re}\int
   (y \cdot\nabla \ol{Q_k}) (f''(Q_k)\cdot \ve_k^2) dy \nonumber \\
  & +\calo\((P+\|R\|_{L^2} + e^{-\frac{\delta}{T-t}}) Mod
          + P^2 \|R\|_{L^2}
               + D^3
               +e^{-\frac{\delta}{T-t}}\),
\end{align}
where the renormalized remainder $\ve_k$ is defined by
\begin{align} \label{Rj-ej}
   R_k(t,x) = \lbb_k^{-\frac d2} \ve_k\(t, \frac{x-\a_k}{\lbb_k}\) e^{i\theta_k},
\end{align}
and
\begin{align} \label{f''Qk-vek2}
   f''(Q_k)\cdot \ve_k^2:=
   \frac 2d(1+\frac 2d) |Q_k|^{\frac 4d-2} Q_k |\ve_k|^2
   + \frac 1d(1+\frac 2d) |Q_k|^{\frac 4d-2} \ol{Q_k} \ve_k^2
   +\frac 1d(\frac 2d-1) |Q_k|^{\frac 4d-4} Q_k^3 \ol{\ve_k}^2.
\end{align}
\end{theorem}

\begin{remark}
The identity
\eqref{lbb2g+g2-f''Qve2} is important in the derivation of the
monotonicity of (modified) localized virial functionals (See Theorems \ref{Thm-Loc-virial} and \ref{Thm-Modi-Loc-virial}).
In fact, this specific algebraic identity provides the necessary  quadratic terms
to perform the localized coercivity of linearized operators.
\end{remark}

The proof of Theorem \ref{Thm-Mod} is similar to that of \cite[Proposition 4.3]{SZ20}
and thus is postponed to the Appendix for the simplicity of exposition.

We end this section with the preliminary estimates of remainder related to the quantity $D$.

\begin{lemma}  \label{Lem-D-R}
(Preliminary estimates of remainder)
We have that for $t$ close to $T$,
\begin{align} \label{R-D}
   \|R(t)\|_{H^1} \leq C \frac{D(t)}{T-t}, \ \
   \|R(t)\|_{L^2} \|\na R(t)\|_{L^2} \leq  \frac{D^2(t)}{T-t},
\end{align}
and
\begin{align} \label{R-D-Lp}
    \|R(t)\|^p_{L^p} \leq C (T-t)^{-d(\frac p2 - 1)} D^p(t), \ \ 2\leq p<\9.
\end{align}
Moreover,  we have
\begin{align}  \label{D-0}
    D(t)=o(1),\ \ for\ t\ close\ to\ T.
\end{align}
\end{lemma}

{\bf Proof.}
Estimate \eqref{R-D} follows directly from the definition of $D$,
and \eqref{R-D-Lp} follows from the Gagliardo-Nirenberg inequality \eqref{G-N} below.
Estimate \eqref{D-0} is a consequence of the preliminary estimate of remainder in \eqref{R-t-0}.
\hfill $\square$

\section{Localized mass, energy and localized virial functionals} \label{Sec-Energy-Locvirial}

This section mainly contains the crucial controls of localized mass,
energy and (modified) localized virial functionals.

We first deal with the sharp estimate of localized mass
and its relationship with the remainder in Subsection \ref{Subsec-loc-mass-energy}.
The coercivity type control of energy will be also proved there,
which allows to upgrade the convergence rate of remainder to the first order,
i.e., $D=\calo(T-t)$.
Then, Subsection \ref{Subsec-Modif-Loc-Vir} contains the key monotonicity of
(modified) localized virial functionals,
which yield the refined space time estimate of remainder and the precise asymptotic order of parameter $\g$.
By virtue of these estimates,
in Subsection \ref{Subsec-Energy-Quant}
we are able to prove the $H^1$ dispersion of remainder along a sequence and
obtain the energy quantization,
which is the key towards the derivation of refined energy estimate.

\subsection{Localized mass and energy}  \label{Subsec-loc-mass-energy}

Recall from \eqref{Mj-def-Mod} that
\begin{align} \label{Mj-def}
 M_k=2{\rm Re} \< R_{k}, {U_{k}}\>+\int |R|^2\Phi_kdx,
\end{align}
where $R_k=R\Phi_k$
and $\{\Phi_k\}$ are the localization functions given by \eqref{phi-local}.
Let $\wt R:= v-\sum_{k=1}^K S_k$
with $S_k$ given by \eqref{Sj-blowup}, $1\leq k\leq K$.

The main estimates of localized mass
are formulated in Theorem \ref{Thm-Loc-Mass} below.

\begin{theorem} (Control of localized mass) \label{Thm-Loc-Mass}
For every $1\leq k\leq K$ and for any $t$ close to $T$,
we have
\begin{align}
    |M_k(t)| \leq& C
                    \int_t^T \int_s^T \|\wt R\|_{H^1}^2 drds +  Ce^{-\frac{\delta}{T-t}}
      \leq C(T-t)^{2+\zeta},  \label{Mj-t2+}  \\
    |M_k(t)| \leq& C \|\nabla\Phi\|_{L^\infty}\int_{t}^{T} \frac{D^2(s)}{T-s}ds
       + Ce^{-\frac{\delta}{T-t}}. \label{Mj-D}
\end{align}
Moreover, for any $\ve >0$ as in {\rm Case (I)} and {\rm Case (II)}, we have
that for $t$ close to $T$
\begin{align}
  \bigg|\sum_{k=1}^{K}\frac{M_k(t)}{\lbb_k^2(t)} \bigg|
 \leq&   \frac{ C \ve}{(T-t)^2} \int_{t}^{T}\frac{D^2(s)}{T-s}ds
         + Ce^{-\frac{\delta}{T-t}},   \label{sum-Mj-lbbj2-D}
\end{align}
where   $C, \delta>0$ are independent of $\ve$ and $t$.

Assume additionally that for $t$ close to $T$
\begin{align} \label{gk-w2t}
   \g_k(t) = \omega_k^2(T-t)  + o(T-t),
\end{align}
then we also have
\begin{align} \label{Mj-glbb4-D}
  \bigg|\sum\limits_{k=1}^K \frac{\g_k(t)}{\lbb_k^4(t)} M_k(t) \bigg|
   \leq \frac{C\ve}{(T-t)^3} \int_t^T \frac{D^2(s)}{T-s} ds + Ce^{-\frac{\delta}{T-t}}.
\end{align}
\end{theorem}

\begin{remark}
One advantage of the estimate of localized mass is
that it allows to  control the unstable direction $Q$ in the scalar \eqref{Scal-def} below,
while the (almost) orthogonality in the previous section
controls the remaining five unstable directions.
\end{remark}

\begin{remark}
Estimate \eqref{R-t-0}
and the simple inequality $|M_k| \leq C (\|R\|_{L^2} + \|R\|_{L^2}^2)$
merely give the crude estimate $M_k(t)=o(1)$,
which is insufficient to run the upgradation procedure.

In order to explore higher convergence rate, the keypoint here is that
two more  orders $(T-t)^{2+}$ can be gained by
the local virial identities.
This convergence rate is indeed almost sharp to justify
the integrability in \eqref{Mj-t2+} below.
Moreover,
estimates \eqref{Mj-D}-\eqref{Mj-glbb4-D}
relate the localized mass and the remainder together
and enable us to upgrade their estimates simultaneously
by establishing Gronwall type inequalities.
\end{remark}

{\bf Proof of Theorem \ref{Thm-Loc-Mass}.}
Using the geometrical decomposition \eqref{v-dec} and  Lemma \ref{Lem-inter-est} in Appendix
we have that for $1\leq k\leq K$,
\begin{align*}
\int |v |^2\Phi_kdx
=& \int |U|^2\Phi_kdx +\int |R |^2\Phi_kdx +2{\rm Re}\int \overline{U} R\Phi_k dx  \nonumber \\
=& \int |U |^2\Phi_kdx
   + M_k
   + \calo\(e^{-\frac{\delta}{T-t}}\|R\|_{L^2}\),
\end{align*}
where $\delta>0$.
This yields that
\begin{align} \label{Mj-v-U}
M_k =\int |v|^2\Phi_kdx-\int |U|^2\Phi_kdx+ \calo\(e^{-\frac{\delta}{T-t}}\|R\|_{L^2}\).
\end{align}

Note that,
by \eqref{nav-naS-iint-0+-Thm}, there exists a sequence $\{t_n\}$ to $T$
such that
\begin{align}  \label{wtR-H1-0}
\lim_{t_n\rightarrow T}\|\wt R(t_n)\|_{H^1}=0.
\end{align}
Then,
we split
\begin{align} \label{v-U-Phij-t}
&\int |v(t)|^2\Phi_kdx-\int |U(t)|^2\Phi_kdx \nonumber \\
 =& \(\int |v(t)|^2\Phi_k dx-\int |v(t_n)|^2\Phi_k dx\)
    +\(\int |v(t_n)|^2\Phi_k dx-\int |S(t_n)|^2\Phi_k dx\)   \nonumber \\
  &+\(\int |S(t_n)|^2\Phi_k dx-\int |U(t)|^2\Phi_k dx \),
\end{align}
where $S=\sum_{k=1}^K S_k$.
Note that, by the asymptotic \eqref{v-S-H1-o1-Thm}
and the conservation law of mass,
\begin{align} \label{v-S-t-T}
        \lim_{t_n\rightarrow T} \bigg|\int |v(t_n)|^2\Phi_k dx
     -\int |S(t_n)|^2\Phi_k dx \bigg|
     \leq& C \lim\limits_{t_n\to T} \|v(t_n) - S(t_n)\|_{L^2}
           \( \|v(t_n)\|_{L^2} +  \sum\limits_{k=1}^K \| S_k(t_n)\|_{L^2}\)  \nonumber \\
     \leq& C \lim_{t_n\rightarrow T}\|v(t_n)-S(t_n)\|_{L^2}=0,
\end{align}
and by Lemma \ref{Lem-inter-est} in Appendix below,
\begin{align}
\lim_{t_n\rightarrow T}\int |S(t_n)|^2\Phi_k dx=\|Q\|_{L^2}^2, \quad
\int |U(t)|^2\Phi_k dx=\|Q\|_{L^2}^2+ \calo\(e^{-\frac{\delta}{T-t}}\),
\end{align}
which yield that
\begin{align} \label{S-U-t-T}
   \limsup \limits_{ t_n\to T}
    \bigg|\int |S(t_n)|^2\Phi_k dx - \int |U(t)|^2 \Phi_k dx \bigg|
     = \calo\( e^{-\frac{\delta}{T-t}}\).
\end{align}
Thus, plugging \eqref{v-U-Phij-t}, \eqref{v-S-t-T} and \eqref{S-U-t-T} into \eqref{Mj-v-U}
we obtain
\begin{align} \label{Mj-vt-Phij}
    |M_k(t)|
    \leq \limsup\limits_{t_n\to T} \bigg|\int |v(t)|^2\Phi_kdx-\int |v(t_n)|^2\Phi_k dx\bigg| + C e^{-\frac{\delta}{T-t}}.
\end{align}

To estimate the R.H.S. above,
we use the local virial identities
(see \cite[Lemma 3.6]{M93})
\begin{align}
  \frac{d}{dt}\int |v|^2\Phi_kdx
   =&2{\rm Im}\int \ol{v}\nabla v\cdot  \nabla\Phi_k dx, \label{dt-v2-Phij}  \\
  \frac{d}{dt} ({\rm Im}\int \ol{v}\nabla v\cdot  \nabla\Phi_k dx)
   =&2{\rm Re}\int \nabla^2\Phi_k (\nabla v, \nabla \bar{v}) dx
      -\frac{2}{2+d}\int \Delta \Phi_k |v|^{2+\frac{4}{d}}dx
    -\frac 12 \int\Delta^2\Phi_k|v|^2 dx. \label{dtt-v2-Phij}
\end{align}

Since by the Gagliardo-Nirenberg inequality \eqref{G-N} below
and the conservation law of mass,
\begin{align*}
    \|\wt R\|_{L^{2+\frac 4d}}^{2+\frac 4d}
    \leq C \|\wt R\|_{L^2}^{\frac 4d} \|\na \wt R\|_{L^2}^2
    \leq C (\|v\|_{L^2}^{\frac 4d} + \|S\|_{L^2}^{\frac 4d})\|\na \wt R\|_{L^2}^2
    \leq C \|\na \wt R\|_{L^2}^2.
\end{align*}
Taking into account the exponential decay of $S$ on the support of $\na \Phi_k$
we get from \eqref{dtt-v2-Phij}
\begin{align} \label{dtt-v2-Phij-bdd}
        \bigg|\frac{d}{dt} ({\rm Im}\int \ol{v}\nabla v\cdot  \nabla\Phi_k dx)\bigg|
   \leq& C \int_{ |x-x_k|\geq 4\sigma,1\leq k\leq K}
           |\nabla S+ \na \widetilde{R}|^2 +  |S+\widetilde{R}|^{2+\frac{4}{d}}+ |S+\widetilde{R}|^2 dx \nonumber \\
   \leq& C   \(\|\widetilde{R}\|_{H^1}^2+e^{-\frac{\delta}{T-t}}\).
\end{align}
Similarly,
\begin{align} \label{dtt-v2-Phij-bdd*}
         \lim_{t_n\to T} \bigg|{\rm Im}\int \ol{v}(t_n)\nabla v(t_n)\cdot  \nabla\Phi_k dx \bigg|
   \leq&  \lim_{t_n\to T}\int_{ |x-x_k|\geq 4\sigma,1\leq k\leq K} |S(t_n)+\wt R(t_n)||\na S(t_n) + \na \wt R(t_n)| dx  \nonumber \\
    \leq& C\lim_{t_n\to T} \( \|\widetilde{R}(t_n)\|_{H^1}^2+e^{-\frac{\delta}{T-t_n}}\) =0.
\end{align}
Thus,
integrating \eqref{dtt-v2-Phij} from $t$ to $t_n$
and using \eqref{dtt-v2-Phij-bdd} and
the boundary estimate \eqref{dtt-v2-Phij-bdd*}
we obtain
\begin{align} \label{m1}
   \bigg|{\rm Im}\int \ol{v}(t)\nabla v(t)\cdot  \nabla\Phi_k dx \bigg|
   \leq C \(\int_t^T \|\wt R\|_{H^1}^2 ds + e^{-\frac{\delta}{T-t}}\).
\end{align}
Plugging this into \eqref{dt-v2-Phij}
and then integrating both sides again we obtain
\begin{align}
   \bigg|\int |v(t_n)|^2 \Phi_k dx
   - \int |v(t)|^2 \Phi_k dx \bigg|
   \leq C \(\int_t^T \int_s^T \|\wt R\|_{H^1}^2 drds + e^{-\frac{\delta}{T-t}}\),
\end{align}
which, via \eqref{Mj-vt-Phij}, yields
the first inequality in \eqref{Mj-t2+}.
The second inequality follows from \eqref{nav-naS-iint-0+-Thm}.

Using \eqref{v-dec}, \eqref{Mj-vt-Phij}, \eqref{dt-v2-Phij}
and the integration by parts formula to move the derivative to $U$
yields
\begin{align}
 |M_k(t)|
 \leq& \bigg|\int_{t}^{T} 2{\rm Im}\int \ol{v}\nabla v\cdot  \nabla\Phi_k dx ds\bigg| + C e^{-\frac{\delta}{T-t}} \nonumber \\
 \leq& C\|\nabla\Phi\|_{L^\infty}\int_{t}^{T}\int_{ |x-x_k|\geq 4\sigma,1\leq k\leq K}
       |\na U| |U| + (|\na U|+|U|)|R| + |R||\na R| dxds
    +  Ce^{-\frac{\delta}{T-t}}.
\end{align}
Then, using  \eqref{R-D} and the exponential decay of $Q$ in \eqref{Q-decay} we obtain
\begin{align}
 |M_k(t)|\leq& C\|\nabla\Phi\|_{L^\infty}\int_{t}^{T} \|R\|_{L^2}
    \|\nabla R\|_{L^2}ds
    + C e^{-\frac{\delta}{T-t}} \nonumber \\
 \leq& C \|\nabla\Phi\|_{L^\infty}\int_{t}^{T} \frac{D^2(s)}{T-s}ds
       + C e^{-\frac{\delta}{T-t}},
\end{align}
which yields \eqref{Mj-D}.

Next, we prove estimate \eqref{sum-Mj-lbbj2-D}.
In {\rm Case (II)}, estimate \eqref{sum-Mj-lbbj2-D} follows immediately
from \eqref{lbb-a-t-0} and \eqref{Mj-D}, as $\|\na \Phi\|_{L^\9} \leq C \sigma^{-1} \leq C \ve$.
Below we deal with {\rm Case (I)}.

The observation here is that,
a refined exponential estimate can be gained for the sum of $M_k$, i.e.,
\begin{align} \label{sum-Mj-exp}
  \sum\limits_{k=1}^K M_k(t) = \calo\(e^{-\frac{\delta}{T-t}}\), \ \ for\ t\ close\ to\ T.
\end{align}
Actually,
on one hand, by the conservation law of mass, \eqref{wtR-H1-0} and Lemma \ref{Lem-inter-est} below,
\begin{align} \label{v-L2}
   \|v(t)\|^2_{L^2}
   = \lim\limits_{t_n\to T} \|v(t_n)\|^2_{L^2}
   = \lim\limits_{t_n\to T} \sum\limits_{k=1}^K \|S_k(t_n)\|_{L^2}^2
   = K\|Q\|_{L^2}^2.
\end{align}
On the other hand,
using \eqref{v-dec} to expand the mass $\|v\|_{L^2}^2$ around the blow-up profile $U$
we get
\begin{align} \label{v-L2-expan}
   \|v\|^2_{L^2}
   = \|U\|_{L^2}^2 + 2 {\rm Re} \<U, R\> + \|R\|_{L^2}^2
   = K\|Q\|^2_{L^2} + \sum\limits_{k=1}^K M_k + \calo\(e^{-\frac{\delta}{T-t}}\),
\end{align}
where we also used Lemma \ref{Lem-inter-est}
to decouple different bubbles in the last step.
Thus, combining the two identities \eqref{v-L2} and \eqref{v-L2-expan} we obtain \eqref{sum-Mj-exp}, as claimed.

Then, since
$\lbb_k = \wt \lbb_k \omega_k (T-t)$,
due to  \eqref{calp-wtcalp} and \eqref{wtcalp-wtR-esti} with $L$ replaced by $T-t$,
we see that
\begin{align} \label{lbbk-wt2}
     \bigg|\frac{1}{\lbb_k^2} - \frac{1}{\omega^2(T-t)^2} \bigg|
    =& \bigg|\frac{(\omega^2-\omega_k^2) + (1-\wt \lbb_k^2)\omega_k^2}{\wt \lbb_k^2 \omega_k^2 \omega^2 (T-t)^2}\bigg|
    \leq C \frac{|\omega - \omega_k|+|\wt \lbb_k -1|}{(T-t)^2},
\end{align}
where $\omega$ is as in {\rm Case (I)}.
Taking into account $|\omega_k-\omega|\leq \ve$ in {\rm Case (I)} and
the estimate of $\wt \lbb_k$ in \eqref{wtcalp-wtR-esti}
with $S_L$ replaced by $S$
we get
\begin{align} \label{lbbj2-w2T2}
     \bigg|\frac{1}{\lbb_k^2} - \frac{1}{\omega^2(T-t)^2} \bigg|
     \leq C \frac{\ve + \|v-S\|_{L^2}}{(T-t)^2}.
\end{align}

Thus,
combining \eqref{sum-Mj-exp} and \eqref{lbbj2-w2T2} together
we arrive that
\begin{align} \label{Mj-lbbj2-bdd}
     \bigg| \sum\limits_{k=1}^K \frac{M_k}{\lbb_k^2} \bigg|
   =& \bigg| \sum\limits_{k=1}^K \(\frac{1}{\lbb_k^2} - \frac{1}{\omega^2(T-t)^2}\) M_k
     + \frac{1}{\omega^2(T-t)^2} \sum\limits_{k=1}^K M_k \bigg|  \nonumber \\
   \leq&  C \frac{\ve + \|v-S\|_{L^2}}{(T-t)^2} \sum_{k=1}^K |M_k| + C e^{-\frac{\delta}{T-t}} \nonumber \\
   \leq& C \frac{\ve + \|v-S\|_{L^2}}{(T-t)^2}
          \int_t^T \frac{D^2(s)}{T-s} ds + C e^{-\frac{\delta}{T-t}},
\end{align}
where we also used \eqref{Mj-D} in the last step,
and $C, \delta>0$ are independent of $\ve$.
Thus, taking $t$ close to $T$ such that $ \|v-S\|_{L^2} \leq \ve$
we obtain \eqref{sum-Mj-lbbj2-D}.

It remains to prove \eqref{Mj-glbb4-D}.
In {\rm Case (II)},
using \eqref{lbb-a-t-0}, \eqref{Mj-D} and \eqref{gk-w2t} we have
\begin{align*}
  \bigg| \sum\limits_{k=1}^K \frac{\g_k}{\lbb_k^4} M_k \bigg|
  \leq \frac{C}{(T-t)^3} \sum\limits_{k=1}^K  |M_k|
  \leq  \frac{C\|\na \Phi\|_{L^\9}}{(T-t)^3}
         \sum\limits_{k=1}^K  \int_t^T \frac{D^2(s)}{T-s} ds + Ce^{-\frac{\delta}{T-t}}.
\end{align*}
Then, taking into account  $\|\na \Phi\|_{L^\9} \leq C \sigma^{-1} \leq C \ve$
we obtain \eqref{Mj-glbb4-D}.

Regarding {\rm Case (I)}
we use similar arguments as in the proof of \eqref{sum-Mj-lbbj2-D}.
For simplicity, we set
$\lbb_0 = \omega (T-t)$ and $\g_0 = \omega^2 (T-t)$,
where $\omega$ is as in {\rm Case (I)}.
By \eqref{lbb-a-t-0} and \eqref{gk-w2t},
\begin{align}
   & \lbb_k = \wt \lbb_k \omega_k(T-t), \ \ with \ \ \wt \lbb_k = 1 + o(1), \label{lbbk-wtlbbk-o1} \\
   & \g_k =\wt \g_k + \omega_k^2 (T-t), \ \ with\ \ \wt \g_k = o(T-t). \label{gk-wtgk-ot}
\end{align}
Then, in view of \eqref{sum-Mj-exp}, we infer that
\begin{align} \label{Mj-glbb4-D-esti}
    \bigg| \sum\limits_{k=1}^K \frac{\g_k}{\lbb_k^4} M_k\bigg|
    = \bigg| \sum\limits_{k=1}^K \(\frac{\g_k}{\lbb_k^4} - \frac{\g_0}{\lbb_0^4}\) M_k
        +  \frac{\g_0}{\lbb_0^4}  \sum\limits_{k=1}^K M_k\bigg|
    \leq C \(  \sum\limits_{k=1}^K \bigg|\frac{\g_k}{\lbb_k^4} - \frac{\g_0}{\lbb_0^4}\bigg| |M_k|
           + e^{-\frac{\delta}{T-t}} \).
\end{align}
Note that,
by the straightforward computations,
\begin{align*}
  \bigg|\frac{\g_k}{\lbb_k^4} - \frac{\g_0}{\lbb_0^4}\bigg|
  =& (T-t)^{-3}\(\wt \lbb_k^{-4}\omega_k^{-4} \omega^{-4}\)
    \bigg| \frac{\wt\g_k}{T-t} \omega^4
        + \omega_k^2  \omega^4
     - \omega^2 \wt \lbb_k^4 \omega_k^4 \bigg|,
\end{align*}
and
\begin{align*}
  \omega_k^2  \omega^4
     - \omega^2 \wt \lbb_k^4 \omega_k^4
  =  \omega_k^2 \omega^2
       \((1-\wt \lbb_k)(1+\wt \lbb_k)(1+\wt \lbb^2_k)\omega^2
        + (\omega - \omega_k)(\omega + \omega_k) \wt \lbb_k^4\).
\end{align*}
Then,
using \eqref{lbbk-wtlbbk-o1} and \eqref{gk-wtgk-ot}
and taking $t$ close to $T$ such that
$|\frac{\wt\g_k}{T-t}| + |1-\wt \lbb_k| \leq \ve$
we obtain
\begin{align*}
  \bigg|\frac{\g_k}{\lbb_k^4} - \frac{\g_0}{\lbb_0^4}\bigg|
  \leq \frac{C \ve}{(T-t)^3},
\end{align*}
where $C>0$ is independent of $\ve$.

Therefore, inserting this into \eqref{Mj-glbb4-D-esti}
and using \eqref{Mj-D}
we obtain \eqref{Mj-glbb4-D} and finish the proof.
\hfill $\square$

Theorem \ref{Thm-Energy} below contains the coercivity type control of energy.
\begin{theorem} (Coercivity of energy) \label{Thm-Energy}
There exist $C_1, C_2>0$
such that
for $t$ close to $T$,
\begin{align}  \label{energy-esti}
     & \sum_{k=1}^{K}\frac{|\beta_{k}|^2 }{2\lbb_{k}^2} \|Q\|_{2}^2
        + \sum_{k=1}^{K}  \frac{\g_k^2}{8\lbb_k^2} \|yQ\|_{L^2}^2
        + C_1 \frac{D^2}{(T-t)^2}  \nonumber \\
\leq& E(v)
     + \sum_{k=1}^K \frac{|\beta_k|^2}{2\lbb_k^2} M_k
      + \sum_{k=1}^{K}\frac{M_k}{2\lbb_{k}^2}
      + C_2  \(\sum\limits_{k=1}^K \frac{M_k^2}{(T-t)^2} + e^{-\frac{\delta}{T-t}}\) .
\end{align}
\end{theorem}

\begin{remark}
At this stage,
the exact value of energy $E(v)$ is unclear,
due to the low convergence rate in \eqref{v-S-H1-o1-Thm}.
The precise identification of $E(v)$ will be proved in Subsection \ref{Subsec-Energy-Quant} below.
As we shall see, the energy indeed
admits the quantization into the sum of the energies
of pseudo-conformal blow-up solutions.
The key ingredients for this  fact are
the $H^1$ dispersion of remainder
and the refined estimates of geometrical parameters along a sequence,
which in turn rely on the monotonicity of localized virial functional in Subsection \ref{Subsec-Modif-Loc-Vir} below.
\end{remark}

{\bf Proof.}
The arguments presented below
follow the lines as in the proof of \cite[Lemma 5.14]{SZ20}
and give more precise estimates of errors.

According to \eqref{Mj-def}, for $1\leq k\leq K$,
\ben
\frac{{1}}{\lbb_{k}^2} {\rm Re} \int \ol{U_k} R_k
       + \frac12 |R|^2\Phi_kdx
       -\frac{1}{2\lbb_{k}^2}M_k=0.
\enn
Thus we first reformulate the energy
\begin{align}
E(v)
&= E(v)+\sum_{k=1}^{K}\frac{{1}}{\lbb_{k}^2} {\rm Re} \int \ol{U_k} R_k
       + \frac12 |R|^2\Phi_kdx
       -\sum_{k=1}^{K}\frac{1}{2\lbb_{k}^2}M_k.
\end{align}

Then, expanding the energy $E(v)$ around the main blow-up profile $U$ up to the second order of remainder
we obtain that, similarly to \cite[(5.76)]{SZ20},
\begin{align}
E(v)
=& \sum_{k=1}^K\( \frac{|\beta_{k}|^2}{2\lambda_{k}^2} \|Q\|_{L^2}^2
+ \frac{\gamma_{k}^2}{8\lambda_{k}^2} \|yQ\|_{L^2}^2\) \nonumber \\
 & -  \sum_{k=1}^{K} {\rm Re} \int \({\Delta U_{k}-\frac{1}{\lbb_{k}^2} U_{k}+|U_{k}|^{\frac 4d}U_{k}}\) \overline{R_{k}} dx   \nonumber\\
&+\frac{1}{2}{\rm Re} \int |\nabla R|^2+\sum_{k=1}^{K}\frac{1}{\lbb_{k}^2}|R|^2  \Phi_k
     -(1+\frac2d)|U|^{\frac 4d} |R|^2
     - \frac 2d |U|^{\frac 4d-2}U^2\overline{R}^2 dx  \nonumber\\
&   -\sum_{k=1}^{K}\frac{M_k}{2\lbb_{k}^2}+ \calo\( \sum\limits_{l=3}^{2+\frac 4d} \int |U|^{2+\frac 4d -l} |R|^l dx +e^{-\frac{\delta}{T-t}}\) \nonumber\\
   =:&\sum_{j=0}^{2}E_j-\sum_{k=1}^{K}\frac{M_k}{2\lbb_{k}^2}+\calo\( \sum\limits_{l=3}^{2+\frac 4d} \int |U|^{2+\frac 4d -l} |R|^l dx +e^{-\frac{\delta}{T-t}}\), \label{Ev-expand.0}
\end{align}
where  the exponential decay error arises from the interactions between different bubbles
due to Lemma \ref{Lem-inter-est},
while the remaining errors come  from the remainders of orders higher than two
in the Taylor expansion of nonlinearity.
Note that,
this error can be bounded easily by using \eqref{R-t-0}, \eqref{R-D-Lp}
and the fact that $\|U\|_{L^\9} \leq C(T-t)^{-\frac d2}$:
\begin{align} \label{err-esti-energy}
   \sum\limits_{l=3}^{2+\frac 4d} \int |U|^{2+\frac 4d -l} |R|^l dx
   \leq& C \sum\limits_{l=3}^{2+\frac 4d} (T-t)^{-\frac d2 (2+\frac 4d -l)} \|R\|_{L^l}^l  \nonumber \\
   \leq& C \sum\limits_{l=3}^{2+\frac 4d} (T-t)^{-2} D^l(t)
   \leq C (T-t)^{-2} D^3(t).
\end{align}

For the linear terms of $R$ in $E_1$, we infer from \eqref{equa-Qk}
and the renormalized variable $\ve_k$ in \eqref{Rj-ej} that
\begin{align} \label{Ev-linear-esti.0}
   E_1 = &  -\sum\limits_{k=1}^K \frac{1}{\lambda_{k}^2} {\rm Im} \int(\gamma_{k}\Lambda Q_{k} -2\beta_{k}\cdot\nabla Q_{k}) \ol{\varepsilon_k} dy
   - \sum\limits_{k=1}^K \frac{1}{\lbb_{k}^2} {\rm Re} \int |\beta_{k} - \frac {\gamma_{k}}{ 2} y|^2 Q_{k} \ol{\varepsilon_k} dy.
\end{align}
Note that, by the almost orthogonality in  Lemma \ref{Lem-almost-orth},
\begin{align} \label{Ev-linear-esti.1}
  \frac{1}{\lambda_{k}^2} {\rm Im} \int(\gamma_{k}\Lambda Q_{k} -2\beta_{k}\cdot\nabla Q_{k}) \ol{\varepsilon_k} dx
  = \calo\( e^{-\frac{\delta}{T-t}} \|R\|_{L^2}\).
\end{align}
Moreover, using Lemma \ref{Lem-almost-orth} again
and \eqref{Mj-def} we have
\begin{align} \label{Ev-linear-esti.2}
   \frac{1}{\lbb_{k}^2} {\rm Re} \int |\beta_{k} - \frac {\gamma_{k}}{ 2} y|^2 Q_{k} \ol{\varepsilon_k} dy
   =&  \frac{|\beta_k|^2}{\lbb_k^2} {\rm Re} \int Q_k \ol{\ve_k} dy
     + \calo\(e^{-\frac{\delta}{T-t}} \|R\|_{L^2}\)  \nonumber \\
   =& \frac{|\beta_k|^2}{2\lbb_k^2} M_k
     + \calo\( \frac{|\beta_k|^2}{\lbb_k^2} \|R\|_{L^2}^2 + e^{-\frac{\delta}{T-t}}\|R\|_{L^2}\).
\end{align}
Inserting \eqref{Ev-linear-esti.1} and \eqref{Ev-linear-esti.2} into \eqref{Ev-linear-esti.0} we obtain that
\begin{align}  \label{Ev-linear-esti}
    E_1 = - \sum\limits_{k=1}^K \frac{|\beta_k|^2}{2\lbb_k^2} M_k
     + \calo \(\sum\limits_{k=1}^K \frac{|\beta_k|^2}{\lbb_k^2} \|R\|_{L^2}^2 + e^{-\frac{\delta}{T-t}}\|R\|_{L^2}\).
\end{align}

Regarding the quadratic terms of $R$ in $E_2$,
we claim that for some $C>0$,
\begin{align} \label{Ev-quard-esti}
 E_2
        \geq& C \frac{D^2(t)}{(T-t)^2}
      + \calo \( \sum\limits_{k=1}^K \frac{M_k^2}{(T-t)^2} + e^{-\frac{\delta}{T-t}}\) .
\end{align}

In order to prove \eqref{Ev-quard-esti},
as in \cite[(5.80)]{SZ20},
using the renormalized variable $\wt \ve_k$ defined by
\begin{align} \label{wtvek-def}
     R(t,x) = \lbb_k^{-\frac d2} \wt \ve_k\(t,\frac{x-\a_k}{\lbb_k}\) e^{i\theta_k},
\end{align}
we have
\begin{align} \label{R-quad-decoup-2}
E_2
  =& \sum\limits_{k=1}^K  \frac{1}{\lbb_{k}^2}
     {\rm Re}\int (|\nabla \wt \ve_{k}|^2+|\wt \ve_{k}|^2)\phi_{A}
     -(1+\frac2d) Q^{\frac 4d}|\wt \ve_{k}|^2
     - \frac 2d  Q^{\frac 4d-2}\overline{Q_{k}}^2\wt \ve_{k}^2 dy  \nonumber  \\
&+\sum\limits_{k=1}^K \frac{1}{\lbb_{k}^2} \int (|\nabla \wt \ve_{k}|^2+|\wt \ve_{k}|^2)(\Phi_k(\lbb_{k}y+\alpha_{k})
-\phi_{A}(y))dy
  +\calo\(e^{-\frac{\delta}{T-t}}\|R\|_{L^2}^2\) \nonumber \\
=&: \sum\limits_{k=1}^K E_{21,k} + \sum\limits_{k=1}^K  E_{22,k}
+ \calo\(e^{-\frac{\delta}{T-t}}\|R\|_{L^2}^2\),
\end{align}
where $\phi_A$ is the localized function defined as in Lemma \ref{Lem-coer-f-local} below.

Since by \eqref{Uj-Qj-Q},
\begin{align} \label{Qk-Q-P}
   Q_k = Q+ \calo(P\<y\>^4 Q),
\end{align}
where $P$ is given by \eqref{P-def}
and satisfies $P=o(1)$, due to \eqref{lbb-a-t-0} and \eqref{P-t-0},
we compute
\begin{align}
   E_{21,k}
   \geq& \frac{1}{\lbb^2_k}{\rm Re} \int (|\na \wt \ve_k|^2 + |\wt \ve_k|^2) \phi_A
      -(1+\frac 2d) Q^\frac 4d |\wt \ve_k|^2
      - \frac 2d Q^ \frac 4d \wt \ve_k^2 dy
      +\calo\(\frac{P}{\lbb_k^2}\|R\|_{L^2}^2\),
\end{align}
which along with the local coercivity of linearized operators in Lemma \ref{Lem-coer-f-local}
yields that
\begin{align} \label{E21k-scalvek}
    E_{21,k}
   \geq&  \wt c_k\int \(|\nabla R|^2+\frac{1}{\lbb_{k}^2}|R|^2\)\phi_{A, k} dx
          - C \frac{1}{\lbb^2_k} {\rm Scal}(\wt \ve_k) -  \frac{C\ve D^2}{(T-t)^2},
\end{align}
where $\phi_{A,k}(x):= \phi_A(\frac{x-\a_k}{\lbb_k})$, $\wt c_k>0$,
we also used \eqref{lbbk-t-approx} and
$P(t) \leq \ve$ for $t$ close to $T$ in the last step.

In order to estimate the scalar ${\rm Scal}(\wt \ve_k)$,
we infer from \eqref{Scal-def}, \eqref{wtvek-def} and \eqref{Qk-Q-P}
and the almost orthogonality in Lemma \ref{Lem-almost-orth} that
\begin{align} \label{scal-wtve}
  \frac{1}{\lbb^2_k} {\rm Scal}(\wt \ve_k) = \frac{1}{\lbb_k^2} \<\wt \ve_{k,1}, Q\>^2
   + \calo\(e^{-\frac{\delta}{T-t}}\|R\|^2_{L^2}
   + \frac{P^2}{\lbb_k^2}\|R\|^2_{L^2}\),
\end{align}
where $\wt \ve_{k,1} = {\rm Re} \wt\ve_k$.
Note that, by \eqref{Qk-Q-P},
\begin{align*}
   \<\wt \ve_{k,1}, Q\> =& {\rm Re} \<\wt\ve_k, Q_k\> + \calo(P\|R\|_{L^2})
                        =  {\rm Re} \<R, U_k\> + \calo(P\|R\|_{L^2}) \nonumber \\
                        =& \frac 12 M_k  - \frac 12 \int |R|^2 \Phi_k dx
                            + \calo\(P \|R\|_{L^2} + e^{-\frac{\delta}{T-t}}\|R\|_{L^2}\).
\end{align*}
This yields that
\begin{align}  \label{wtveQ-Mk}
   \<\wt \ve_{k,1}, Q\>^2
   = \calo\(M_k^2  + \|R\|_{L^2}^4 +  P^2\|R\|_{L^2}^2 + e^{-\frac{\delta}{T-t}}\|R\|^2_{L^2}\).
\end{align}
Hence, we may take $t$ even closer to $T$
such that $\|R(t)\|_{L^2}^2 \leq \ve$ to get
\begin{align} \label{wtvek-Mk}
   \frac{1}{\lbb_k^2} \<\wt \ve_{k,1}, Q\>^2
   = \calo\(\frac{M_k^2}{(T-t)^{2}} + \frac{\ve D^2}{(T-t)^2} + e^{-\frac{\delta}{T-t}}\|R\|^2_{L^2}\).
\end{align}
Plugging \eqref{scal-wtve} and \eqref{wtvek-Mk} into \eqref{E21k-scalvek} we obtain
\begin{align} \label{E21k}
   E_{21,k}\geq \wt c_k \int (|\na R|^2 + \frac{|R|^2}{\lbb_k^2}) \phi_{A,k} dx
               + \calo\( \frac{M_k^2}{(T-t)^{2}} + \frac{\ve D^2}{(T-t)^2} + e^{-\frac{\delta}{T-t}}\|R\|^2_{L^2}\).
\end{align}

Regarding $E_{22,k}$, we have from \cite[(5.82)]{SZ20} that
\begin{align} \label{E22k}
 E_{22,k}
 \geq& {\wt c} \int \(|\nabla R|^2+\frac{|R|^2}{\lbb_{k}^2}\)(\Phi_k-\phi_{A,k})dx
    + \calo\( e^{-\frac{\delta}{T-t}} \|R\|_{H^1}^2\),
\end{align}
where $\wt c= \min\{\frac 12, \wt c_k, 1\leq k\leq K\}>0$.

Thus, plugging \eqref{E21k} and \eqref{E22k} into \eqref{R-quad-decoup-2}
and taking $\ve$ small enough we obtain \eqref{Ev-quard-esti}, as claimed.

Therefore,
combining \eqref{Ev-expand.0}, \eqref{err-esti-energy}, \eqref{Ev-linear-esti} and \eqref{Ev-quard-esti} altogether
we conclude that
\begin{align} \label{D-E-Mj}
    &\sum_{k=1}^{K} \( \frac{|\beta_{k}|^2 }{2\lbb_{k}^2} \|Q\|_{2}^2
      +\frac{\g_{k}^2}{8\lbb_{k}^2} \|yQ\|_{L^2}^2 \)
      + C_1 \frac{D^2}{(T-t)^2} \nonumber  \\
\leq& E(v)
      + \sum\limits_{k=1}^K \frac{|\beta_k|^2}{2\lbb_k^2} M_k
      + \sum_{k=1}^{K}\frac{M_k }{2\lbb_{k}^2}  \nonumber \\
    &  + C_2 \(\sum\limits_{k=1}^K \frac{|\beta_k|^2}{\lbb_k^2} \|R\|_{L^2}^2 + \frac{D^3}{(T-t)^2}
       + \sum\limits_{k=1}^K \frac{M_k^2}{(T-t)^2} + e^{-\frac{\delta}{T-t}}\),
\end{align}
where $C_1, C_2>0$.
Taking into account \eqref{P-t-0}, \eqref{lbbk-t-approx}  and \eqref{D-0}
we may take $t$  closer to $T$ such that
\begin{align*}
    C_2 \(\frac{|\beta_k|^2}{\lbb_k^2} \|R\|^2_{L^2} + \frac{D^3}{(T-t)^2} \)
    \leq \frac{C_1}{2} \frac{D^2}{(T-t)^2}.
\end{align*}
Plugging this into \eqref{D-E-Mj}
we obtain \eqref{energy-esti} and finish the proof.
\hfill $\square$

In particular,
because
$E(v)$ is a constant independent of $t$ due to the conservation law of energy,
estimate \eqref{lbb-a-t-0} and Theorems \ref{Thm-Mod}, \ref{Thm-Loc-Mass} and \ref{Thm-Energy} yield the
improved estimates of
remainder and geometrical parameters $\beta$ and $\g$,
which improve the previous ones in \eqref{P-t-0} and \eqref{R-t-0}.

\begin{corollary} (First upgradation of estimates) \label{Cor-beta-gam-D-t-1}
For $t$ close to $T$,
we have
\begin{align}
   & |\b(t)|+|\g(t)|+D(t)\leq C(T-t),  \label{beta-gam-D-t-1}\\
   & Mod(t) \leq C (T-t)^2.  \label{mod-t-2}
\end{align}
In particular,
\begin{align}
   & P(t) \leq C(T-t), \label{P-t-1} \\
   & \|R(t)\|_{L^2} \leq C (T-t),\ \ \|\na R(t)\|_{L^2} \leq C.   \label{R-t-1}
\end{align}
\end{corollary}

\subsection{(Modified) localized virial functionals}  \label{Subsec-Modif-Loc-Vir}

In this subsection,
two types of localized virial functionals adapted to the multi-bubble case
are introduced.
The main results are Theorems \ref{Thm-Loc-virial} and \ref{Thm-Modi-Loc-virial} below,
containing the key monotonicity properties of these functionals.

We first analyze the localized virial functional defined by \eqref{I-def}
to obtain the refined space time estimate of remainder.
This allows to obtain the precise asymptotic order and positivity of parameter $\g$.
More importantly,
it reveals the $H^1$ dispersion of remainder along a sequence,
which is the key towards the derivation of energy quantization in Subsection \ref{Subsec-Energy-Quant} below.

Then, we prove the monotonicity of another modified localized virial functional
defined by \eqref{wtI-def}.
This enables us to cancel the bad $\calo(1)$ terms in the energy estimate
and thus leads to an improved estimate \eqref{D-iter-2}
in Subsection \ref{Subsec-D-2} below.

Let us first consider the localized virial functional.
Let $\chi(x)=\psi(|x|)$ be a smooth radial function on $\R^d$,
where $\psi$ satisfies
$\psi'(r) = r$ if $r\leq 1$,
$\psi'(r) = 2- e^{-r}$ if $r\geq2$,
and
\be\label{chi}
 \bigg|\frac{\psi^{'''}(r)}{\psi^{''}(r)} \bigg|\leq C,
\ \ \frac{\psi'(r)}{r}-\psi^{''}(r) \geq0.
\ee
Let $\chi_A(x) :=A^2\chi(\frac{x}{A})$, $A>0$.
The localized virial functional is defined by
\begin{align} \label{I-def}
   \mathcal{L}:=\sum_{k=1}^K \half {\rm Im}
            \int (\nabla\chi_A) \(\frac{x-\alpha_{k}}{\lambda_k}\)\cdot\nabla R \ol{R} \Phi_kdx
            - \sum_{k=1}^K \frac{\g_k}{4\lbb_k}\|xQ\|_{L^2}^2.
\end{align}

\begin{remark}
The localization functions $\{\Phi_K\}$ are introduced  in \eqref{I-def} in this particular way
mainly to ensure the key monotonicity property in \eqref{dt-I} below.
\end{remark}

\begin{theorem}  \label{Thm-Loc-virial}
(Monotonicity of localized virial functional)
We have that for $t$ close to $T$,
\begin{align}   \label{dt-I}
\frac{d\mathcal{L}}{dt}
   \geq   C \sum_{k=1}^{K} \frac{1}{\lambda_k^3} \int |\nabla \varepsilon_k|^2 e^{-\frac{|y|}{A}}+|\varepsilon_k|^2 dy
          -\sum_{k=1}^{K}\frac{ M_k}{\lbb_k^3}
          + \calo(Er),
\end{align}
where $C>0$ is independent of $t$,
$\ve_k$ is the renormalized remainder given by \eqref{Rj-ej} and
\begin{align} \label{Error}
   Er(t)= (T-t) + (T-t)^{-2}\(D^2(t)+ \sum\limits_{k=1}^K |M_k(t)|\)
               + e^{-\frac{\delta}{T-t}}.
\end{align}
\end{theorem}

{\bf Proof.}
The Morawetz type functional on the R.H.S. of \eqref{I-def}
can be estimated similarly as in the proof of \cite[Lemma 5.12]{SZ20}.

To be precise,
straightforward computations show that
\begin{align} \label{dt-calL.1}
    & \frac{d}{dt} \bigg(\sum_{k=1}^K \half {\rm Im}
            \int \nabla\chi_A (\frac{x-\alpha_{k}}{\lambda_k})\cdot\nabla R\ol{R}\Phi_kdx \bigg) \nonumber \\
   =&\sum_{k=1}^K\frac{1}{2}{\rm Im}
       \bigg< \partial_t \(\nabla\chi_A\(\frac{x-\alpha_{k}}{\lambda_{k}}\)\)\cdot\nabla R, {R}_k \bigg>
     +\sum_{k=1}^K \frac{1}{2\lambda_{k}}   {\rm Im}
       \bigg< \Delta\chi_A \(\frac{x-\alpha_{k}}{\lambda_{k}}\)R_k, \pa_t {R} \bigg>  \nonumber \\
   &+   \sum_{k=1}^K \frac{1}{2}    {\rm Im}
       \bigg< \nabla\chi_A\(\frac{x-\alpha_{k}}{\lambda_{k}}\)\cdot ( \nabla R_k + \na R \Phi_k), \partial_t  R \bigg> \nonumber \\
   =:& \sum\limits_{k=1}^K \call_{t,k1} + \call_{t,k2} + \call_{t,k3}.
\end{align}

Using direct computations, \eqref{mod-t-2} and \eqref{P-t-1}
we have that, similarly to \cite[(5.49)]{SZ20},
\begin{align} \label{calL-k1}
   \call_{t,k1}
   \leq C  \lbb_k^{-2} (Mod_k + P) \|\na R\|_{L^2} \|R\|_{L^2}
   \leq C \frac{D^2}{(T-t)^2}.
\end{align}

Moreover, in order to treat $\call_{t,k2}$,
we infer from equation \eqref{equa-NLS} that
\begin{align} \label{equa-R}
   i \partial_t R + \Delta R + |v|^{\frac 4d} v- |U|^{\frac 4d} U = -\eta,
\end{align}
where
\begin{align} \label{eta-def}
   \eta = i\partial_t U + \Delta U + |U|^{\frac 4d} U.
\end{align}
Note that,
by equation \eqref{equa-Ut} and Lemma \ref{Lem-inter-est},
\begin{align} \label{eta-Mod-esti}
   \|\eta\|_{L^2} \leq C\(\frac{Mod}{(T-t)^2} + e^{-\frac{\delta}{T-t}}\), \ \
  \|\na \eta\|_{L^2} \leq C\(\frac{Mod}{(T-t)^3} + e^{-\frac{\delta}{T-t}}\),
\end{align}
Then, in view of the expansion \eqref{f-expan}, equation \eqref{equa-R} yields that
\begin{align}  \label{call-k2.0}
   \call_{t,k2} = -\frac{1}{2\lbb_k} {\rm Re}
                  \bigg<\Delta \chi_A\(\frac{x-\a_k}{\lbb_k}\) R_k,
                   \Delta R + f'(U) \cdot R + f''(U,R)\cdot R^2 + \eta\bigg>,
\end{align}
where $f'(U)\cdot R$ and $f''(U,R)\cdot R^2$ are given by \eqref{f'UR} and \eqref{f''UR-R2} below, respectively.
Note that, similarly to \cite[(5.56)]{SZ20}, we have
\begin{align}   \label{call-k2.1}
   & -\frac{1}{2\lbb_k} {\rm Re}
                  \bigg<\Delta \chi_A\(\frac{x-\a_k}{\lbb_k}\) R_k, \Delta R\bigg>   \nonumber \\
   =& - \frac{1}{4\lbb_k^3} {\rm Re} \int \Delta^2  \chi_A \(\frac{x-\a_k}{\lbb_k}\) |R_k|^2 dx
     + \frac{1}{2\lbb_k} {\rm Re} \int \Delta \chi_A \(\frac{x-\a_k}{\lbb_k}\) |\na R_k|^2 dx \nonumber \\
     & + \calo\(\|R\|_{H^1}^2 + \|\na R\|_{L^2}\|R\|_{L^2}\).
\end{align}
We also infer from Lemma \ref{Lem-inter-est} that
\begin{align}   \label{call-k2.2}
   {\rm Re} \bigg<\Delta \chi_A\(\frac{x-\a_k}{\lbb_k}\) R_k,  f'(U) \cdot R\bigg>
   =   {\rm Re} \bigg<\Delta \chi_A\(\frac{x-\a_k}{\lbb_k}\) R_k,  f'(U_k) \cdot R_k \bigg>
      + \calo(e^{-\frac{\delta}{T-t}}).
\end{align}
Moreover, since by \eqref{f''UR-R2},
\begin{align} \label{f''UR-R2-esti}
   |f''(U, R)\cdot R^2| \leq C (|U|^{\frac 4d -1} + |R|^{\frac 4d -1}) |R|^2,
\end{align}
taking into account
$\|U(t)\|_{L^\9} \leq C(T-t)^{-\frac d2}$,
\eqref{R-D-Lp}, \eqref{beta-gam-D-t-1} and \eqref{eta-Mod-esti} we get
\begin{align}  \label{call-k2.3}
       & \bigg|\frac{1}{2\lbb_k} {\rm Re} \bigg<\Delta \chi_A\(\frac{x-\a_k}{\lbb_k}\) R_k,  f''(U,R)\cdot R^2 + \eta\bigg> \bigg| \nonumber \\
   \leq& C (T-t)^{-1}\( (T-t)^{-2+\frac d2} \|R\|_{L^3}^3 + \|R\|_{L^{2+\frac 4d}}^{2+\frac 4d} + \|\eta\|_{L^2}\|R\|_{L^2}\)  \nonumber \\
     \leq& C \(\frac{D^2}{(T-t)^2} + \frac{Mod \|R\|_{L^2}}{(T-t)^3} + e^{-\frac{\delta}{T-t}}\).
\end{align}
Thus, plugging \eqref{call-k2.1}, \eqref{call-k2.2} and \eqref{call-k2.3} into  \eqref{call-k2.0}
we obtain
\begin{align} \label{call-k2}
   \call_{t,k2}
   =& -\frac{1}{4\lambda_k^3}{\rm Re}\int \Delta^2\chi_A \(\frac{x-\alpha_k}{\lambda_k} \)|R_k|^2 dx
      +\frac{1}{2\lambda_k}{\rm Re}\int \Delta\chi_A\(\frac{x-\alpha_k}{\lambda_k}\)|\nabla R_k|^2 dx  \nonumber \\
   &-\frac{1}{2\lambda_{k}}
   {\rm Re} \bigg< \Delta\chi_A\(\frac{x-\alpha_{k}}{\lambda_{k}}\)R_k, f'(U_k)\cdot R_k \bigg>
   + \calo\(\frac{D^2}{(T-t)^2} + \frac{Mod \|R\|_{L^2}}{(T-t)^3} + e^{-\frac{\delta}{T-t}}\).
\end{align}

Regarding $\call_{t,k3}$ we use equation \eqref{equa-R} and \eqref{f-expan} again to get
\begin{align} \label{call-k4.0}
   \call_{t,k3} =& -\frac{1}{2} {\rm Re}
                  \bigg<\na \chi_A\(\frac{x-\a_k}{\lbb_k}\) (\na R_k+ \na R \Phi_k),
                   \Delta R + f'(U) \cdot R + f''(U,R)\cdot R^2 + \eta\bigg>.
\end{align}
Similarly to \cite[(5.62)]{SZ20},
we have
\begin{align} \label{call-k4.1}
   & - {\rm Re} \bigg< \nabla\chi_A\(\frac{x-\alpha_{k}}{\lambda_{k}}\)\cdot
        (\nabla R_k + \nabla R \Phi_k ),  \Delta {R} \bigg>   \nonumber  \\
 =&  {\rm Re} \int \frac{2}{\lbb_k} \na^2 \chi_A \(\frac{x-\alpha_{k}}{\lambda_{k}}\) (\na R_k, \na \ol{R_k})
     - \frac{1}{\lbb_k} \Delta \chi_A \(\frac{x-\alpha_{k}}{\lambda_{k}}\) |\na R_k|^2 dx  + \calo\(\frac{D^2}{(T-t)^2}\).
\end{align}
We also have,
via Lemma \ref{Lem-inter-est},
\begin{align} \label{call-k4.2}
    & \frac{1}{2} {\rm Re}  \bigg<\na \chi_A\(\frac{x-\a_k}{\lbb_k}\) (\na R_k+ \na R \Phi_k), f'(U) \cdot R \bigg> \nonumber \\
   =& {\rm Re}  \bigg<\na \chi_A\(\frac{x-\a_k}{\lbb_k}\) \cdot \na R_k, f'(U_k) \cdot R_k \bigg>
      + \calo\(e^{-\frac{\delta}{T-t}}\).
\end{align}
Moreover, using \eqref{R-D-Lp}, \eqref{eta-Mod-esti}, \eqref{f''UR-R2-esti}
and integration by parts formula we get
\begin{align} \label{call-k4.3}
   & \bigg|\frac{1}{2} {\rm Re}  \bigg<\na \chi_A\(\frac{x-\a_k}{\lbb_k}\) (\na R_k+ \na R \Phi_k), f''(U,R) \cdot R^2 + \eta \bigg> \bigg| \nonumber \\
   \leq& C \int (|\na R|+|R|)(|U|^{\frac 4d-1} + |R|^{\frac 4d-1}) |R|^2 dx
        + C (\lbb_k^{-1} \|R\|_{L^2} \|\eta\|_{L^2} + \|R\|_{L^2} \|\na \eta\|_{L^2}) \nonumber \\
   \leq& C \((T-t)^{-2+\frac d2} \|R\|_{H^1} \|R\|_{L^4}^2
            + \|R\|_{H^1} \|R\|_{L^{2+\frac 8d}}^{1+\frac 4d}
            +  \frac{Mod\|R\|_{L^2}}{(T-t)^3}   + e^{-\frac{\delta}{T-t}} \) \nonumber \\
   \leq&  C\bigg(\frac{D^2}{(T-t)^2} + \frac{Mod\|R\|_{L^2}}{(T-t)^3}+ e^{-\frac{\delta}{T-t}} \bigg).
\end{align}
Thus, estimates \eqref{call-k4.0}-\eqref{call-k4.3} together yield that
\begin{align}   \label{calL-k3}
    \call_{t,k3} =&\frac{1}{\lambda_k}{\rm Re}\int \nabla^2\chi_A\(\frac{x-\alpha_k}{\lambda_k}\)(\nabla R_k,\nabla \ol{R_k}) dx
    -\frac{1}{2\lambda_k}{\rm Re}\int \Delta\chi_A\(\frac{x-\alpha_k}{\lambda_k}\)|\nabla R_k|^2 dx \nonumber  \\
   &- \bigg< \nabla\chi_A\(\frac{x-\alpha_{k}}{\lambda_{k}}\)\cdot \nabla R_k,f^\prime(U_k)\cdot R_k \bigg>
     + \calo\(\frac{D^2}{(T-t)^2} + \frac{Mod \|R\|_{L^2}}{(T-t)^3} +e^{-\frac{\delta}{T-t}}\).
\end{align}

Therefore, combining \eqref{calL-k1}, \eqref{call-k2} and \eqref{calL-k3} altogether
we lead to
\begin{align}\label{lv2}
 &\frac{d}{dt} \bigg(\sum_{k=1}^K \half {\rm Im}
            \int \nabla\chi_A \(\frac{x-\alpha_{k}}{\lambda_k}\)\cdot\nabla R\ol{R}\Phi_kdx \bigg) \nonumber \\
 =& \sum\limits_{k=1}^K \frac{1}{\lbb_k} {\rm Re} \int \na^2 \chi_A\(\frac{x-\alpha_{k}}{\lambda_k}\) (\na R_k, \na \ol{R_k}) dx
    -  \sum\limits_{k=1}^K \frac{1}{4\lbb_k^3} {\rm Re} \int \Delta^2 \chi_A \(\frac{x-\a_k}{\lbb_k}\) |R_k|^2 dx \nonumber \\
  & -  \sum\limits_{k=1}^K {\rm Re} \bigg\< \frac{1}{2\lbb_k} \Delta \chi_A\(\frac{x-\a_k}{\lbb_k}\) R_k
       + \na \chi_A\(\frac{x-\a_k}{\lbb_k}\)\cdot \na R_k, f'(U_k)\cdot R_k \bigg\> \nonumber \\
  &+ \calo\(\frac{D^2}{(T-t)^2} + \frac{Mod \|R\|_{L^2}}{(T-t)^3} +e^{-\frac{\delta}{T-t}} \).
\end{align}
By the integration by parts formula,
the third term on the R.H.S. above equals to
\begin{align*}
   \sum\limits_{k=1}^K {\rm Re} \int \na \chi_A(\frac{x-\a_k}{\lbb_k}) \cdot \na \ol{U_k}
       \ (f''(U_k)\cdot R_k^2) dx,
\end{align*}
where $f''(U_k)\cdot R_k^2$ is as in \eqref{f''Qk-vek2} with $U_k$, $R_k$ replacing $Q_k$ and $\ve_k$, respectively.
Thus, using the renormalized remainder $\ve_k$ given by \eqref{Rj-ej}
we arrive at
\begin{align}
&\frac{d}{dt} \bigg(\sum_{k=1}^K \half {\rm Im}
            \int \nabla\chi_A \(\frac{x-\alpha_{k}}{\lambda_k}\)\cdot\nabla R\ol{R}\Phi_kdx \bigg) \nonumber \\
  = &   \sum_{k=1}^{K} \frac{1}{\lambda^3_k} {\rm Re} \int \nabla^2\chi_A(y)\(\nabla \ve_k,\nabla \ol{\ve_k}\) dy
           - \sum_{k=1}^{K} \frac{1}{4\lambda_k^3}\int \Delta^2\chi_A(y)|\ve_k|^2 dy
            \nonumber \\
 &+ \sum_{k=1}^{K}\frac{1}{\lbb_k^3} {\rm Re} \int
   (\nabla\chi_A (y )\cdot\nabla \ol{Q_k}) (f''(Q_k) \cdot \ve_k^2)  dy
     + \calo\(\frac{D^2(t)}{(T-t)^2}
      + \frac{Mod\|R\|_{L^2}}{(T-t)^3}
      +e^{-\frac{\delta}{T-t}}\).
\end{align}
This gives the control of Morawetz type functional in \eqref{I-def}.

Regarding the second term on the R.H.S. of \eqref{I-def},
using \eqref{beta-gam-D-t-1} we compute
\begin{align}   \label{gam-lbb-dt}
   \frac{d}{dt} \(\frac{\g_k}{\lbb_k}\)
=\frac{\lbb_k^2\dot{\g_k}+\g_k^2}{\lbb_k^3}-\frac{\g_k}{\lbb_k}\frac{\lbb_k\dot{\lbb}_k+\g_k}{\lbb_k^2}
=\frac{\lbb_k^2\dot{\g_k}+\g_k^2}{\lbb_k^3}+\calo\(\frac{Mod}{\lbb_k^2}\),
\end{align}
which along with \eqref{lbb2g+g2-f''Qve2}, \eqref{P-t-1} and \eqref{R-t-1} yields
\begin{align}  \label{dt-I.2}
     \frac{d}{dt} \(\sum_{k=1}^K \frac{\g_k}{4\lbb_k}\|xQ\|_2^2 \)
  =& \sum_{k=1}^{K}\frac{1}{\lbb_k^3}M_k(t) - \sum_{k=1}^{K}\frac{1}{\lbb_k^{3}}\int |R(t)|^2\Phi_kdx  \nonumber \\
   &  + \sum_{k=1}^{K} \frac{1}{\lbb_k^{3}}{\rm Re}
      \int   (1+\frac 2d)|Q_k|^{\frac {4}{d}}|\varepsilon_k|^2
             + \frac 2d |Q_k|^{\frac 4d-2} Q_k^2 \ol{\ve_k}^2\ dy     \nonumber \\
    & + \sum_{k=1}^{K} \frac{1}{\lbb_k^{3}}{\rm Re}\int
   (y \cdot\nabla \ol{Q_k})\ (f''(Q_k)\cdot \ve_k^2) \ dy \nonumber \\
    & +(T-t)^{-3}\calo\((T-t)Mod + (T-t)^2\|R\|_{L^2}   + D^3
               + e^{-\frac{\delta}{T-t}}\).
\end{align}

Thus,
combining \eqref{lv2} and \eqref{dt-I.2} altogether
and then using the inequality
\begin{align} \label{R2phik-vek2}
    \int |R|^2 \Phi_k dx
    \geq \int |R|^2 \Phi^2_k dx
    = \int |\ve_k|^2 dy
\end{align}
we arrive at, if $\ve_k= \ve_{k,1} + i \ve_{k,2}$,
\begin{align} \label{dt-local-virial}
\frac{d\mathcal{L}}{dt}
   \geq&
   \sum\limits_{k=1}^K
   \frac{1}{\lambda_k^3}
   \bigg(\int \nabla^2\chi \(\frac{y}{A} \)(\nabla \varepsilon_k,\nabla \ol{\varepsilon_k}) dy
          + \int |\varepsilon_k|^2dy
          -\int (1+\frac 4d)Q^\frac{4}{d}\varepsilon_{k,1}^2+Q^\frac{4}{d}\varepsilon_{k,2}^2  dy \nonumber \\
&  -\frac{1}{4A^2}\int \Delta^2\chi (\frac{y}{A} )|\varepsilon_k|^2 dy\bigg)
 + \sum\limits_{k=1}^K  \frac 2d  \frac{1}{\lambda_k^3}\int (A\nabla\chi (\frac{y}{A} )-y ) \cdot \nabla Q
Q^{\frac 4d -1} \((1+\frac 4d)\varepsilon_{k,1}^2+ \varepsilon_{k,2}^2\)dy \nonumber \\
& -\sum_{k=1}^{K}\frac{M_k }{\lbb_k^3}
  + \calo(\wt{Er}),
\end{align}
where the error term
\begin{align}  \label{wtError}
   \wt{Er}
   :=(T-t)^{-3}\( (T-t)Mod
               + (T-t)^2 \|R\|_{L^2}
               + (T-t)D^2
               + D^3
               + e^{-\frac{\delta}{T-t}}\).
\end{align}

Since
\begin{align}
\int \na^2 \chi\(\frac yA\)(\na \ve_k, \na \ol{\ve_k}) dy
\geq \int \psi''\(\bigg|\frac yA\bigg|\) |\na \ve_k|^2 dy,
\end{align}
applying Lemma \ref{Lem-coer-f-local} with $\phi(x) :=\psi''(|x|)$
we lead to
\begin{align*}
   \frac{d\call}{dt}
   \geq& \sum\limits_{k=1}^K \frac{1}{\lbb_k^3}
   \( C  \int\psi''(\frac yA) |\na \ve_k|^2 +  |\ve_k|^2 dy - \frac{1}{4A^2} \int \Delta^2 \chi(\frac yA)|\ve_k|^2 dy \) \nonumber \\
   & + \sum\limits_{k=1}^K  \frac 2d  \frac{1}{\lambda_k^3}\int (A\nabla\chi (\frac{y}{A} )-y ) \cdot \nabla Q
Q^{\frac 4d -1} \((1+\frac 4d)\varepsilon_{k,1}^2+ \varepsilon_{k,2}^2\)dy  \nonumber \\
   & -\sum\limits_{k=1}^K \frac {M_k }{\lbb_k^3}
     + \calo\(\sum\limits_{k=1}^K \frac{1}{\lbb_k^3} {\rm Scal}(\ve_k) + \wt{Er}\).
\end{align*}
Taking into account for $A$ large enough
\begin{align*}
    \frac{1}{4A^2} \bigg|\Delta^2\chi \bigg(\frac{y}{A} \bigg)\bigg| \leq \frac {1}{4}  C  \psi'' \bigg(\bigg|\frac{y}{A}\bigg| \bigg),  \ \
     \frac 2d(2+\frac 4d) \bigg|A\nabla\chi (\frac{y}{A} )-y \bigg|\bigg|\na QQ^{\frac 4d -1}\bigg| \leq \frac {1}{4} C \psi'' \bigg(\bigg|\frac{y}{A}\bigg| \bigg),
\end{align*}
we arrive at
\begin{align}\label{dtI-esti}
\frac{d\call}{dt}
   \geq&  C \sum_{k=1}^{K} \frac{1}{\lambda_k^3} \int |\nabla \varepsilon_k|^2 e^{-\frac{|y|}{A}}+|\varepsilon_k|^2 dy
          -\sum_{k=1}^{K}\frac{M_k }{\lbb_k^3}  +  \calo\(\sum\limits_{k=1}^K \frac{1}{\lbb_k^3} {\rm Scal}(\ve_k) + \wt{Er}\).
\end{align}

Note that,
by estimates \eqref{Mod-bdd} and \eqref{R-t-1},
\begin{align}
(T-t)Mod\leq C\((T-t)^3D+(T-t)D^2 + (T-t)\sum_{k=1}^{K}|M_k|+e^{-\frac{\delta}{T-t}}\).
\end{align}
Moreover,
by Cauchy's inequality,
for any $\ve>0$,
\begin{align} \label{tRL2-esti}
(T-t)^2\|R\|_{L^2}\leq \varepsilon \|R\|^2_{2}+ \ve^{-1} (T-t)^4.
\end{align}
Hence, we infer from \eqref{lbb-a-t-0} and \eqref{wtError}-\eqref{tRL2-esti} that
\begin{align} \label{wtEr-Er}
   \wt{Er} \leq C \(\ve (T-t)^{-3} \|R\|_{L^2}^2 + Er\)
           \leq C \(\ve \sum\limits_{k=1}^K \frac{1}{\lbb_k^3} \int |\ve_k|^2 dy +  Er\),
\end{align}
where $Er$ is given by \eqref{Error}.

Moreover, similarly to \eqref{scal-wtve} and \eqref{wtveQ-Mk},
we have
\begin{align} \label{scalvek-esti}
   {\rm Scal}(\ve_k) = \calo(M_k^2 + \|R\|_{L^2}^4+ P^2 \|R\|_{L^2}^2 + e^{-\frac{\delta}{T-t}}\|R\|_{L^2}^2).
\end{align}
This yields that
\begin{align} \label{lbb3scalvek}
    \sum_{k=1}^{K} \frac{1}{\lbb_k^3} {\rm Scal}(\ve_k)
    \leq C (T-t)^{-2} \(D^2 + \sum\limits_{k=1}^K |M_k|\) + Ce^{-\frac{\delta}{T-t}}
    =\calo(E_r).
\end{align}

Therefore, plugging \eqref{wtEr-Er} and \eqref{lbb3scalvek}  into \eqref{dtI-esti}
and then taking $\ve$ sufficiently small
we obtain \eqref{dt-I}.
The proof  of Theorem \ref{Thm-Loc-virial} is complete.
\hfill $\square$

By virtue of
Theorem \ref{Thm-Loc-virial}, we derive the refined space time estimate of remainder and
the precise  asymptotic order of  parameter $\g$,
which improve the previous estimates in \eqref{R-t-1} and  \eqref{P-t-0} respectively
and, more importantly,
allow to obtain the energy quantization in the next section.

\begin{theorem} (Refined estimates of $R$ and $\g$) \label{Thm-gamj-wj2-t}
There exists $C>0$ such that for $t$ close to $T$,
\begin{align}  \label{Rj-lbbj3-bdd}
   \sum\limits_{k=1}^K \int_{t}^{T}\frac{\|R_k\|_{L^2}^2}{\lbb_k^3}ds\leq C.
\end{align}
Moreover, for $1\leq k\leq K$,
\begin{align} \label{gamj-wj-t}
\g_k (t) -\omega_k^2(T-t)=o(T-t),\ \ as\ t\ close\ to\ T.
\end{align}
\end{theorem}

\begin{remark} \label{Rem-R-integ}
The space time estimate \eqref{Rj-lbbj3-bdd} is important to derive
the precise asymptotic behavior of $\g_k$ in \eqref{gamj-wj-t}
and to control the remainder near the singularities (see Lemma \ref{Lem-R-H1-In} below).
Note that, the previous estimate \eqref{beta-gam-D-t-1} is insufficient to get \eqref{Rj-lbbj3-bdd}.
Hence, roughly speaking,  \eqref{Rj-lbbj3-bdd} upgrades the convergence rate of $\|R\|_{L^2}$
from  $T-t$ to $(T-t)^{1+}$.

\end{remark}

\begin{remark} \label{Rem-g-asymp}
Compared with  \eqref{P-t-0},
\eqref{gamj-wj-t} gives the precise leading order of the asymptotic behavior of $\g_k$.
This enables to identify the exact value of energy in the next section.
It also yields the positivity of $\g_k$
for $t$ close to $T$,
which is important to derive the monotonicity of
modified localized virial functional in
Theorem \ref{Thm-Modi-Loc-virial} below.
\end{remark}

{\bf Proof of Theorem \ref{Thm-gamj-wj2-t}.}
Integrating \eqref{dt-I} from $t$ to $\tilde t$ we have
\begin{align}  \label{Rj-lbbj3-bdd.0}
\sum_{k=1}^{K}\int_{t}^{\tilde t}\frac{\|R_k\|_{L^2}^2}{\lbb_k^3}ds
  \leq & C\(|\mathcal{L}(\tilde  t)-\mathcal{L}(t)|
        +\sum_{k=1}^{K}\int_{t}^{\tilde t}\frac{|M_k|}{\lbb_k^3} ds
         +\int_{t}^{\tilde t}
            Er\ ds\),
\end{align}
where the error term $Er$ is given by \eqref{Error}.

Note that,
by estimates \eqref{lbbk-t-approx} and \eqref{beta-gam-D-t-1},
\begin{align} \label{I-t}
|\mathcal{L}(t)|\leq C\(\|\na R\|_{L^2}\|R\|_{L^2}+ \sum\limits_{k=1}^K \bigg|\frac{\g_k}{\lbb_k}\bigg|\)\leq C<\9.
\end{align}
Moreover, by \eqref{Mj-t2+},
\begin{align} \label{Mklbb3-integ}
\sum_{k=1}^{K}\int_{t}^{T}\frac{ |M_k| }{\lbb_k^3} ds
     \leq&C (T-t)^{\zeta} ,
\end{align}
and by \eqref{Mj-t2+} and \eqref{beta-gam-D-t-1},
\begin{align} \label{Error-int-esti}
\int_{t}^{T} Er\ ds\leq C(T-t).
\end{align}

Thus,
plugging estimates \eqref{I-t}-\eqref{Error-int-esti} into \eqref{Rj-lbbj3-bdd.0}
and letting $\tilde t\to T$ we obtain \eqref{Rj-lbbj3-bdd}.

Next we prove \eqref{gamj-wj-t}.
Using \eqref{gam-lbb-dt} we get
\begin{align}
   \int_t^T \bigg|\frac{d}{ds}\(\frac{\g_k}{\lbb_k}\)\bigg|ds
   \leq \int_t^T \frac{|\lbb_k^2 \dot \g_k + \g_k^2|}{\lbb_k^3} + C\frac{Mod}{\lbb_k^2} ds.
\end{align}
Then, using \eqref{lbb2g+g2-Mod.0}, \eqref{lbb-a-t-0}, \eqref{Mj-t2+},
Corollary \ref{Cor-beta-gam-D-t-1} and \eqref{Rj-lbbj3-bdd}
we obtain
\begin{align*}
    \int_t^T \bigg|\frac{d}{ds}\(\frac{\g_k}{\lbb_k}\)\bigg|ds
    \leq& C \int_t^T \frac{\|R\|_{L^2}^2 + |M_k| + (T-s)^3}{(T-s)^3} ds
          + C \int_t^T \frac{Mod}{\lbb_k^2} ds \nonumber \\
    \leq& C \(\int_t^T \frac{\|R\|_{L^2}^2 + (T-s)^{2+\zeta}}{(T-s)^3} ds +T-t\)
    \leq C <\9.
\end{align*}
In particular, this yields that for some $c_k \in \bbr$,
\begin{align}  \label{gamj-lbbj-cj}
\lim_{t\rightarrow T}\frac{\g_k(t)}{\lbb_k(t)}=c_k.
\end{align}

Below, we claim that
\begin{align} \label{cj-wj}
    c_k = \omega_k, \ \ 1\leq k\leq K.
\end{align}
Hence,
\eqref{gamj-lbbj-cj} and \eqref{cj-wj} together give that
\begin{align} \label{glbb-omega-o1}
\frac{\g_k(t)}{\lbb_k(t)}-\omega_k=o(1), \ \ for\ t\ close\ to\ T,
\end{align}
which along with \eqref{lbb-a-t-0} yields the desirable estimate \eqref{gamj-wj-t}.

It remains to prove \eqref{cj-wj}.
For this purpose,
we see that
\begin{align*}
\dot{\lbb_k}=\frac{\lbb_k\dot{\lbb_k}+\g_k}{\lbb_k}-\frac{\g_k}{\lbb_k}
= \calo \(\frac{Mod}{\lbb_k}\)
-\frac{\g_k}{\lbb_k}.
\end{align*}
Since $\lim_{t\to T} \lbb_k(t) = 0$,
we get
\begin{align} \label{lbb-integ-esti.0}
   \lbb_k(t) = \int_t^T \frac{\g_k}{\lbb_k} ds
               + \calo\(\int_t^T \frac{Mod}{\lbb_k} ds\).
\end{align}
Note that, by \eqref{gamj-lbbj-cj},
\begin{align} \label{lbb-integ-esti.1}
   \int_t^T \frac{\g_k}{\lbb_k} ds
   = c_k(T-t) + o(T-t).
\end{align}
Moreover, by \eqref{mod-t-2},
\begin{align} \label{lbb-integ-esti.2}
  \int_{t}^{T}\frac{Mod}{\lbb_k}ds
  \leq C \int_t^T (T-s)ds
  \leq C (T-t)^2
  =o(T-t).
\end{align}
Thus, combining \eqref{lbb-integ-esti.0}-\eqref{lbb-integ-esti.2}
together we conclude that
\begin{align}
    \lbb_k(t) = c_k(T-t)+o (T-t),
\end{align}
which along with the estimate of $\lbb_k$ in \eqref{lbb-a-t-0} yields \eqref{cj-wj},
thereby finishing the proof.
\hfill $\square$

We close this subsection with the monotonicity of modified localized virial functional below
\begin{align} \label{wtI-def}
     \mathscr{L}
           :=\sum_{k=1}^K  \frac{\g_k}{2\lbb_k} {\rm Im}
            \int (\nabla\chi_A) \(\frac{x-\alpha_{k}}{\lambda_k}\)\cdot\nabla R\ol{R}\Phi_kdx
            - \sum_{k=1}^K  \frac{\g_k^2}{8\lbb_k^2} \|xQ\|_2^2,
\end{align}
where $\chi_A$ is as in \eqref{I-def}.

\begin{remark}
The modified localized virial functional \eqref{wtI-def} is introduced
mainly to derive
the refined estimate \eqref{Rj-lbbj3-refine} below for the second upgradation purpose.
Actually, it enables to cancel the bad $\calo(1)$ terms in the refined energy estimate \eqref{energy-esti-refined}
and thus allows to  further upgrade the convergence rate to the second order,
see Subsection \ref{Subsec-D-2} below.
\end{remark}

The main monotonicity of modified localized virial functional is formulated in Theorem \ref{Thm-Modi-Loc-virial}.

\begin{theorem} \label{Thm-Modi-Loc-virial}
(Monotonicity of modified localized virial functional)
We have that for $t$ close to $T$,
\begin{align}  \label{dt-wtI}
\frac{d \mathscr{L}}{dt}
\geq& C \sum_{k=1}^{K} \frac{\g_k}{\lambda_k^4} \int |\nabla \varepsilon_k|^2e^{-\frac{|y|}{A}}+|\varepsilon_k|^2 dy
       -\sum_{k=1}^{K}\frac{\g_k}{\lbb_k^4}M_k
    + \calo(Er),
\end{align}
where $C>0$ is independent of $t$,
and the error term $Er$ is given by \eqref{Error}.
\end{theorem}

\begin{remark}
Note that,
the positivity of $\g_k$, implied by \eqref{gamj-wj-t},
is important to derive the  monotonicity
of  modified localized virial functional.
\end{remark}

{\bf Proof of Theorem \ref{Thm-Modi-Loc-virial}.}
We first note that,
by \eqref{lbb-a-t-0} and \eqref{gamj-wj-t},
\begin{align}  \label{gamj2-lbbj2-dt}
  \frac{d}{dt} \(\frac{\g_k^2}{\lbb_k^2}\)
=\frac{2\g_k}{\lbb_k^4} \(\lbb_k^2\dot{\g_k}+\g_k^2\)-\frac{2\g^2_k}{\lbb^2_k}
\frac{\lbb_k\dot{\lbb_k}+\g_k}{\lbb_k^2}
=\frac{2\g_k}{\lbb_k^4}\(\lbb_k^2\dot{\g_k}+\g_k^2\)+\calo\(\frac{Mod}{\lbb_k^2}\).
\end{align}
Using \eqref{Mod-bdd} and \eqref{lbb2g+g2-f''Qve2}
we get
\begin{align} \label{lv3}
   \frac{d}{dt} \(\sum_{k=1}^K \frac{\g_k^2}{8\lbb_k^2}\|xQ\|_2^2\)
=&\sum_{k=1}^{K}\frac{\g_k}{\lbb_k^4}M_k - \sum_{k=1}^{K}\frac{\g_k}{\lbb_k^4}\int |R|^2\Phi_kdx  \nonumber \\
 & +\sum_{k=1}^{K}\frac{\g_k}{\lbb_k^4}{\rm Re}\int
   (1+\frac 2d)|Q_k|^{\frac {4}{d}}|\varepsilon_k|^2
     + \frac 2d |Q_k|^{\frac 4d-2} Q_k^2 \ol{\varepsilon_k}^2 dy \nonumber \\
 & + \sum_{k=1}^{K}\frac{\g_k}{\lbb_k^4}{\rm Re}\int
   (y \cdot\nabla \ol{Q_k}) (f''(Q_k)\cdot \ve_k^2) dy + \calo\(\wt{Er}\),
\end{align}
where the error term $\wt{Er}$ is as in \eqref{wtError}.

Moreover,
arguing as in the proof of \eqref{lv2} (see also \cite[Lemma 5.12]{SZ20}) we get
\begin{align}  \label{dtwtI.1}
 &\frac{d}{dt} \bigg(\sum_{k=1}^K \frac{\g_k}{2\lbb_k} {\rm Im}
            \int \nabla\chi_A \(\frac{x-\alpha_{k}}{\lambda_k}\)\cdot\nabla R\ol{R}\Phi_kdx \bigg) \nonumber \\
 =&   \sum_{k=1}^{K} \frac{\g_k}{\lambda^4_k} {\rm Re} \int \nabla^2\chi_A(y)(\nabla \ve_k,\nabla \ol{\ve_k}) dy
          - \sum_{k=1}^{K} \frac{\g_k}{4\lambda_k^4}\int \Delta^2\chi_A(y)|\ve_k|^2 dy
            \nonumber \\
 & + \sum_{k=1}^{K} \frac{\g_k}{\lbb^4_k}  {\rm Re}\int
   (\nabla\chi_A (y)\cdot\nabla \ol{Q_k}) (f''(Q_k) \cdot \ve_k^2)  dy \nonumber \\
   &   + \calo\(\frac{D^2(t)}{(T-t)^2}
      + \frac{Mod\|R\|_{L^2}}{(T-t)^3}
      +e^{-\frac{\delta}{T-t}}\) .
\end{align}

Therefore, combining   \eqref{lv3} and \eqref{dtwtI.1} altogether
and using similar arguments as those below \eqref{dt-I.2}
we obtain \eqref{dt-wtI}.
The proof is complete.
\hfill $\square$

\subsection{Energy quantization} \label{Subsec-Energy-Quant}

This subsection deals with the key energy quantization,
that is,
the energy of multi-bubble blow-up solutions to \eqref{equa-NLS}
indeed admits the quantization into the sum
of the energies of pseudo-conformal blow-up solutions.
The main result is stated in Theorem \ref{Thm-Energy-Ident}.

\begin{theorem} (Energy quantization) \label{Thm-Energy-Ident}
Let $v$ be the multi-bubble blow-up solution to \eqref{equa-NLS} satisfying
\eqref{v-S-H1-o1-Thm} and \eqref{nav-naS-iint-0+-Thm}.
Then, we have
\begin{align} \label{energy-ident}
E(v)=\sum_{k=1}^K\frac{\omega_k^2}{8}\|yQ\|_{L^2}^2 = \sum_{k=1}^K E(S_k),
\end{align}
where $S_k$ is the pseudo-conformal blow-up solution given by \eqref{Sj-blowup}, $1\leq k\leq K$.
\end{theorem}

The important consequence of Theorems \ref{Thm-Energy} and \ref{Thm-Energy-Ident}
is the following refined energy estimate.

\begin{corollary} \label{Cor-energy-refined}
(Refined energy estimate)
There exist $C_1, C_2>0$
such that
for any $t$ close to $T$,
\begin{align}  \label{energy-esti-refined}
      \sum_{k=1}^{K}\frac{|\beta_{k}|^2 }{2\lbb_{k}^2} \|Q\|_{2}^2
        + C_1 \frac{D^2}{(T-t)^2}
\leq& \frac{\|yQ\|_{L^2}^2}{8} \sum\limits_{k=1}^K \(\omega_k^2 - \frac{\g_k^2}{\lbb_k^2}\)
     + \sum_{k=1}^K \frac{|\beta_k|^2}{2\lbb_k^2} M_k
    + \sum_{k=1}^{K}\frac{M_k}{2\lbb_{k}^2}  \nonumber \\
    &  + C_2 \(\sum\limits_{k=1}^K \frac{M_k^2}{(T-t)^2} + e^{-\frac{\delta}{T-t}}\).
\end{align}
\end{corollary}

The remaining of this subsection is devoted to the proof of Theorem \ref{Thm-Energy-Ident}.
The key ingredient of proof is the $H^1$ dispersion of remainder
as formulated in Proposition \ref{Prop-D-beta-tn} below.
We shall first prove Lemmas \ref{Lem-H1-out} and \ref{Lem-R-H1-In},
which give the average $H^1$ dispersion of remainder.

Let $B_{2\sigma}(x_k):=\{x\in \bbr^d: |x-x_k|\leq 2\sigma\}$, $1\leq k\leq K$,
where $\sigma$ is as in \eqref{sep-xj-0}. Let $\Omega_{2\sigma}:= \cup_{k=1}^K B_{2\sigma}(x_k)$
be the interior region near the singularities
and $\wt \Omega_{2\sigma}:= \bbr^d \setminus \Omega_{2\sigma}$ be the exterior region.
Let $\wt R := v-S$ with $S= \sum_{k=1}^K S_k$.

\begin{lemma} ($H^1$ dispersion away from the  singularities) \label{Lem-H1-out}
For $t$ close to $T$, we have
\begin{align}  \label{R-H1-out}
   \|\na R(t)\|_{L^2(\wt \Omega_{2\sigma})}
   \leq  \|\na \wt R(t)\|_{L^2(\wt \Omega_{2\sigma})}  + \calo(e^{-\frac{\delta}{T-t}}).
\end{align}
In particular,
\begin{align}  \label{R-H1-out-tn-0}
\lim_{t\to T}\frac{1}{T-t}\int_{t}^{T}\frac{1}{T-s}\int_{s}^{T}\|\nabla R(r)\|^2_{L^2(\wt \Omega_{2\sigma})}drds=0.
\end{align}
\end{lemma}

{\bf Proof.}
We note from  \eqref{v-dec} that
\begin{align*}
   R= S-U+\wt R,
\end{align*}
which   yields that
\begin{align}
   \|\na R\|_{L^2(\wt \Omega_{2\sigma})}
   \leq   \|\na U\|_{L^2(\wt \Omega_{2\sigma})}
           +\|\na S\|_{L^2(\wt \Omega_{2\sigma})} + \|\na \wt R\|_{L^2(\wt \Omega_{2\sigma})}.
\end{align}

Note that,
by \eqref{lbb-a-t-0}, $|\a_k-x_k|=o(1)$,
we may take $t$ sufficiently close to $T$ such that $|\a_k -x_k|<\sigma$,
and so $|x-\a_k|>\sigma$ for any $x\in \wt \Omega_{2\sigma}$.
Then,
by the explicit expression of $U_k$ in \eqref{Uj-Qj-Q}, \eqref{lbbk-t-approx} and
\eqref{Q-decay},
\begin{align}
  \|\na U\|^2_{L^2(\wt \Omega_{2\sigma})}
  \leq \sum\limits_{k=1}^K \int_{|x-\a_k|\geq \sigma} |\na U_k|^2 dx
  \leq \sum\limits_{k=1}^K \lbb_k^{-2} \int_{|y|\geq \frac{\sigma}{\lbb_k}} (\na Q_k)^2(y) dy
  \leq C e^{-\frac{\delta}{T-t}},
\end{align}
where $C,\delta>0$.
Similarly, we have
\begin{align}
  \|\na S\|^2_{L^2(\wt \Omega_{2\sigma})} \leq C e^{-\frac{\delta}{T-t}}.
\end{align}
Thus,
putting the above estimates altogether we prove \eqref{R-H1-out}.
Estimate \eqref{R-H1-out-tn-0} then follows from \eqref{R-H1-out}
and   \eqref{nav-naS-iint-0+-Thm}.
\hfill $\square$

The more delicate interior region near the  singularities  is treated in Lemma \ref{Lem-R-H1-In} below.

\begin{lemma} ($H^1$ dispersion in average near the  singularities) \label{Lem-R-H1-In}
For $t$ close to $T$, we have
\begin{align} \label{R-H1-in-integ}
   \lim_{t \to T} \frac{1}{T-t} \int_{t}^{T} \frac{1}{T-s} \int_s^T \|\nabla R(r)\|^2_{L^2(\Omega_{2\sigma})} dr ds =0.
\end{align}
\end{lemma}

{\bf Proof.}
Let $\chi$ be a smooth radial cutoff function
such that $\nabla \chi(x)=x$ for $|x|\leq 2\sigma$,
$\nabla^2\chi$ is supported in $|x|\leq 3\sigma$ and
is positive semidefinite on $\R^d$.
Let $\chi_k :=\chi(x-x_k)$, $1\leq k\leq K$.
We use the local virial functional $ \mathbb{L}$  defined by
\begin{align} \label{I-virial-in}
   \mathbb{L}:=\sum_{k=1}^{K}{\rm Im}\int\nabla\chi_k\cdot \nabla R \ol{R}\Phi_k dx,
\end{align}
where $\Phi_k$, $1\leq k\leq K$, are the localization functions given by \eqref{phi-local}.

Note that, by \eqref{R-t-1},
\begin{align} \label{estr}
| \mathbb{L}(t)|\leq C\|\nabla R\|_{L^2}\|R(t)\|_{L^2}\leq C\|R(t)\|_{L^2}.
\end{align}

Moreover,
using the integration by parts formula we have
\begin{align*}
   \frac{d\mathbb{L}}{dt}
   = \sum\limits_{k=1}^K {\rm Re} \<\Delta \chi_k R_k, i\partial_t R\>
     + \sum\limits_{k=1}^K {\rm Re} \< \na\chi_k \cdot (\na R_k + \na R \Phi_k), i\partial_t R\>,
\end{align*}
where $R_k = R\Phi_k$.
Then, in view of  equation \eqref{equa-R},
we obtain
\begin{align} \label{dI-I1-I4}
      \frac{d \mathbb{L}}{dt}
   =&-\sum_{k=1}^{K}{\rm Re} \< \Delta \chi_k R_k, \Delta R\>
     -\sum_{k=1}^{K}{\rm Re} \< \na \chi_k\cdot (\na R_k+ \na R \Phi_k), \Delta R\> \nonumber \\
    & - \sum_{k=1}^{K}{\rm Re} \< \Delta \chi_k R_k + \na \chi_k \cdot (\na R_k + \na R \Phi_k), \eta\> \nonumber \\
   &- \sum_{k=1}^{K}{\rm Re} \< \Delta \chi_k R_k + \na \chi_k \cdot (\na R_k + \na R \Phi_k), |v|^{\frac 4d} v - |U|^{\frac 4d}U\> \nonumber \\
   =&: \sum_{k=1}^{K} ( \mathbb{L}_{t,k1} + \mathbb{L}_{t,k2} + \mathbb{L}_{t,k3} + \mathbb{L}_{t,k4}).
\end{align}
Below we estimate the four terms
$\mathbb{L}_{t,ki}$, $1\leq i\leq 4$, separately.

First,
since the supports of $\Delta \chi_k$ and $R_j$ are disjoint for any $j\not =k$, we see that,
\begin{align} \label{I1j}
 \mathbb{L}_{t,k1} = -{\rm Re} \<\Delta \chi_k R_k, \Delta R_k \>
        =& \int \Delta \chi_k|\nabla R_k|^2dx  -\frac12\int\Delta^2 \chi_k|R_k|^2dx \nonumber \\
        =& \int \Delta \chi_k |\na R_k|^2 dx + \calo(\|R\|_{L^2}^2).
\end{align}

Next, we claim that
\begin{align} \label{I2j}
   \mathbb{L}_{t,k2} =&2 {\rm Re} \int\nabla^2 \chi_k(\nabla R_k, \nabla \overline{R_k})dx
            -\int \Delta \chi_k|\nabla R_k|^2dx
      + 2 {\rm Re} \<\na \chi_k \cdot \na R, \na \Phi_k \cdot \na R \>
     +\calo(\|R\|_{L^2}\|\nabla R\|_{L^2}).
\end{align}

In order to prove \eqref{I2j},
we note from $ \na R \Phi_k = \na R_k - R \na \Phi_k$ that
\begin{align} \label{I2j.0}
    \mathbb{L}_{t,k2} = - 2 {\rm Re} \<\na \chi_k \cdot \na R_k, \Delta R\>
             + {\rm Re} \<\na \chi_k \cdot  \na \Phi_k R, \Delta R\>
           =: \mathbb{L}_{k,21} + \mathbb{L}_{k,22}.
\end{align}
Using the integration by parts formula we compute
\begin{align*}
   \mathbb{L}_{k,21}  = 2 {\rm Re}  \int \na^2 \chi_k (\na R_k, \na \ol{R}) dx
              - 2 {\rm Re} \int \Delta \chi_k \na R_k \cdot \na \ol{R} dx
              - 2 \sum\limits_{i,h=1}^d {\rm Re} \<\partial_i \chi_k \partial_h R_k, \partial_{ih} R\>.
\end{align*}
Since
the supports of $\na^2 \chi_k$ and $\Delta \chi_k$
are disjoint with $R_j$, $j\not = k$,
we infer that
\begin{align} \label{I21j}
    \mathbb{L}_{k,21} = 2{\rm Re} \int \na^2 \chi_k(\na R_k, \na \ol{R_k}) dx
              -2  \int \Delta \chi_k |\na R_k|^2 dx
              - 2 \sum\limits_{i,h=1}^d {\rm Re} \<\partial_i \chi_k \partial_h R_k, \partial_{ih} R\>.
\end{align}

Moreover, we use the integration by parts formula again to compute
\begin{align*}
   \mathbb{L}_{k,22}
   =& - {\rm Re}  \int \na \chi_k \cdot \na \Phi_k |\na R|^2 dx
      + \calo(\|R\|_{L^2} \|\na R\|_{L^2}) \nonumber \\
   =& {\rm Re} \int \Delta \chi_k \Phi_k |\na R|^2 dx
      + 2 \sum\limits_{i,h=1}^d {\rm Re} \< \partial_i \chi_k\Phi_k\partial_h R,   \partial_{ih} R\>
      + \calo(\|R\|_{L^2} \|\na R\|_{L^2}).
\end{align*}
Since on the support of $\Delta \chi_k$,
$\Phi_k =1$ and $\Phi_j =0$ for $j\not =k$,
it follows that
\begin{align} \label{I22j}
   \mathbb{L}_{k,22}
    = {\rm Re} \int \Delta \chi_k |\na R_k|^2 dx
      + 2 \sum\limits_{i,h=1}^d {\rm Re} \<\partial_i \chi_k \Phi_k \partial_h R,  \partial_{ih} R\>
      + \calo(\|R\|_{L^2} \|\na R\|_{L^2}).
\end{align}
Note  that
\begin{align} \label{I21j-I22j}
    - {\rm Re} \<\partial_i \chi_k \partial_h R_k, \partial_{ih} R\>
     + {\rm Re} \<\partial_i \chi_k \Phi_k  \partial_h R, \partial_{ih} R \>
    =  {\rm Re} \<\partial_i \chi_k  \partial_i R, \partial_h \Phi_k \partial_{h} R\>
      + \calo(\|R\|_{L^2} \|\na R\|_{L^2}).
\end{align}
Thus, plugging \eqref{I21j}, \eqref{I22j} and \eqref{I21j-I22j} into \eqref{I2j.0}
we obtain \eqref{I2j},
as claimed.

Regarding the third term $\mathbb{L}_{t,k3}$,
we move the derivative to $\eta$ and use \eqref{mod-t-2} and \eqref{eta-Mod-esti} to get
\begin{align} \label{I3j}
   |\mathbb{L}_{t,k3}| =\calo(\|R\|_{L^2} \|\eta\|_{H^1})
             = \calo\(\frac{Mod\|R\|_{L^2}}{(T-t)^3} + e^{-\frac{\delta}{T-t}}\|R\|_{L^2}\)
             =\calo\(\frac{\|R\|_{L^2}}{T-t}\).
\end{align}

It remains to treat the last term $\mathbb{L}_{t,k4}$
on the R.H.S. of \eqref{dI-I1-I4}.
First, since
\begin{align} \label{fv-fU}
  \bigg||v|^{\frac 4d}v-|U|^{\frac 4d}U\bigg|
  \leq C\sum_{l=1}^{1+\frac{4}{d}}|U|^{1+\frac{4}{d}-l}|R|^{l},
\end{align}
using \eqref{R-D-Lp} and \eqref{R-t-1} we get
\begin{align}   \label{I4j.1}
   |{\rm Re} \<\Delta \chi_k R_k, |v|^{\frac 4d}v-|U|^{\frac 4d}U \>|
   = \calo\(\frac{\|R\|^2_{L^2}+D^3}{(T-t)^2}\)
   = \calo\(\frac{\|R\|_{L^2}}{T-t} + \frac{D^3}{(T-t)^2}\).
\end{align}

Moreover,
we claim that
\begin{align}  \label{I4j.2}
  |{\rm Re}\<\na \chi_k\cdot (\na R_k + \na R \Phi_k), |v|^{\frac 4d}v-|U|^{\frac 4d}U\>|
  =\calo \(\frac{\|R\|_{L^2}}{T-t}+ \frac{D^3}{(T-t)^2}+ e^{-\frac{\delta}{T-t}}\).
\end{align}
This along with \eqref{I4j.1} yields that
\begin{align} \label{I4j}
    \mathbb{L}_{t,k4} =\calo \( \frac{\|R\|_{L^2}}{T-t} + \frac{D^3}{(T-t)^2} + e^{-\frac{\delta}{T-t}}\).
\end{align}

In order to prove \eqref{I4j.2},
we use  \eqref{fv-fU} and Lemma \ref{Lem-inter-est} in Appendix to bound
\begin{align} \label{nachij-fvfU-esti.0}
       {\rm L.H.S.\ of\ } \eqref{I4j.2}
       \leq& C\bigg( \int |\nabla \chi_k| |U_k|^{\frac{4}{d}}|R_k\nabla R_k|dx
          +\int |\nabla \chi_k| |U_k|^{\frac{4}{d}-1}|R_k^2\nabla R_k|dx \nonumber \\
      & + \sum\limits_{l=3}^{1+\frac 4d} \int |\na \chi_k| |U|^{1+\frac 4d -l} |R|^l |\na R| dx
        + \sum\limits_{l=3}^{1+\frac 4d} \int |\na \chi_k| |U|^{1+\frac 4d -l} |R|^{l+1} dx
               + e^{-\frac{\delta}{T-t}}
      \bigg).
\end{align}

Similarly to \eqref{I4j.1},
we use \eqref{R-D-Lp} to bound
\begin{align}
    & \sum\limits_{l=3}^{1+\frac 4d}\int |\na \chi_k| |U|^{1+\frac 4d -l} |R|^{l+1} dx
    \leq  C \frac{\|R\|_{L^2}^2 + D^3}{(T-t)^2}
    \leq C \(\frac{\|R\|_{L^2}}{T-t} + \frac{D^3}{(T-t)^2}\), \label{nachi-U-D} \\
   & \sum\limits_{l=3}^{1+\frac 4d}\int |\na \chi_k|  |U|^{1+\frac 4d -l} |R|^l |\na R| dx
    \leq    C (T-t)^{-\frac d2 -2 + \frac d2 l} \|\na R\|_{L^2} \|R\|_{L^{2l}}^l
    \leq  C \frac{D^3}{(T-t)^{2}}.
\end{align}

Next, we deal with the first two integrations on the R.H.S. of
\eqref{nachij-fvfU-esti.0}.

Note that for $t$ close to $T$,
$|\a_k(t)-x_k|\leq \sigma$
and  $\nabla \chi(x)=x$ for $|x|\leq 2\sigma$.
Using the renormalized variable $\ve_k$ defined by \eqref{Rj-ej}, \eqref{lbbk-t-approx},
\eqref{P-t-1}, \eqref{R-t-1}
and \eqref{Q-decay} below we get
\begin{align}
\int |\nabla \chi_k| &|U_k|^{\frac{4}{d}}|R_k\nabla R_k|dx
\leq \frac{C}{(T-t)^3}\int |\nabla \chi(\lbb_k y+\a_k-x_k)| Q^{\frac{4}{d}}(y)|\varepsilon_k\nabla \varepsilon_k|dy  \nonumber \\
&\leq C \frac{|\lbb_k |+|\a_k-x_k|}{(T-t)^3}\int_{|y|\leq \frac{\sigma}{\lbb_k}} (1+|y|) Q^{\frac{4}{d}}(y)|\varepsilon_k\nabla \varepsilon_k|dy
  +\frac{C}{(T-t)^3}\int_{|y|\geq \frac{\sigma}{\lbb_k}}  Q^{\frac{4}{d}}(y)|\varepsilon_k\nabla \varepsilon_k|dy \nonumber \\
&\leq C\(\frac{\|R\|_{L^2}}{T-t}+e^{-\frac{\delta}{T-t}}\).
\end{align}
Moreover, by \eqref{R-t-1} and the Gagliardo-Nirenberg inequality \eqref{G-N},
\begin{align} \label{nachiR2R}
\int |\nabla \chi_k| |U_k|^{\frac{4}{d}-1}|R^2_k\nabla R_k|dx
\leq& C (T-t)^{-2+\frac d2} \|R_k\|_{L^4}^2 \|\na R_k\|_{L^2}  \nonumber \\
\leq& C (T-t)^{-2+\frac d2} \|R_k\|_{L^2}^{2-\frac d2} \|\na R_k\|_{L^2}^{1+\frac d2}
\leq C \frac{\|R\|_{L^2}}{T-t}.
\end{align}

Thus, inserting estimates \eqref{nachi-U-D}-\eqref{nachiR2R} into \eqref{nachij-fvfU-esti.0}
we obtain \eqref{I4j.2}, as claimed.

Now, combining estimates  \eqref{dI-I1-I4}, \eqref{I1j}, \eqref{I2j}, \eqref{I3j}
and \eqref{I4j} altogether we conclude
\begin{align}  \label{naRk-dL}
        2 {\rm Re} \int \na^2 \chi_k(\na R_k, \na \ol{R_k}) dx
   \leq&  \frac{d\mathbb{L}}{dt} + 2|{\rm Re} \< \na \chi_k\cdot \na R, \na \Phi_k\cdot \na R\>|
      + C\(\frac{\|R\|_{L^2}}{T-t}+\frac{D^3}{(T-t)^2} + e^{-\frac{\delta}{T-t}}\).
\end{align}

Note that,
since the support of $\na \Phi_k$ is in the exterior region $\wt \Omega_{2\sigma}$,
Lemma \ref{Lem-H1-out} yields
\begin{align} \label{L-o1}
   |{\rm Re} \< \na \chi_k\cdot \na R, \na \Phi_k\cdot \na R\>|
   \leq C \|\na R\|^2_{L^2(\wt \Omega_{2\sigma})}
   \leq C(\| \na \wt R\|^2_{L^2} + e^{-\frac{\delta}{T-t}}).
\end{align}
Taking into account
$\na^2\chi_k(\na R_k, \na \ol{R_k})\geq  |\na R_k|^2$ for $x\in B_{2\sigma}(x_k)$,
$\|\na R\|^2_{L^2(\Omega_{2\sigma})} = \sum_{k=1}^K \|\na R_k\|^2_{L^2(B_{2\sigma}(x_k))}$,
integrating \eqref{naRk-dL} from $t$ to $T$ and using
\eqref{beta-gam-D-t-1} and the boundary estimate \eqref{estr}
we get
\begin{align}
\int_{t}^{T}\|\nabla R(s)\|^2_{L^2(\Omega_{2\sigma})}ds
   \leq C \(\|R(t)\|_{L^2}+  \int_{t}^{T}\frac{\|R(s)\|_{L^2}}{T-s} ds
         + \int_t^T \| \na \wt R\|^2_{L^2} ds
         +  o(T-t)\),
\end{align}
which along with  the Cauchy inequality
and the space time estimate \eqref{Rj-lbbj3-bdd} yields
\begin{align} \label{aver.1}
\frac{1}{T-t} \int_{t}^{T}\|\nabla R(s)\|^2_{L^2(\Omega_{2\sigma})}ds
\leq& C \(\frac{\|R(t)\|_{L^2}}{T-t}+ \left(\int_{t}^{T}\frac{\| R(s)\|^2_{L^2}}{(T-s)^3}ds\right)^{\half}
         + \frac{1}{T-t} \int_t^T \| \na \wt R\|^2_{L^2} ds \)  + o(1) \nonumber \\
\leq& C \(\frac{\|R\|_{L^2}}{T-t}  + \frac{1}{T-t} \int_t^T \| \na \wt R\|^2_{L^2} ds\)  + o(1).
\end{align}
Furthermore,
taking the average of \eqref{aver.1} again
and then using \eqref{nav-naS-iint-0+-Thm} and \eqref{Rj-lbbj3-bdd}
we arrive at
\begin{align}
      &  \frac{1}{T-t} \int_t^T \frac{1}{T-s} \int_s^T  {\|\na R\|_{L^2(\Omega_{2\sigma})}^2} dr ds   \nonumber \\
  \leq& C\(\frac{1}{T-t} \int_t^T \frac{\|R\|_{L^2}}{T-s} ds
           +  \frac{1}{T-t} \int_t^T \frac{1}{T-s} \int_s^T \| \na \wt R\|^2_{L^2} drds \)
         + o(1) \nonumber \\
  \leq& C \(\int_t^T \frac{\|R\|_{L^2}^2}{(T-s)^3} ds\)^{\frac 12} + o(1)
  \to 0, \ \ as\ t\to T.
\end{align}

Therefore,
the proof of Lemma \ref{Lem-R-H1-In} is complete.
\hfill $\square$

By virtue of Lemmas \ref{Lem-H1-out} and \ref{Lem-R-H1-In},
we have the average estimates of $\|R\|_{L^2}$ and $\beta$ below.

\begin{proposition} \label{Prop-R-beta-tn}
There exists a universal constant $C>0$
such that for $t$ close to $T$,
\begin{align} \label{R-beta-integ-t3}
  \frac{1}{T-t} \int_t^T \frac{\|R\|_{L^2}^2 + |\beta|^2}{(T-s)^2} ds
  \leq C <\9.
\end{align}
\end{proposition}

{\bf Proof.}
By direct calculations,
\begin{align} \label{dt-betaj-lbbj}
\(\frac{(\a_k-x_k)\cdot \b_k}{\lbb_k}\)_t
&=\frac{(\lbb_k\dot{\a_k}-2\b_k)\cdot \b_k+2|\b_k|^2}{\lbb_k^2}
+\frac{(\lbb_k^2\dot{\b_k}+\b_k\g_k)\cdot (\a_k-x_k)}{\lbb_k^3}
  -\frac{(\lbb_k\dot{\lbb_k}+\g_k)(\a_k-x_k)\cdot \b_k}{\lbb_k^3}  \nonumber \\
&=2\frac{|\b_k|^2}{\lbb_k^2}
   +\calo\(\bigg|\frac{\b_k}{\lbb^2_k}\bigg|+\bigg|\frac{\a_k-x_k}{\lbb_k^3}\bigg|
   +\bigg|\frac{(\a_k-x_k)\b_k}{\lbb_k^3}\bigg|\)Mod.
\end{align}
Note that, by \eqref{lbb-a-t-0}, \eqref{beta-gam-D-t-1} and \eqref{mod-t-2},
\begin{align} \label{betajlbbj-o1}
\bigg|\frac{(\a_k-x_k)\cdot \b_k}{\lbb_k}\bigg|\leq C|\a_k-x_k| = o(T-t),
\end{align}
and
\begin{align} \label{betaj-aj-Mod-o1}
\(\bigg|\frac{\b_k}{\lbb^2_k}\bigg|+\bigg|\frac{\a_k-x_k}{\lbb_k^3}\bigg|
+\bigg|\frac{(\a_k-x_k)\cdot \b_k}{\lbb_k^3}\bigg| \)Mod = o(1).
\end{align}
Hence, integrating \eqref{dt-betaj-lbbj} from $t$ to $T$
and using \eqref{lbbk-t-approx}, \eqref{betajlbbj-o1} and \eqref{betaj-aj-Mod-o1} we get
\begin{align} \label{beta-t2-ot}
    \int_{t}^{T}\frac{|\b_k|^2}{(T-s)^2}ds = o(T-t),\ \ 1\leq k\leq K.
\end{align}

Moreover, in view of \eqref{lbbk-t-approx} and \eqref{Rj-lbbj3-bdd}, we get
\begin{align} \label{RL2-t3-bdd}
 \frac{1}{T-t} \int_{t}^{T}\frac{\|R\|_{L^2}^2}{(T-s)^2}ds
\leq C\sum_{k=1}^{K}\int_{t}^{T}\frac{\|R_k\|_2^2}{\lbb_k^3}ds\leq C<\9.
\end{align}

Thus, \eqref{beta-t2-ot} and \eqref{RL2-t3-bdd} yield
\eqref{R-beta-integ-t3} immediately.
\hfill $\square$

As a consequence of Lemmas \ref{Lem-H1-out}, \ref{Lem-R-H1-In} and Proposition \ref{Prop-R-beta-tn},
we are now able to upgrade the estimates of
remainder and parameter $\beta$ along a sequence,
which in particular improve the previous ones in \eqref{beta-gam-D-t-1}.

\begin{proposition} \label{Prop-D-beta-tn} ($H^1$ dispersion along a sequence)
There exists a sequence $\{t_n\}$ to $T$  such that
\begin{align}  \label{D-beta-o1+-tn}
\lim_{n\to +\infty}\frac{D(t_n)+ |\b(t_n)|}{T-t_n}=0.
\end{align}
In particular,
\begin{align} \label{R-naR-beta-tn}
   \|R(t_n)\|_{L^2}  + |\b(t_n)| = o(T-t_n),\ \
   \|\na R(t_n)\|_{L^2} = o(1).
\end{align}
\end{proposition}

Now, we are ready to prove the key energy quantization in Theorem \ref{Thm-Energy-Ident}.

{\bf Proof of Theorem \ref{Thm-Energy-Ident}.}
We mainly focus on the proof of the first equality in \eqref{energy-ident},
as the second one follows from straightforward computations
by using the  expression of $S_k$ in \eqref{Sj-blowup}.
It suffices to prove this for a subsequence,
due to the energy conservation law of energy.

We infer from \eqref{Ev-expand.0}, \eqref{err-esti-energy} and \eqref{Ev-linear-esti} that
\begin{align}
E(v)
=& \sum_{k=1}^K\( \frac{|\beta_{k}|^2}{2\lambda_{k}^2} \|Q\|_{L^2}^2
+ \frac{\gamma_{k}^2}{8\lambda_{k}^2} \|yQ\|_{L^2}^2
 -\frac{|\b_k|^2}{2\lbb_k^2}M_k
-\frac{M_k}{2\lbb_k^2}\)  \nonumber \\
&+\frac{1}{2}{\rm Re} \int |\nabla R|^2+\sum_{k=1}^{K}\frac{1}{\lbb_{k}^2}|R|^2  \Phi_k
     -(1+\frac2d)|U|^{\frac 4d} |R|^2
     - \frac 2d |U|^{\frac 4d-2}U^2\overline{R}^2 dx  \nonumber\\
&  + \calo\( \sum_{k=1}^{K}\frac{|\b_k|^2}{\lbb_k^2}\|R\|^2_{L^2}+ \frac{D^3}{(T-t)^2} +e^{-\frac{\delta}{T-t}}\)  \nonumber \\
=:&\sum_{i=1}^{3}J_i. \label{33}
\end{align}

Take the sequence $\{t_n\}$ as in Proposition \ref{Prop-D-beta-tn}.
Thus, by \eqref{lbb-a-t-0}, \eqref{glbb-omega-o1} and \eqref{R-naR-beta-tn},
we have that along the sequence $\{t_n\}$,
\be
J_1=\sum_{k=1}^{K}\frac{\omega_k^2}{8}\|yQ\|_{L^2}^2+o(1),
\ee
and
\be
J_2 =\calo\(\|\na R(t_n)\|_{L^2}^2 + \frac{\|R(t_n)\|_{L^2}^2}{(T-t_n)^2}\) = o(1).
\ee
Using \eqref{lbb-a-t-0}, \eqref{beta-gam-D-t-1} and \eqref{R-naR-beta-tn} again
we infer that along the sequence $\{t_n\}$,
\begin{align}
    J_3 = o(1).
\end{align}

Thus, taking into account the conservation law of energy we arrive at
\begin{align}
E(v(t))=\lim_{n\rightarrow\infty}E(v(t_n))=\sum_{k=1}^{K}\frac{\omega_k^2}{8}\|yQ\|_{L^2}^2.
\end{align}

Therefore, the proof of Theorem \ref{Thm-Energy-Ident} is complete.
\hfill $\square$

\section{Upgradation of the convergence rate of remainder} \label{Sec-Upgrad-Remainder}

This section is mainly devoted to the crucial upgradation of the convergence rate of remainder.
First in Subsection \ref{Subsec-D-2},
we upgrade the convergence rate to the second order,
by virtue of the refined energy estimate,
the monotonicity of modified localized virial functional
and the relationship \eqref{Mj-D}-\eqref{Mj-glbb4-D} between
localized mass and remainder.
Then, in order to further upgrade the convergence rate,
we prove the monotonicity of modified generalized energy in Subsection \ref{Subsec-Modf-Gen-Energy}.
Eventually, in Subsection \ref{Subsec-Upgra-Exp},
we complete the final upgradation
to the exponential decay rate.

\subsection{Upgradation to the second order}  \label{Subsec-D-2}

The main result in this step is Theorem \ref{Thm-D-t2} below,
which in  particular improves the previous estimates in \eqref{R-t-1} and \eqref{R-naR-beta-tn}
for any $t$ close to $T$.

\begin{theorem}  \label{Thm-D-t2}
There exists $C>0$ such that for $t$ close to $T$,
\begin{align} \label{D-t-2}
   D(t) \leq C(T-t)^2.
\end{align}
\end{theorem}

{\bf Proof.}
The proof proceeds in two steps below.

{\bf Step $1$. Upgradation to the order $1+$.}
We first use \eqref{wtI-def} and the precise asymptotic of $\g_k$ in \eqref{glbb-omega-o1}
to derive from \eqref{wtI-def} that
\begin{align*}
   \lim\limits_{\wt t\to T} \mathscr{L} (\wt t)  - \mathscr{L} (t)
   = \frac{\|yQ\|_{L^2}^2}{8} \sum\limits_{k=1}^K \(\frac{\g_k^2(t)}{\lbb_k^2(t)} - \omega_k^2\)
     + \calo(\|R(t)\|_{L^2} \|\na R(t)\|_{L^2}).
\end{align*}
Using \eqref{glbb-omega-o1} again we also have
that for some $\wt c>0$,
\begin{align*}
   \tilde c \int_t^T \frac{1}{\lbb_k^3} \|R_k\|_{L^2}^2 ds
   \leq \int_t^T \frac{\g_k}{\lbb_k^4} \|R\|_{L^2}^2 ds.
\end{align*}
Then, integrating both sides of  \eqref{dt-wtI}
we obtain
\begin{align}  \label{Rj-lbbj3-refine}
  \tilde c \sum_{k=1}^{K}\int_{t}^{T} \frac{\|R_{k}\|_{L^2}^2}{\lbb_{k}^3} ds
   \leq& \frac{\|yQ\|_{L^2}^2}{8} \sum_{k=1}^{K}\(\frac{\g_k^2(t)}{\lbb_k^2(t)}-\omega^2_k\)
         +C\|R(t)\|_{L^2}\|\nabla R(t)\|_{L^2} \nonumber \\
       &  + \bigg|\int_{t}^{T}\sum_{k=1}^{K}\frac{\g_k}{\lbb_k^4}M_kds\bigg|
          + C \int_{t}^{T} Er ds.
\end{align}
where $M_k$ and $Er$ are given by \eqref{Mj-def} and \eqref{Error}, respectively, and $\wt c, C>0$.

Thus, combining \eqref{Rj-lbbj3-refine} and the refined energy estimate \eqref{energy-esti-refined} altogether and using
\eqref{Mj-t2+}, Corollary \ref{Cor-beta-gam-D-t-1} and
the inequality
\begin{align}  \label{R-D*}
\|R\|_{L^2}\|\nabla R\|_{L^2}\leq C\|R\|_{L^2}\leq \delta \frac{D^2}{(T-t)^2}+ \frac{C^2}{\delta} (T-t)^2,
\end{align}
with $\delta$ small enough
we obtain that
there exists a positive constant $C>0$ such that
\begin{align} \label{D-iter-2}
&\sum_{k=1}^{K}\frac{|\beta_{k}(t)|^2 }{\lbb_{k}^2} \|Q\|_{2}^2
     +  \frac{D^2(t)}{(T-t)^2}
     +  \sum_{k=1}^{K}\int_{t}^{T} \frac{\|R_{k}\|_{L^2}^2}{\lbb_{k}^3}  ds   \nonumber \\
\leq& C\( (T-t)^2
      +\bigg|\sum_{k=1}^{K}\frac{M_k}{\lbb_{k}^2}\bigg|
      +\bigg|\sum\limits_{k=1}^K \int_{t}^{T}\frac{\g_k}{\lbb_k^4}{M_k}ds\bigg|
       +\int_{t}^{T}   \frac{D^2 + \sum_{k=1}^K |M_k|}{(T-s)^2} ds \).
\end{align}
Note that,
by \eqref{nav-naS-iint-0+-Thm}, \eqref{Mj-t2+} and \eqref{gamj-wj-t},
\begin{align} \label{intlbb4Mk}
   \bigg|\int_{t}^{T}\sum_{k=1}^{K}\frac{\g_k}{\lbb_k^4}M_kds\bigg|
   \leq& C \int_t^T \frac{1}{(T-s)^3} \int_s^T \int_r^T \|\wt R(\tau)\|_{L^2}^2 d\tau dr ds
         + C \int_t^T \frac{1}{(T-s)^3}e^{-\frac{\delta}{T-s}} ds \nonumber \\
   \leq& C \int_t^T \frac{1}{(T-s)^2} \int_s^T \frac{1}{T-r} \int_r^T \|\wt R(\tau)\|_{L^2}^2 d\tau dr ds
         + C e^{-\frac{\delta}{T-t}} \nonumber  \\
   \leq& C (T-t)^{\zeta}.
\end{align}
Using \eqref{Mj-t2+} again we get that
\begin{align} \label{lbbMk-intt2D2}
  \bigg|\sum\limits_{k=1}^K \frac{M_k}{\lbb_k^2} \bigg|
  + \int_t^T \frac{D^2+ \sum_{k=1}^K |M_k|}{(T-s)^2} ds
  \leq C(T-t)^{\zeta}.
\end{align}
Thus, plugging \eqref{intlbb4Mk} and \eqref{lbbMk-intt2D2} into \eqref{D-iter-2}
we obtain that for some $C>0$,
\begin{align} \label{D-t-1+}
D(t)\leq C (T-t)^{1+\frac \zeta 2}.
\end{align}
Note that,
estimate \eqref{D-t-1+} already improves the previous one \eqref{R-t-1}
and \eqref{R-naR-beta-tn} for any $t$ close to $T$.

Let us also mention that,
even though we have the high order $(T-t)^2$ on the R.H.S. of \eqref{D-iter-2},
the estimate \eqref{Mj-t2+} of localized mass $M_k$ actually restricts
the upgradation to the lower order $(T-t)^{1+\frac \zeta 2}$.

In order to further upgrade the convergence rate,
the key idea in the next step is to establish a Gronwall type inequality
by relating together the localized mass and the quantity $D$.

{\bf Step $2$. Upgradation to the second order.}
The goal in this step is to upgrade the convergence rate to the second order $(T-t)^2$, i.e., for some $C>0$,
\ben
\frac{D^2(t)}{(T-t)^2}\leq C(T-t)^2.
\enn

For this purpose, we shall establish  a Gronwall type estimate from \eqref{D-iter-2}.
Precisely, for $t$ close to $T$ such that $T-t\leq \varepsilon$,
using \eqref{lbbk-t-approx}, \eqref{Mj-D}-\eqref{Mj-glbb4-D}  and \eqref{D-iter-2}
we obtain the Gronwall type inequality
\begin{align}
  \frac{D^2(t)}{(T-t)^2}
   \leq & C_1\bigg( (T-t)^2+ \varepsilon\int_{t}^{T}\frac{D^2(s)}{(T-s)^3}ds
       + \ve \int_{t}^{T}\frac{1}{(T-s)^{3}}\int_{s}^{T}\frac{D^2(r)}{T-r}dr ds
     \bigg), \label{D-iter-iint}
\end{align}
where $C_1>0$ is independent of $\ve$.

Then, plugging \eqref{D-t-1+} into \eqref{D-iter-iint}
we get
\begin{align} \label{D-t-1+-0}
   \frac{D^2(t)}{(T-t)^2}
   \leq C_1 (T-t)^{2} + \frac{2C_1^2 \ve}{\zeta}  (T-t)^{\zeta}.
\end{align}

We use the induction arguments
and suppose that for some $n\geq 1$,
\begin{align} \label{D-itern-23}
   \frac{D^2(t)}{(T-t)^2}
   \leq \(\sum\limits_{k=0}^{n-1} C_1^{k+1} \ve^k\)  (T-t)^{2}
       + \frac{2^{n+1}C_1^{n+1} \ve^{n}}{\zeta^{n}}  (T-t)^{\zeta}.
\end{align}
Note that, \eqref{D-t-1+-0} verifies \eqref{D-itern-23}
at the preliminary step $n=1$.
Then,
plugging \eqref{D-itern-23} into \eqref{D-iter-iint}
and using straightforward computations
we see that estimate \eqref{D-itern-23} is still valid
with $n+1$ replacing $n$.
Thus, the induction arguments yield that,
for $\ve$ sufficiently small such that
$\frac{C_1 \ve}{\zeta} \leq \frac 14$,
\begin{align} \label{D-itern-23*}
   \frac{D^2(t)}{(T-t)^2}
   \leq \lim\limits_{n\to \9}
        \bigg( \(\sum\limits_{k=0}^{n-1} C_1^{k+1} \ve^k\)  (T-t)^{2}
       + \frac{2^{n+1}C_1^{n+1} \ve^{n}}{\zeta^{n}}  (T-t)^{\zeta} \bigg)
   \leq  2C_1  (T-t)^{2}.
\end{align}

Therefore, the proof of Theorem \ref{Thm-D-t2} is complete.
\hfill $\square$

As a consequence of Theorems \ref{Thm-Mod}, \ref{Thm-Loc-Mass} and \ref{Thm-D-t2}
and estimate \eqref{D-iter-2},
we obtain the improved estimates of localized mass and geometrical parameters below
\begin{corollary}  \label{Cor-Mk-Mod-4}
There exists $C>0$ such that for every $1\leq k\leq K$ and any $t$ close to $T$,
\begin{align}
   & |M_k(t)| + Mod(t) \leq C  (T-t)^4, \label{Mk-Mod-t-4} \\
    & |\beta_k(t)| \leq C  (T-t)^2.   \label{beta-t-2}
\end{align}
\end{corollary}

\subsection{Modified generalized energy}   \label{Subsec-Modf-Gen-Energy}

The last key ingredient to fulfill the upgradation to exponential decay rate is the modified generalized energy used in \cite{SZ20},
defined by
\begin{align} \label{gen-energy-def}
 \mathscr{I}:= &\frac{1}{2}\int |\nabla R|^2dx+\frac{1}{2}\sum_{k=1}^K\int\frac{1}{\lambda_{k}^2} |R|^2 \Phi_kdx
           -{\rm Re}\int F(u)-F(U)-f(U)\ol{R}dx \nonumber \\
&+\sum_{k=1}^K\frac{\gamma_{k}}{2\lambda_{k}}{\rm Im} \int (\nabla\chi_A) \(\frac{x-\alpha_{k}}{\lambda_{k}}\)\cdot\nabla R \ol{R}\Phi_kdx
 \ \ (:= \mathscr{I}^{(1)} + \mathscr{I}^{(2)}),
\end{align}
where $F(z):= \frac{d}{2d+4} |z|^{2+ \frac 4d}$
and  $f(z):= |z|^{\frac 4d}z$, $z\in \mathbb{C}$.

The key monotonicity of modified generalized energy
is formulated in Theorem \ref{Thm-I-mono} below.

\begin{theorem}  \label{Thm-I-mono}
(Monotonicity of modified generalized energy)
For any $\ve >0$ as in {\rm Case (I)} and {\rm Case (II)},
we have that for $t$ close to $T$,
\begin{align} \label{dIt-mono-case1}
\frac{d \mathscr{I}}{dt}
\geq& C_1 \sum_{k=1}^{K}\frac{1}{\lambda_{k}}
     \int \(|\nabla R_{k}|^2+\frac{1}{\lbb_{k}^2} |R_{k}|^2 \)
     e^{-\frac{|x-\alpha_{k}|}{A\lbb_k}}dx \nonumber \\
    & - C_2 \( \ve \frac{D^2(t)}{(T-t)^3}
            +\frac{D^2(t)}{(T-t)^2}
            + \sum_{k=1}^K \frac{|M_k(t)|}{(T-t)^2}
            + e^{-\frac{\delta}{T-t}}\),
\end{align}
where $C_1, C_2>0$ are independent of $\ve$ and $t$.
\end{theorem}

{\bf Proof.}
The proof is similar to that of \cite[Theorem 5.9]{SZ20}
but requires more delicate estimates of the error terms.
Because the Morawetz type functional $\mathscr{I}^{(2)}$ in \eqref{gen-energy-def} has been
already analyzed in \eqref{dtwtI.1},
we will focus on the first functional $\mathscr{I}^{(1)}$ in \eqref{gen-energy-def}.

We use Taylor's expansion and equation \eqref{equa-R} to compute
(see \cite[(5.33)]{SZ20})
\begin{align} \label{dt-mscrI}
 \frac{d \mathscr{I}^{(1)}}{dt}
=& -\sum_{k=1}^K \dot{\lbb_{k}}\lbb_{k}^{-3} {\rm Im}\int|R|^2\Phi_kdx
  -\sum_{k=1}^K \lbb_{k}^{-2} {\rm Im} \<f^\prime(U)\cdot R, R_k\>  \nonumber\\
  &  -{\rm Re} \< f''(U,R)\cdot R^2,  \partial_t {U} \>
    - \sum_{k=1}^K \lbb_{k}^{-2}{\rm Im} \< R \nabla\Phi_k, \na R\> \nonumber\\
  &-\sum_{k=1}^K \lbb_{k}^{-2} {\rm Im} \< f''(U,R)\cdot R^2, R_k\>
-{\rm Im} \bigg<\Delta R - \sum_{k=1}^K \lbb_{k}^{-2} R_k +f(u)-f(U), \eta\bigg>  \nonumber\\
=:&  \sum\limits_{l=1}^6 \mathscr{I}^{(1)}_{t,l},
\end{align}
where $f'(U)\cdot R$, $f''(U,R)\cdot R^2$ are given by \eqref{f'UR} and \eqref{f''UR-R2}, respectively.

Note that,
since $|\frac{\lbb_k \dot{\lbb_k} +\g_k}{\lbb_k^4}|\leq \frac{Mod}{\lbb_k^4} \leq C$
due to \eqref{Mk-Mod-t-4},
we compute
\begin{align} \label{mscrI.1}
  \mathscr{I}^{(1)}_{t,1}
  =  \sum_{k=1}^{K} \(\frac{\g_k}{\lbb_k^4}\int |R|^2\Phi_kdx-\frac{\lbb_k \dot{\lbb_k} +\g_k}{\lbb_k^4}\int |R|^2\Phi_kdx\)
  = \sum_{k=1}^{K} \frac{\g_k}{\lbb_k^4}\int |R|^2\Phi_kdx
     + \calo(D^2).
\end{align}

Moreover, we have from \cite[(5.35)]{SZ20} that
\begin{align} \label{I1t2t3-esti}
   \mathscr{I}^{(1)}_{t,2} +  \mathscr{I}^{(1)}_{t,3}
   =&- \sum\limits_{k=1}^K \frac{\gamma_k}{\lambda_k^2} {\rm Re}
      \int (1+\frac 2d)|U_k|^{\frac {4}{d}}|R_k|^2
     + \frac 2 d |U_k|^{\frac 4d-2} \ol{U_k}^2 R_k^2\  dx \nonumber \\
 &- \sum\limits_{k=1}^K\frac{\gamma_k}{\lambda_k} {\rm Re}
    \int \(\frac{x-\alpha_k}{\lambda_k}\)\cdot\nabla \ol{U_k}
           \ (f''(U_k)\cdot R_k^2) dx
       + \calo(E_2),
\end{align}
where $f''(U_k)\cdot R_k^2$ is defined as in \eqref{f''Qk-vek2}
with $U_k$ and $R_k$ replacing $Q_k$ and $\ve_k$, respectively,
and
\begin{align}
  E_2  \leq&\sum\limits_{k=1}^K \frac{|\beta_k|}{|\lbb_k|} \int |U_k|^{\frac 4d-1} |\na U_k| |R_k|^2 dx
            + C\sum\limits_{k=1}^K\int |U_k|^{\frac 4d} |R_k|^2 dx \nonumber \\
          & +C\frac{Mod(t)}{(T-t)^4}\|R\|^2_{L^2}
            + C\sum\limits_{l\geq 3} (T-t)^{\frac d2 l - 4 -d} \|R\|_{L^l}^l.
\end{align}
Using \eqref{R-D-Lp}, \eqref{gamj-wj-t}, Theorem \ref{Thm-D-t2} and Corollary \ref{Cor-Mk-Mod-4} we have
\begin{align} \label{E2-D2}
   E_2 \leq& C \((T-t)^{-2} D^2
              + (T-t)^{-\frac d2 (\frac 4d-1)} D^2
            + \sum\limits_{l= 3}^{1+\frac 4d} (T-t)^{\frac d2 l - 4-d} (T-t)^{-d(\frac l2-1)} D^l \)  \nonumber \\
      \leq& C  \frac{D^2}{(T-t)^{2}}.
\end{align}

Regarding $\mathscr{I}^{(1)}_{t,4}$,
we consider {\rm Case (I)} and {\rm Case (II)} separately.
In {\rm Case (I)},
since $\sum_{k=1}^K \na \Phi_k = 0$,
as in \eqref{lbbk-wt2},
we bound $\mathscr{I}^{(1)}_{t,4}$ by
\begin{align}
   |\mathscr{I}^{(1)}_{t,4}|
   =& \bigg| \( \frac{1}{\lbb_k^2} - \frac{1}{\omega^2(T-t)^2} \)
      {\rm Im} \<R\na \Phi_k, \na R\>  \bigg|  \nonumber \\
   \leq& \sum\limits_{k=1}^K \frac{|\omega- \omega_k||\omega+\omega_k|}{\omega^2\omega_k^2 (T-t)^2} \|R\|_{L^2} \|\na R\|_{L^2}
   \leq C\ve \frac{D^2}{(T-t)^3},
\end{align}
where we have used \eqref{lbbk-t-approx} and the inequality
$|\omega - \omega_k| \leq \ve$ in {\rm Case (I)}.

In {\rm Case (II)},
since $\|\na \Phi\|_{L^\9} \leq C \sigma^{-1} \leq C \ve$,
we  bound it easily  by
\begin{align}
    |\mathscr{I}^{(1)}_{t,4}|
    \leq \sum\limits_{k=1}^K  \frac{1}{\lbb_k^2} \|\na \Phi_k\|_{L^\9} \|R\|_{L^2} \|\na R\|_{L^2}
    \leq C \ve \frac{D^2}{(T-t)^3}.
\end{align}

Moreover, using \eqref{R-D-Lp}, \eqref{f''UR-R2-esti} and \eqref{D-t-2} we infer that
\begin{align}
   |\mathscr{I}^{(1)}_{t,5}|
   \leq& C(T-t)^{-2} \(\int |U|^{\frac 4d-1}|R|^3 dx + \|R\|_{L^{2+\frac 4d}}^{2+\frac 4d} \) \nonumber \\
   \leq& C(T-t)^{-2} \((T-t)^{-2}  D^3 + (T-t)^{-2} D^{2+\frac 4d}\)
   \leq  C \frac{D^2}{(T-t)^{2}}.
\end{align}

Finally, regarding $\mathscr{I}^{(1)}_{t,6}$,
since
\begin{align*}
   |f(u) - f(U)| \leq C (|U|^{\frac 4d} + |R|^{\frac 4d})|R|
   \leq C((T-t)^{-2}+ |R|^{\frac 4d})|R|,
\end{align*}
using the integration by parts formula and \eqref{eta-Mod-esti} we have
\begin{align*}
  |\mathscr{I}^{(1)}_{t,6}|
  \leq& C(\|\na \eta\|_{L^2} \|\na R\|_{L^2} + (T-t)^{-2}\|R\|_{L^2}\|\eta\|_{L^2} + \|R\|_{H^1}^{1+\frac 4d}\|\eta\|_{L^2}) \nonumber \\
  \leq& C\(\frac{Mod}{(T-t)^4} D+ \(\frac{D}{T-t}\)^{1+\frac 4d} \frac{Mod}{(T-t)^2}+e^{-\frac{\delta}{T-t}}\).
\end{align*}
Note that, by \eqref{Mod-bdd} and \eqref{D-t-2},
\begin{align*}
   Mod  \leq C \(\sum\limits_{k=1}^K |M_k| + (T-t)^2 D  + e^{-\frac{\delta}{T-t}}\).
\end{align*}
Taking into account Theorem \ref{Thm-D-t2} we obtain
\begin{align} \label{mscrI.6}
    |\mathscr{I}^{(1)}_{t,6}|
    \leq C\(\frac{D^2}{(T-t)^2} + \sum\limits_{k=1}^K \frac{|M_k|}{(T-t)^2} +  e^{-\frac{\delta}{T-t}}\).
\end{align}

Now, combining  estimates \eqref{dt-mscrI}-\eqref{mscrI.6} altogether
and using $\ve_k$ in \eqref{Rj-ej}
we conclude that
\begin{align} \label{dt-mscrI-esti}
 \frac{d \mathscr{I}^{(1)}}{dt}
 =  &\sum_{k=1}^{K}\frac{\gamma_k}{\lambda_k^4} \|\ve_k\|_{L^2}^2
     -\sum_{k=1}^{K}\frac{\gamma_k}{\lambda_k^4} {\rm Re}\int (1+\frac 2d)|Q_k|^\frac 4d|\ve_k|^2
      +\frac 2 d   |Q_k|^{\frac 4d -2} \ol{Q_k}^2 \ve_k^2  dy   \nonumber  \\
&-\sum_{k=1}^{K}\frac{\gamma_k}{\lambda_k^4} {\rm Re}
      \int y \cdot\nabla \ol{Q_k}
        \ (f''(Q_k)\cdot \ve_k^2)  dy \nonumber \\
 &  +\calo\(\frac{D^2}{(T-t)^2} + \ve \frac{D^2}{(T-t)^3}
        + \sum\limits_{k=1}^K \frac{|M_k|}{(T-t)^2}
        +  e^{-\frac{\delta}{T-t}}\).
\end{align}

Therefore, combining \eqref{dtwtI.1} and \eqref{dt-mscrI-esti}
and arguing
as those below \eqref{dt-local-virial} we prove \eqref{dIt-mono-case1}.
\hfill $\square$

\subsection{Upgradation to the exponential order} \label{Subsec-Upgra-Exp}

We are now ready to fulfill the final upgradation to the exponential decay rate.
The main result is formulated in Theorem \ref{Thm-D-exp} below.

\begin{theorem}  \label{Thm-D-exp}
(Exponential decay of remainder)
There exist $C, \delta>0$ such that for $t$ close to $T$,
\begin{align} \label{R-P-exp}
   D(t)  \leq Ce^{-\frac{\delta}{T-t}}.
\end{align}
\end{theorem}

{\bf Proof.}
On one hand,
expanding the modified generalized energy $\mathscr{I}$ up to the second order yields
\begin{align} \label{gen-energy-expan-2order}
  \mathscr{I}
  = & \frac{1}{2} {\rm Re} \int |\nabla R|^2
         + \sum_{k=1}^K  \frac{1}{\lambda_{k}^2}|R|^2\Phi_k
         - (1+\frac 2d) |U|^{\frac 4d} |R|^2
         - \frac{2}{d}|U|^{\frac 4d -2} U^2 \ol{R}^2 dx
    +\calo\(\wh Er\),
\end{align}
where the error term
\begin{align}
  \wh Er=&\sum\limits_{l=3}^{2+\frac 4d} |U|^{2+\frac 4d-l}|R|^l dx
      + \|\na R\|_{L^2} \|R\|_{L^2}.
\end{align}

Using \eqref{R-D-Lp} and taking $t$ close to $T$ such that $T-t\leq \ve$ we get
\begin{align} \label{wtEr}
    \wh Er
     \leq& \sum\limits_{l=3}^{2+\frac 4d} (T-t)^{-\frac d2(2+\frac 4d -l)} (T-t)^{-d(\frac l2 -1)} D^l + \frac{D^2}{T-t}
     \leq  C \frac{D^2}{T-t}
     \leq C\ve \frac{D^2}{(T-t)^2}.
\end{align}

Moreover, arguing as in the proof of \eqref{Ev-quard-esti}
we deduce that the main order on the R.H.S. of \eqref{gen-energy-expan-2order} is bounded from below by
\begin{align} \label{D2t2-Mkt2}
   \tilde c \frac{D^2}{(T-t)^2} + \calo\(\sum\limits_{k=1}^K \frac{M_k^2}{(T-t)^2} + e^{-\frac{\delta}{T-t}}\).
\end{align}
Note that,
by \eqref{Mj-D} and \eqref{D-t-2},
for $t$ close to $T$ such that $T-t\leq \ve$,
\begin{align}
   \frac{M_k^2}{(T-t)^2}
   \leq& \frac{C}{(T-t)^2} \(\int_t^T \frac{D^2}{T-s}ds\)^2 + C e^{-\frac{\delta}{T-t}} \nonumber \\
   \leq& C (T-t)^4 \int_t^T \frac{D^2}{(T-s)^3} ds  + C e^{-\frac{\delta}{T-t}} \nonumber \\
   \leq& C \ve  \int_t^T \frac{D^2}{(T-s)^3} ds + C e^{-\frac{\delta}{T-t}}.
\end{align}

Thus,  taking $\ve$   small enough
we conclude from \eqref{gen-energy-expan-2order}, \eqref{wtEr} and \eqref{D2t2-Mkt2} that
\begin{align} \label{I-lowbdd}
   \mathscr{I}(t)
    \geq  \tilde c \frac{D^2(t)}{(T-t)^2}
         - C \ve \int_t^T \frac{D^2}{(T-s)^3} ds - C e^{-\frac{\delta}{T-t}}.
\end{align}

On the other hand,
for $t$ close to $T$ such that $T-t<\ve$, we have
\begin{align*}
  \frac{D^2}{(T-t)^2} \leq \ve \frac{D^2}{(T-t)^3}, \ \
  \frac{1}{(T-t)^2} \int_t^T \frac{D^2}{T-s}ds
  \leq \int_t^T \frac{D^2}{(T-s)^3} ds.
\end{align*}
Then, Theorem \ref{Thm-I-mono} and \eqref{Mj-D} yield that for
$t$ close to $T$,
\begin{align}  \label{dIdt-lowbdd}
\frac{d\mathscr{I}}{dt}
\geq  - C\(\varepsilon \frac{D^2(t)}{(T-t)^3}
     +  \int_t^T \frac{D^2}{(T-s)^3} ds  +e^{-\frac{\delta}{T-t}}\).
\end{align}

Thus, combining \eqref{I-lowbdd} and \eqref{dIdt-lowbdd}
and using the fundamental theorem of calculus
we obtain
\begin{align} \label{dIdt-esti-0}
    \tilde  c \frac{D^2(t)}{(T-t)^2}
 \leq  \mathscr{I}(\tilde t)
       + C\int_{t}^T \bigg(\varepsilon \frac{D^2(s)}{(T-s)^3}
       +  \int_s^T \frac{D^2}{(T-r)^3} dr \bigg) ds
        + C e^{-\frac{\delta}{T-t}}.
\end{align}
Since  by \eqref{D-t-2}, \eqref{gen-energy-expan-2order} and \eqref{wtEr},
\begin{align} \label{dIdt-esti-1}
    {\mathscr{I}}(\tilde  t)\leq C \frac{D^2(\tilde  t)}{(T-\tilde  t)^2} \leq C (T-\tilde  t)^2 \to 0,\ \ as\ \tilde  t\to T.
\end{align}
Letting $\tilde  t \to T$ in \eqref{dIdt-esti-0}
and using  estimate \eqref{dIdt-esti-1}
we obtain a new Gronwall type inequality
\begin{align} \label{D-iter-exp}
     \frac{D^2(t)}{(T-t)^2}
 \leq  C_2 \int_{t}^{T} \bigg(\varepsilon \frac{D^2(s)}{(T-s)^3}
      +   \int_s^T \frac{D^2}{(T-r)^3} dr \bigg)  ds
      +C_2 e^{-\frac{\delta}{T-t}}.
\end{align}
The keypoint here is that,
the second order $(T-t)^2$ in \eqref{D-iter-iint}
now has been replaced by the much faster exponential rate $e^{-\frac{\delta}{T-t}}$.

Now, we claim that
there exist $C, \delta >0$ such that for $t$ close to $T$,
\begin{align} \label{D-t-exp}
   D(t) \leq C e^{-\frac{\delta}{T-t}}.
\end{align}

In order to prove \eqref{D-t-exp},
we take $t$ close to $T$ such that
$T-t \leq \ve$ and
$e^{-\frac{\delta}{T-t}} \leq (T-t)e^{-\frac{\delta}{2(T-t)}}$.
Then, inserting \eqref{D-t-2} into \eqref{D-iter-exp} we get
\begin{align}  \label{D-iter-exp-1}
   \frac{D^2}{(T-t)^2}
   \leq C^2_2 \ve (T-t)^2
       + C^2_2 (T-t) e^{-\frac{\delta}{2(T-t)}}.
\end{align}

Again, we use the induction arguments and assume that for some $n\geq 1$,
\begin{align}  \label{D-iter-exp-n}
   \frac{D^2}{(T-t)^2}
   \leq \sum\limits_{k=0}^{n-1} 2^{k} C_2^{k+2} \ve^k (T-t) e^{-\frac{\delta}{2(T-t)}}
        + C_2^{n+1} \ve^n (T-t)^2.
\end{align}
Note that,
\eqref{D-iter-exp-n} holds at the initial step $n=1$
due to \eqref{D-iter-exp-1}.
Moreover,
plugging \eqref{D-iter-exp-n} into \eqref{D-iter-exp}
and using straightforward computations
we see that \eqref{D-iter-exp-n} is also valid
with $n+1$ replacing  $n$.

Thus, the induction arguments yield
that \eqref{D-iter-exp-n} holds for all $n\geq 1$.
In particular,
taking   $\ve$ small enough such that
$C_2 \ve \leq \frac 12$
we obtain \begin{align*}
   \frac{D^2(t)}{(T-t)^2}
   \leq \lim\limits_{n\to \9}
        \bigg( \(\sum\limits_{k=0}^{n-1} 2^kC_2^{k+1} \ve^k\)  (T-t)e^{-\frac{\delta}{2(T-t)}}
       + C_2^{n+1} \ve^{n}  (T-t)^{2} \bigg)
   \leq  2C_2(T-t)e^{-\frac{\delta}{2(T-t)}},
\end{align*}
which yields \eqref{D-t-exp},
thereby completing the proof of Theorem \ref{Thm-D-exp}.
\hfill $\square$

As a consequence of Theorems \ref{Thm-Mod}, \ref{Thm-Loc-Mass} and \ref{Thm-D-exp},
we obtain the exponential decay estimates of
localized mass and modulation equations,
which in particular improve the previous estimate \eqref{Mk-Mod-t-4}.

\begin{corollary} \label{Cor-Mk-Mod-exp}
(Exponential decay of localized mass and modulation equations)
There exist $C, \delta>0$ such that for $t$ close to $T$,
\begin{align} \label{Mk-Mod-exp}
   \sum\limits_{k=1}^K|M_k(t)| + Mod(t) \leq Ce^{-\frac{\delta}{T-t}}.
\end{align}
\end{corollary}

We conclude this section with the upgradation to the exponential decay rate of geometrical parameters.
Let us set
\begin{align} \label{calp-0k}
   \calp_{0,k}
   = (\lbb_{0,k}, \a_{0,k}, \beta_{0,k}, \g_{0,k}, \theta_{0,k})
   =(\omega_k(T-t), x_k, 0, \omega_k^2(T-t), \omega_k^{-2}(T-t)^{-1} +\vartheta_k),
\end{align}
which are the geometrical parameters corresponding to the pseudo-conformal blow-up solutions $S_k$,
$1\leq k\leq K$.
Let
$\calp_k=(\lbb_k, \a_k, \beta_k, \g_k, \theta_k)$ be the modulation parameters as in Theorem \ref{Thm-geometri-dec}
and set
\begin{align*}
   |\calp_k - \calp_{0,k}|:=  |\lbb_k-\omega_k(T-t)|+|\alpha_k(t)-x_k|+|\b_k(t)|
                              +|\g_k(t)-\omega_k^2(T-t)|+|\t_k(t)-\omega_k^{-2}(T-t)^{-1}-\vartheta_k|.
\end{align*}

\begin{corollary}  \label{Cor-P-exp}
(Exponential decay of geometrical parameters)
There exist $C, \delta>0$ such that for $t$ close to $T$,
\begin{align} \label{R-P-exp}
   \sum\limits_{k=1}^K |\calp_k(t)- \calp_{0,k}(t)| \leq Ce^{-\frac{\delta}{T-t}}.
\end{align}
\end{corollary}

{\bf Proof.}
{\it $(i)$ Estimates of $\lambda_k$ and $\g_k$.}
We compute that,
by \eqref{lbbk-t-approx}, \eqref{beta-gam-D-t-1} and \eqref{Mk-Mod-exp},
\begin{align} \label{dt-glbb*}
  \bigg|\frac{d}{dt}\(\frac{\gamma_{k}}{\lambda_{k}} \)\bigg|
 =\frac{|\lambda_{k}^2\dot{\gamma_{k}}-\lambda_{k}\dot{\lambda_{k}}\gamma_{k}|}{\lambda_{k}^3}
 \leq \frac{|\lambda_{k}^2\dot{\gamma}_{k} + \g^2_{k}|}{\lbb_k^3}
      + \bigg|\frac{\g_k}{\lbb_k}\bigg| \frac{|\g_k + \lambda_{k} \dot{\lbb_{k}}|}{\lambda_{k}^2}
\leq C \frac{Mod(t)}{\lambda_{k}^3}\leq C e^{-\frac{\delta}{3(T-t)}}.
\end{align}
Then, taking into account \eqref{glbb-omega-o1}
we infer that
\begin{align}  \label{gamlbb-1}
 \bigg| \(\frac{\gamma_k}{\lambda_k} \)(t)-\omega_k \bigg|
\leq\int_{t}^{T}  \bigg| \frac{d}{dr} \(\frac{\gamma_{k}}{\lambda_{k}} \) \bigg|dr
\leq C e^{-\frac{\delta}{3(T-t)}}.
\end{align}

Moreover,  by \eqref{gamlbb-1},
\begin{align*}
\bigg|\frac{d}{dt}(\lambda_{k} -\omega_k (T-t))\bigg|
 = \bigg| \dot{\lbb_k} + \frac{\g_k}{\lbb_k} + \omega_k - \frac{\g_k}{\lbb_k} \bigg|
 \leq\frac{Mod}{\lambda_{k}}+ C e^{-\frac{\delta}{3(T-t)}}
\leq C e^{-\frac{\delta}{3(T-t)}},
\end{align*}
which along with \eqref{lbb-a-t-0} yields that
\begin{align} \label{lbb-t-exp}
|\lambda_{k}(t) -\omega_k(T-t)|
\leq\int_{t}^{T}  \bigg|\frac{d}{dr}(\lambda_{k}-\omega_k(T-r))  \bigg|dr
\leq C  e^{-\frac{\delta}{3(T-t)}},
\end{align}
thereby yielding the estimate of $\lbb_k$ in \eqref{R-P-exp}.

Similarly,
by \eqref{D-t-exp}, \eqref{Mk-Mod-exp}  and \eqref{gamlbb-1},
\begin{align*}
   \bigg|\frac{d}{dt}(\g_k - \omega_k^2 (T-t))\bigg|
   = \bigg|\dot\g_{k} + \frac{\g_k^2}{\lbb_k^2} + \omega_k^2 - \frac{\g_k^2}{\lbb_k^2}\bigg|
   \leq  \frac{Mod}{\lbb_k^2} + C\bigg|\omega_k - \frac{\g_k}{\lbb_k}\bigg|
   \leq C e^{-\frac{\delta}{3(T-t)}}.
\end{align*}
Taking into account \eqref{gamj-wj-t} we obtain
\begin{align} \label{gam-t-exp}
   |\g_k(t) - \omega_k^2 (T-t)|
   \leq \int_t^{T}  \bigg|\frac{d}{dr} (\g_k(r) - \omega_k^2 (T-r)) \bigg| dr
   \leq C e^{-\frac{\delta}{3(T-t)}}.
\end{align}
which yields the estimate of $\g_k$ in \eqref{R-P-exp}.

{\it $(ii)$ Estimates of $\beta_k$ and $\alpha_k$.}
We use the refined estimate of $\beta_k$ in Corollary \ref{Cor-energy-refined} and \eqref{Mk-Mod-exp} to get
\begin{align}
     |\beta_{k}|^2
   \leq C (T-t)^2\sum\limits_{k=1}^K \bigg|\frac {\g_{k}}{ \lbb_{k}}-\omega_k\bigg| + C e^{-\frac{\delta}{T-t}}.
\end{align}
which along with  \eqref{gamlbb-1}
yields that
\begin{align} \label{b-t-exp*}
    |\beta_k(t)|
 \leq C e^{-\frac{\delta}{6(T-t)}},
\end{align}
thereby implying the estimate of $\beta_k$ in \eqref{R-P-exp}.

Moreover, by \eqref{Mk-Mod-exp} and \eqref{b-t-exp*},
\begin{align}
|\dot{\alpha}_{k}|=
\bigg|\frac{\lambda_k\dot{\alpha}_{k}-2\beta_{k}}{\lambda_{k}}+\frac{2\beta_{k}}{\lambda_{k}}\bigg|
\leq \frac{Mod}{\lambda_{k}}+\frac{2|\beta_{k}|}{\lambda_{k}}
\leq C e^{-\frac{\delta}{7(T-t)}}.
\end{align}
Hence,
taking into account $\lim_{t\to T}\a_k(t) = x_k$
we infer that
\begin{align} \label{a-t-exp}
|\alpha_{k}(t)-x_k|\leq\int_{t}^{T}|\dot{\alpha}_{k}(r)|dr
\leq C e^{-\frac{\delta}{7(T-t)}},
\end{align}
which yields the estimate of $\a_k$ in \eqref{R-P-exp}.

{\it $(iii)$ Estimate of $\theta_k$.}
By \eqref{Mk-Mod-exp}, \eqref{lbb-t-exp} and \eqref{b-t-exp*},
\begin{align}  \label{esti-thetan.0*}
   \bigg|\frac{d}{dt}(\theta_{k} - ({\omega_k^{-2}(T-t)^{-1}}+ \vartheta_{k}))\bigg|
=& \bigg|\frac{\lbb_{k}^2 \dot{\theta}_{k} -1 -|\beta_{k}|^2}{\lbb_{k}^2}
  +\frac{|\beta_{k}|^2}{\lbb_{k}^2}
  + \frac{1}{\lbb_{k}^2}
  - \frac{1}{\omega_k^2(T-t)^2}\bigg| \nonumber \\
\leq& \frac{Mod}{\lbb_k^2} + \frac{|\beta_k|^2}{\lbb_k^2} +
     \frac{|\lbb_{k}-\omega_k(T-t)||\lbb_{k}+\omega_k(T-t)|}{\omega_k^2\lbb_{k}^2(T-t)^2 }
   \leq  C e^{-\frac{\delta}{4(T-t)}}.
\end{align}
In view of \eqref{P-t-0},
$\lim_{t\to T}|\theta_k - (\omega_k^{-2}(T-t)^{-1} + \vartheta_k)| =0$,
we obtain that for $t$  close to $T$,
\begin{align}
 |\theta_{k}- ({\omega_k^{-2}(T-t)^{-1}} + \vartheta_{k})|
 \leq& \int_{t}^{T} \bigg| \frac{d}{dr}(\theta_k - ( {\omega_k^{-2}(T-r)^{-1}} +\vartheta_{k})) \bigg|dr \nonumber
 \leq C e^{-\frac{\delta}{4(T-t)}},
\end{align}
thereby proving the estimate of $\theta_k$ in \eqref{R-P-exp}.

Therefore, the proof of Corollary \ref{Cor-P-exp} is complete.
\hfill $\square$

\section{Proof of main results} \label{Sec-Proof-Main}

In this section, we prove the main results  in Theorems \ref{Thm-Uniq-Blowup},
\ref{Thm-Uniq-Solitons-H1} and \ref{Thm-Uniq-Solitons-Sigma}.

\subsection{Proof of Theorem \ref{Thm-Uniq-Blowup}}   \label{Subsec-Uni-Blowup}

By virtue of Theorem \ref{Thm-D-exp} and Corollary \ref{Cor-P-exp},
we have that
\begin{align} \label{v-Uk-Pk}
   \|v(t) - \sum\limits_{k=1}^K U_k(t)\|_{H^1}
   + \sum\limits_{k=1}^K |\calp_k(t) - \calp_{0,k}(t)| \leq C e^{-\frac{\delta}{T-t}},
\end{align}
where $\calp_k$ and $\calp_{0,k}$ are as in Corollary \ref{Cor-P-exp},
and $C, \delta>0$.

Moreover,
using similar computations as in \eqref{SL-UL-diff} and  \eqref{naUL-naSL-L2}
and using \eqref{v-Uk-Pk} we obtain
\begin{align} \label{Uk-Sk-H1-exp}
    \|\sum\limits_{k=1}^K U_k(t)- \sum\limits_{k=1}^K S_k(t)\|_{H^1}
   \leq C e^{-\frac{\delta}{T-t}}.
\end{align}

Thus, combining \eqref{v-Uk-Pk} and \eqref{Uk-Sk-H1-exp} altogether we obtain
that for $t$ close to $T$
\begin{align} \label{v-Sk-H1-exp}
   \|v(t) - \sum\limits_{k=1}^K S_k(t)\|_{H^1}
   \leq \|v(t) - \sum\limits_{k=1}^K U_k(t)\|_{H^1}
       + \| \sum\limits_{k=1}^K U_k(t) - \sum\limits_{k=1}^K S_k(t)\|_{H^1}
    \leq C e^{-\frac{\delta}{T-t}}.
\end{align}
In particular, this yields that for any $\zeta \in (0,1)$,
\begin{align}
   \|v(t) - \sum\limits_{k=1}^K S_k(t)\|_{H^1} \leq C (T-t)^{3+\zeta}.
\end{align}

Thus, by virtue of  Theorem \ref{Thm-Uniq-Blowup-3+},
we obtain the uniqueness of multi-bubble blow-up solutions
satisfying \eqref{v-S-H1-o1-Thm} and \eqref{nav-naS-iint-0+-Thm}.
In particular,
the unique multi-bubble blow-up solution coincides with the one constructed in \cite{SZ20}.
Thus, in view of \eqref{v-S-sig-exp+-Thm},
we obtain the exponential convergence \eqref{v-S-Sigma-exp-Thm}.

Therefore, the proof of Theorem \ref{Thm-Uniq-Blowup} is complete.  \hfill $\square$

\subsection{Proof of Theorems \ref{Thm-Uniq-Solitons-H1} and \ref{Thm-Uniq-Solitons-Sigma}.}  \label{Subsec-Uniq-Multi-Solitary}

In order to prove Theorem \ref{Thm-Uniq-Solitons-H1},
let us first derive the convergence rate of $u-\sum_{k=1}^K W_k$
in the pseudo-coformal space $\Sigma$ from that in the energy space $H^1$.

\begin{lemma} \label{Lem-uW-H1-Sigma}
Assume that $u$ is a multi-soliton to equation \eqref{equa-NLS} satisfying that
for some $\mu>2$,
\begin{align}\label{u-Wj-H1-3+}
\|u(s)-\sum_{k=1}^KW_k(s)\|_{H^1} = \calo\(\frac{1}{s^{\mu}}\),\ \ as\ s\to +\infty,
\end{align}
where $\{W_k\}$ are the solitary waves given by \eqref{Wj-soliton}
with distinct speeds $\{v_k\}$.
Then, we have
\begin{align} \label{u-Wj-Sigma-1+-H1}
\|x(u(s)-\sum_{k=1}^KW_k(s))\|_{L^2} = \calo\(\frac{1}{s^{\mu-2}}\),\ \ as\ s\to +\infty.
\end{align}
Moreover, in the single soliton case where $K=1$,
if \eqref{u-Wj-H1-3+} holds with $\mu>1$
and the propagation speed $v_1=0$, then we have
\begin{align} \label{u-Wj-Sigma-0+-H1}
\|x(u(s)-W_1(s))\|_{L^2} = \calo\(\frac{1}{s^{\mu-1}}\),\ \ as\ s\to +\infty.
\end{align}
\end{lemma}

{\bf Proof.}
Set $z:=u-W$ with $W:= \sum_{k=1}^K W_k$, $f(z)=|z|^{\frac 4d} z$.
Let $\varphi$ be a smooth radial cut-off function such that $\varphi(x)=|x|^2$ for $|x|\leq 1$,
$\varphi(x)=0$ for $|x|\geq 2$, $0\leq \vf\leq 1$,
and for a universal constant $C>0$
\begin{align}  \label{navf-vf}
|\nabla \varphi(x)|^2\leq C \varphi(x) .
\end{align}
Set $\varphi_A(x) :=A^2\varphi(\frac{x}{A})$ for $A>0$
and
\begin{align} \label{MA-def}
M_A(s):=\int \varphi_A|z(s)|^2dx.
\end{align}

By equations \eqref{equa-NLS} and \eqref{Wj-soliton},
$z$ satisfies the equation
\begin{align} \label{equa-z}
iz_s+\Delta z +f(z) +f(W,z)+H=0,
\end{align}
where $f(W,z)$ denotes the difference
\begin{align} \label{fWz-def}
   f(W,z) : = f(z+W)-f(W)-f(z),
\end{align}
and the remainder $H$ contains the remaining coupling terms
\begin{align}
   H =f(W)-\sum_{k=1}^{K}f(W_k).
\end{align}

Then, using the integration by parts formula we get
\begin{align}  \label{dtMA}
\frac{d}{ds}M_A =-2{\rm Im}\langle\nabla \varphi_A z,\nabla z\rangle+2{\rm Im}\langle \varphi_A z,f(W,z)+H\rangle.
\end{align}

Note that, by H\"older's inequality,  \eqref{u-Wj-H1-3+} and \eqref{navf-vf},
\begin{align} \label{navfz-naz}
|{\rm Im}\langle\nabla \varphi_A z,\nabla z\rangle|\leq C\|z\|_{H^1}M_A^{\half}\leq {C} \frac{M_A^{\half}}{s^{\mu}} ,
\end{align}
where $C$ is a universal positive constant independent of $s$ and $A$.

Moreover,
since
\begin{align}
   |f(W,z) | \leq C (|W|^{\frac 4d}|z| +|W||z|^{\frac 4d} ),
\end{align}
using the exponential decoupling between $W_j$ and $W_k$, $j\not =k$,
we infer that
\begin{align*}
|{\rm Im}\langle \varphi_A z,f(W,z)+H\rangle|
&\leq C\sum_{k=1}^{K}\int \varphi_A |z|(|W_k|^{\frac{4}{d}}|z| +|W_k||z|^{\frac{4}{d}}) dx+ C M_A^{\frac12} e^{-\delta s} \nonumber \\
&\leq C\sum_{k=1}^{K}M_A^{\frac12}\(\|\varphi_A^{\half}W_k^{\frac{4}{d}}\|_{L^4}\|z\|_{H^1}
+\|\varphi_A^{\half}W_k\|_{L^4}\|z\|_{H^1}^\frac{4}{d}+e^{-\delta s}\),
\end{align*}
where $C$ is a universal constant.
Note that, by \eqref{Wj-soliton} and the change of variables
\begin{align*}
\|\varphi_A^{\half}W_k\|_{L^4}+\|\varphi_A^{\half}W_k^{\frac{4}{d}}\|_{L^4}
\leq C(\|xW_k\|_{L^4}+\|xW_k^{\frac{4}{d}}\|_{L^4}) \leq Cs.
\end{align*}
Taking into account \eqref{u-Wj-H1-3+} we get
\begin{align} \label{vfzf-MA}
|{\rm Im}\langle \varphi_A z,f(W,z) +H \rangle| \leq  {C} \frac{M_A^{\frac12}}{s^{\mu-1}}.
\end{align}
Note that,
one decay rate is lost here.

Then, plugging \eqref{navfz-naz} and \eqref{vfzf-MA}   into \eqref{dtMA} we get
\begin{align}
\bigg|\frac{d}{ds}M_A \bigg|\leq  {C} \frac{M_A^\half}{s^{\mu-1}},
\end{align}
which yields that, for any $s, \tilde{s} \in \bbr^+$, $s<\tilde{s}$,
\begin{align}
|M_A^{\frac12}(\tilde s)-M_A^{\frac12}(s)|\leq C( \frac{1}{s^{\mu-2}}- \frac{1}{\tilde{s}^{\mu-2}}),
\end{align}
where $C>0$ is  independent of $s$ and $A$.
Since $\mu>2$, and by \eqref{u-Wj-H1-3+},
$M_A(\tilde{s}) \to 0$ as $\tilde{s}\to \9$,
we let $\tilde{s} \to \9$ and get
\begin{align}
|M_A^{\frac12}(s)|\leq   \frac{C}{s^{\mu-2}}.
\end{align}
Since $C$ is  independent of $s$ and $A$,
we may let $A \to \9$ and use Fatou's lemma to obtain
\begin{align}
\|xz(s)\|_{L^2}\leq  \frac{C}{s^{\mu-2}},
\end{align}
which  yields \eqref{u-Wj-Sigma-1+-H1}.

In the single soliton case $K=1$, if the speed $v_1=0$,
then we have the improved estimate
\begin{align*}
\|\varphi_A^{\half}W_1\|_{L^4}+\|\varphi_A^{\half}W_1^{\frac{4}{d}}\|_{L^4}\leq \|xW_1\|_{L^4}+\|xW_1^{\frac{4}{d}}\|_{L^4}\leq C,
\end{align*}
and thus the above arguments  yield
\begin{align}
|{\rm Im}\langle \varphi_A z,f(W,z) + H\rangle| \leq  {C} \frac{M_A^{\frac12}}{s^{\mu}},
\end{align}
which improves the previous estimate \eqref{vfzf-MA}.
Thus, arguing as above we get \eqref{u-Wj-Sigma-0+-H1}.
\hfill $\square$

{\bf Proof of Theorem \ref{Thm-Uniq-Solitons-H1}.}
First note that,
for any $v, u \in \Sigma$, $v= P_T(u)$ given by \eqref{pseu-conf-transf},
one has
\begin{align}
   & \|v(t)\|_{L^2} = \|u(\frac{1}{T-t})\|_{L^2}, \label{v-u-L2-pct} \\
   & \|xv(t)\|_{L^2}=(T-t)\|yu(\frac{1}{T-t})\|_{L^2}, \label{xv-xu-L2-pct} \\
   & \|\na v(t)\|_{L^2}\leq C\(\frac{1}{T-t}\|\nabla u(\frac{1}{T-t})\|_{L^2}+\|yu(\frac{1}{T-t})\|_{L^2}\). \label{v-u-H1-pct}
\end{align}

Now, let $u_i$ denote the multi-soliton to \eqref{equa-NLS}
satisfying \eqref{u-W-H1-2+-Thm} and
$W= \sum_{k=1}^K W_k$ be the sum of $K$ solitary waves.
Let
$v_i = P_T(u_i)$, $S=P_T(W)$ be the blow-up solutions
corresponding to $u_i$ and $W$, respectively,
$i=1,2$.
We infer from Lemma \ref{Lem-uW-H1-Sigma}, \eqref{u-W-H1-2+-Thm}, \eqref{v-u-L2-pct} and \eqref{v-u-H1-pct} that
\begin{align*}
  \|v_i(t)-S(t)\|_{H^1}
  \leq \frac{1}{(T-t)} \|(u_i-W)(\frac{1}{T-t})\|_{H^1}
        + \|y(u_i-W)(\frac{1}{T-t})\|_{L^2}
  \leq C (T-t)^{\zeta},
\end{align*}
which verifies the condition \eqref{v-S-H1-0+-Thm}.

Thus, by virtue of Theorem \ref{Thm-Uniq-Blowup},
we obtain that $v_1 \equiv v_2$,
which, via the pseudo-conformal transformation \eqref{pseu-conf-transf},
in turn yields that $u_1\equiv u_2$,
and thus the uniqueness of multi-solitions follows.
Moreover, estimate \eqref{u-W-Sigma-exp-Thm} follows from \eqref{v-S-Sigma-exp-Thm}
and estimates \eqref{v-u-L2-pct}-\eqref{v-u-H1-pct}.

Therefore, the proof of Theorem \ref{Thm-Uniq-Solitons-H1} is complete.
\hfill $\square$

Next, we prove Theorem \ref{Thm-Uniq-Solitons-Sigma}.
The  key ingredient is the monotonicity of virial functional in Lemma \ref{Lem-Mono-virial-soliton},
which enables to obtain an improved space time estimate of  gradient
in Corollary \ref{Cor-v-H1-Sig-2nu}.

\begin{lemma} \label{Lem-Mono-virial-soliton}
(Monotonicity of virial functional)
Assume that $u$ is a multi-soliton to \eqref{equa-NLS}
such that
\begin{align}\label{u-Wj-Sigma-12+*}
\|u(s) - \sum\limits_{k=1}^K W_k(s)\|_{\Sigma} = \calo \(\frac{1}{s^{\mu}}\),\ \ for\ s\ large\ enough,
\end{align}
where $\mu>0$, $\{W_k\}$ are the solitary waves given by \eqref{Wj-soliton}
with distinct speeds $\{v_k\}$.
Assume additionally that $v_k \not = 0$, $1\leq k\leq K$.
Let $I$ denote the virial functional
\begin{align}  \label{virial-soliton}
I:={\rm Im}\int x\cdot \nabla z \bar{z}dx,\ \ with\ z:= u-\sum_{k=1}^K W_k.
\end{align}
Then,
there exists $C>0$ such that for $s$ large enough,
\begin{align} \label{dt-I-soliton}
\frac{dI}{ds}\geq  \|\nabla z(s)\|_{L^2}^2-C\(\frac{1}{s^{2\mu+1}}+e^{-\delta s}\).
\end{align}
\end{lemma}

{\bf Proof.} Using equation \eqref{equa-z} we compute
\begin{align} \label{dt-I-Lamz-soliton}
\frac{dI}{ds}&=2{\rm Im}\int \Lambda z \bar{z}_sdx
=-2{\rm Re}\langle \Lambda z, \Delta z+f(z)\rangle-2{\rm Re}\langle \Lambda z, f(W,z)\rangle-2{\rm Re}\langle \Lambda z, H\rangle,
\end{align}
where $\Lambda := \frac d2 {\rm I_d} + x\cdot \na$
and $f(W,z)$ is as in \eqref{fWz-def}.

For the first term on the R.H.S. of \eqref{dt-I-Lamz-soliton},
using the integration by parts formula we compute
\begin{align} \label{dt-I-Lamz-esti.1*}
-2{\rm Re}\langle \Lambda z, \Delta z+f(z)\rangle=2\|\nabla z\|_{L^2}^2-\frac{2d}{d+2}\|z\|^{2+\frac{4}{d}}_{L^{2+\frac{4}{d}}}=4E(z).
\end{align}
Moreover, by the Gagliardo-Nirenberg inequality \eqref{G-N} in Appendix,
\begin{align}
E(z)=\frac{1}{2}\int|\nabla z|^2 dx-\frac{d}{2d+4}\int |z|^{2+\frac{4}{d}}dx
\geq \(\half-C\|z\|_{L^2}^{\frac 4d}\)\|\nabla z\|_{L^2}^2,
\end{align}
where $C>0$.
Then, in view of \eqref{u-Wj-Sigma-12+*},
we may take $s$ large enough such that $C\|z(s)\|_{L^2}^{\frac 4d}\leq \frac14$
and get
\begin{align}  \label{Ez-naz-soltion}
   E(z(s))\geq \frac14\|\nabla z(s)\|_{L^2}^2.
\end{align}
Plugging this into \eqref{dt-I-Lamz-esti.1*} we get
\begin{align} \label{dt-I-Lamz-esti.1}
   -2{\rm Re}\langle \Lambda z, \Delta z+f(z)\rangle
   \geq \|\nabla z(s)\|_{L^2}^2.
\end{align}

Regarding the second inner product on the R.H.S. of \eqref{dt-I-Lamz-soliton},
using the integration by parts formula we
move the derivative onto $W$ and compute
\begin{align}
{\rm Re}\langle \Lambda z, f(z+W)-f(W)-f(z)\rangle
=& {\rm Re}\langle x\cdot \nabla W, f(z+W)-f(W)-f^\prime (W)\cdot z\rangle \nonumber  \\
 &+ \calo\(\sum\limits_{l=3}^{1+\frac 4d} \int  |W|^{2+\frac 4d -l}|z|^l dx\),
\end{align}
which along with the exponential decoupling between different traveling waves yields that
\begin{align} \label{Lamz-fzW-esti.1}
|{\rm Re}\langle \Lambda z, f(z+W)-f(W)-f(z)\rangle|
\leq & C\sum_{k=1}^{K}\sum_{l=2}^{\frac4d}\int |x| |\nabla W_k| |W_k|^{1+\frac4d-l} |z|^l dx  \nonumber \\
     & + C\sum_{k=1}^{K}\sum_{l=3}^{1+\frac4d}\int  |W_k|^{2+\frac4d-l} |z|^l dx
      + Ce^{-\delta s}.
\end{align}

For  the quadratic terms above,
since
$|x -v_k s| \geq \frac{|v_k| s}{2} $
for $|x|\leq \frac{|v_k| s}{2}$,
and $\na W_k, W_k\in L^\9$,
we infer from the exponential decay property of $Q$ that
\begin{align}  \label{Lamz-fzW-esti.2}
      \int |x| |\nabla W_k| |W_k|^{\frac4d-1}|z|^2 dx
\leq& C\(\int_{|x|\geq \frac{|v_k|s}{2}} |x Q(\frac{x-v_ks}{\omega_k})z^2|dx +
          \int_{|x|\leq \frac{|v_k| s}{2}} |x Q(\frac{x-v_k s}{\omega_k})z^2|dx\) \nonumber \\
\leq& C\(\|xz\|_{L^2}\|z\|_{L^2(|x|\geq \frac{|v_k|s}{2})} + e^{-\delta s}\) \nonumber \\
\leq& C\(\frac{\|xz\|_{L^2}^2}{|v_k|s} + e^{-\delta s}\)
\leq C\(\frac{1}{s^{2\mu+1}} + e^{-\delta s}\),
\end{align}
where we also used \eqref{u-Wj-Sigma-12+*}
and $v_k\not =0$ in the last step, $1\leq k\leq K$.

Moreover, for the higher order terms with $l\geq 3$, similarly we have
\begin{align}
\int |x| |\nabla W_k| |W_k|^{1+\frac4d-1}|z|^l dx
      + \int |W_k|^{2+\frac4d-l}|z|^l dx
 \leq C\(\|z\|_{\Sigma}\|z\|^{l-1}_{L^{2(l-1)}(|x|\geq \frac{|v_k|s}{2})} + e^{-\delta s}\).
\end{align}
Note that, by the Gagliardo-Nirenberg inequality \eqref{G-N},
\begin{align}
\|z\|^{l-1}_{L^{2(l-1)}(|x|\geq \frac{|v_k| s}{2})}
\leq C\(\|z\|^{l-1+d-\frac{ld}{2}}_{L^{2}(|x|\geq \frac{|v_k| s}{3})} \|\nabla z\|^{\frac{ld}{2}-d}_{L^{2}(|x|\geq \frac{|v_k| s}{3})}
+ \(\frac 1s\)^{l-1}\|z\|^{l-1}_{L^{2}(|x|\geq \frac{|v_k| s}{3})}\)
\leq \frac{C}{s^{\mu+1}}.
\end{align}
This along with \eqref{u-Wj-Sigma-12+*} yields that, for $3\leq l\leq 1+\frac 4d$,
\begin{align}  \label{Lamz-fzW-esti.3}
\int |x| |\nabla W_k| |W_k|^{1+\frac4d-1}|z|^l dx
      + \int |W_k|^{2+\frac4d-l}|z|^l dx
\leq C \( \frac{1}{s^{2\mu+1}} + e^{-\delta s} \).
\end{align}

Thus, plugging \eqref{Lamz-fzW-esti.2} and \eqref{Lamz-fzW-esti.3} into \eqref{Lamz-fzW-esti.1} we obtain
\begin{align} \label{Lamz-fzW-esti}
|{\rm Re}\langle \Lambda z, f(z+W)-f(W)-f(z)\rangle|\leq C\( \frac{1}{s^{2\mu+1}} + e^{-\delta s}\).
\end{align}

For the last term on the R.H.S. of \eqref{dt-I-Lamz-soliton},
again using the exponential decoupling between different traveling waves we get
\begin{align}  \label{Lamz-H-esti}
|{\rm Re}\langle \Lambda z, H\rangle|\leq \|z\|_{L^2}\|\Lambda H\|_{L^2}\leq Ce^{-\delta s}.
\end{align}

Therefore, plugging \eqref{Ez-naz-soltion}, \eqref{Lamz-fzW-esti} and \eqref{Lamz-H-esti} into \eqref{dt-I-Lamz-soliton}
we prove \eqref{dt-I-soliton}
and finish the proof.
\hfill $\square$

\begin{corollary} \label{Cor-v-H1-Sig-2nu}
Assume that   $u$ is a multi-soliton to equation \eqref{equa-NLS}
such that
\begin{align}\label{u-Wj-Sigma-12+}
\|u(s) - \sum\limits_{k=1}^K W_k(s)\|_{\Sigma} = \calo \(\frac{1}{s^{\mu}}\),\ \ for\ s\ large\ enough,
\end{align}
where $\mu>0$ and $\{W_k\}$ are as in \eqref{u-Wj-Sigma-12+*}.
Then, we have
\begin{align} \label{v-Sj-H1-12+}
\int_{s}^{+\infty}\|\nabla u(r) - \sum\limits_{k=1}^K \na W_k (r)  \|_{L^2}^2dr = \calo\(\frac{1}{s^{2\mu}}\),\ \ for\ s\ large\ enough.
\end{align}

In particular, if $\mu=\frac{1}{2}+\zeta$, $\zeta\in(0,1)$, then we have
\begin{align} \label{v-Sj-H1-1+}
\int_{s}^{+\infty}\|\nabla u(r) - \sum\limits_{k=1}^K \na W_k (r)\|_{L^2}^2dr = \calo\(\frac{1}{s^{1+2\zeta}}\),\ \ for\ s\ large\ enough,
\end{align}
and for $v:=P_Tu$ and $S_k := P_T(W_k)$, $0<T<\9$,
\begin{align}\label{nawtR-aver-0+-soliton}
\frac{1}{T-t}\int_{t}^{T}\|\nabla v(r) - \sum\limits_{k=1}^K \na S_k (r)\|_{L^2}^2dr =\calo((T-t)^{2\zeta}), \quad for \ t\ close\ to\ T .
\end{align}
\end{corollary}

\begin{remark}
Compared with \eqref{u-Wj-Sigma-12+},
\eqref{v-Sj-H1-12+} allows to gain half more convergence rate along some sequence $\{s_n\}$ to infinity,
namely,
\begin{align} \label{v-Sj-H1-12nu+}
\|\nabla u(s_n) - \sum\limits_{k=1}^K \na W_k(s_n) \|_{L^2} = \calo\(\frac{1}{s_n^{\mu+ \frac 12}}\).
\end{align}
\end{remark}

{\bf Proof of Corollary \ref{Cor-v-H1-Sig-2nu}.}
Let $z:=u- \sum_{k=1}^K W_k$,
$\wt R:= v- \sum_{k=1}^K S_k$
and $I$ be the virial functional as in \eqref{virial-soliton}.
Note that, by \eqref{u-Wj-Sigma-12+},
\begin{align}  \label{I-2nu-soliton}
|I(s)|\leq \|z\|_{\Sigma}^2\leq \frac{1}{s^{2\mu}}.
\end{align}
Moreover, integrating \eqref{dt-I-soliton} from $s$ to $\tilde{s}$
and using the boundary estimate \eqref{I-2nu-soliton} we get
\begin{align}
\int_{s}^{\tilde{s}}\|\nabla z\|_{L^2}^2dr
\leq C|I(\tilde{s})-I(s)|+ C\int_{s}^{\tilde{s}}\frac{1}{r^{1+2\mu}} + e^{-\delta r}dr
\leq C\(\frac{1}{s^{2\mu}}+ \frac{1}{\tilde{s}^{2\mu}}\).
\end{align}
Letting $\tilde{s}$ tend to infinity we obtain \eqref{v-Sj-H1-12+}.

In particular, if $\mu = \frac {1}{2}+\zeta$, we have that for some $C>0$,
\begin{align}
\int_{s}^{+\infty}\|\nabla z(r)\|_{L^2}^2dr
\leq \frac{C }{s^{1+2\zeta}},
\end{align}
which yields that, for the  rescaled time $t$ defined by $s=\frac{1}{T-t}$,
\begin{align}  \label{naz-int-0+-soliton}
\int_{t}^{T}\frac{\|\nabla z(\frac{1}{T-r})\|_{L^2}^2}{(T-r)^2}dr\leq C (T-t)^{1+2\zeta}.
\end{align}
Since $\wt R = P_T (z)$,
by \eqref{v-u-H1-pct},
\begin{align}
\|\nabla \widetilde{R}(r)\|_{L^2}^2\leq C\( \frac{\|\nabla z(\frac{1}{T-r})\|_{L^2}^2}{(T-r)^2}
+ \|yz(\frac{1}{T-r})\|_{L^2}^2\).
\end{align}
Thus, by virtue of \eqref{u-Wj-Sigma-12+} with $\mu = \frac 12 +\zeta$ and \eqref{naz-int-0+-soliton} we arrive at
\begin{align}
\int_{t}^{T}\|\nabla\widetilde{R}\|_{L^2}^2dr
\leq C \int_{t}^{T}\frac{\|\nabla z(\frac{1}{T-r})\|_{L^2}^2}{(T-r)^2}dr
+ C \int_{t}^{T}\|z(\frac{1}{T-r})\|_{\Sigma}^2dr
\leq C(T-t)^{1+2\zeta},
\end{align}
which yields \eqref{nawtR-aver-0+-soliton},
thereby finishing the proof of Corollary \ref{Cor-v-H1-Sig-2nu}.
\hfill $\square$

We are now ready to prove Theorem \ref{Thm-Uniq-Solitons-Sigma}.

{\bf Proof of Theorem \ref{Thm-Uniq-Solitons-Sigma}.}
Let $z$ and $\wt R$ be as in Corollary \ref{Cor-v-H1-Sig-2nu}.
Using \eqref{v-u-L2-pct}, \eqref{v-u-H1-pct} and \eqref{u-W-Sigma-12+-Thm}
we have that for $t$ close to $T$,
\begin{align*}
   \|\wt R(t)\|_{L^2} + (T-t) \|\na \wt R(t)\|_{L^2}
   \leq  \|z(\frac{1}{T-t})\|_{H^1} + (T-t)  \|yz(\frac{1}{T-t})\|_{L^2}
   \leq   C(T-t)^{\frac {1}{2}+\zeta},
\end{align*}
which verifies the condition \eqref{v-S-H1-o1-Thm}.
Moreover,
the condition \eqref{nav-naS-iint-0+-Thm}
is now verified by \eqref{nawtR-aver-0+-soliton}.

Thus, by virtue of Theorem \ref{Thm-Uniq-Blowup},
we obtain the uniqueness of multi-bubble blow-up solutions,
which in turn yields the uniqueness of multi-solitons
via the pseudo-conformal transformations.
Therefore, the proof of Theorem \ref{Thm-Uniq-Solitons-Sigma} is complete.
\hfill $\square$

At last, we end this section with the proof of Corollary \ref{Cor-Uniq-Soliton-Single}.

{\bf Proof of Corollary \ref{Cor-Uniq-Soliton-Single}.}
The uniqueness in \eqref{u-W-H1-2+-single-Thm}
is an application of Theorem \ref{Thm-Uniq-Solitons-H1} to the case $K=1$.
Moreover,
\eqref{u-W-H1-1+-single-Thm} can be proved by the improved convergence rate in \eqref{u-Wj-Sigma-0+-H1}
and similar arguments as in the proof of Theorem \ref{Thm-Uniq-Solitons-H1}.
Therefore, the proof of Corollary \ref{Cor-Uniq-Soliton-Single} is complete.
\hfill $\square$

\section{Appendix} \label{Sec-App}
In this appendix we collect the tools used in the previous sections
and the proof of modulation equations in Theorem \ref{Thm-Mod}.

\begin{lemma} \label{Lem-GN}
(\cite[Theorem 1.3.7]{C03})
Let $d\geq 1$ and $2\leq p< \infty$.
Then, there exists $C>0$ such that
\begin{align}  \label{G-N}
\|f\|_{L^p}\leq C \|f\|_{L^2}^{1 - d (\frac 12-\frac 1p)} \|\nabla f\|_{L^2}^{d(\frac 12-\frac 1p)}, \ \ \forall f\in H^1.
\end{align}
In particular,
for any $1<p<\9$,
\begin{align} \label{fLp-fH1}
   \|f\|_{L^p} \leq C \|f\|_{H^1},\ \ \forall f\in L^p.
\end{align}
\end{lemma}

\begin{lemma}(Decoupling estimates \cite[Lemma 3.1]{SZ20}) \label{Lem-inter-est}
For every $1\leq k\leq K$,
set
\begin{align} \label{Gj-gj-g}
    G_k(t,x)
    := \lbb_k^{-\frac d2} g_k(t,\frac{x-\a_k}{\lbb_k}) e^{i\theta_k}, \ \
    with\ \ g_{k}(t,y) := g(y) e^{i(\beta_{k}(t) \cdot y - \frac 14\g_{k}(t) |y|^2)},
\end{align}
where $g \in C_b^2(\bbr^d)$
decays exponentially fast at infinity
\begin{align*}
   |\partial^\nu g(y)| \leq C e^{-\delta |y|}, \ \ |\nu|\leq 2,
\end{align*}
with $C,\delta>0$,
$\calp_k:=(\lbb_k,\alpha_k,\beta_k,\gamma_k,\theta_k) \in C([T^*,T); \mathbb{X})$
satisfies that for $T^*$ close to $T$
\begin{align} \label{aj-xj}
  |(\alpha_k(t)-x_k)\cdot {\bf v_1}|\leq \sigma, \
  |x_k - \a_k(t)| \leq 1,\
  \half\leq\frac{\lambda_{k}(t)}{|\omega_{k}(T-t)|}\leq2, \ \ t\in[T^*,T),
\end{align}
and $|\beta_k|+|\g_k|\leq 1$,
\begin{align} \label{T-M-0}
    C (T-T^*) (1+ \max_{1\leq k\leq K} |x_k|) \leq  1,
\end{align}
where $C$ is sufficiently large but independent of $T$.
Then,
there exist $C, \delta>0$ such that
for any $1\leq k\not =l\leq K$,
$m\in \mathbb{N}$,
for any multi-index $\nu$ with $|\nu|\leq 2$
and for any $T^*$ close to $T$,
\begin{align} \label{Gj-Gl-decoup}
  \int\limits_{\bbr^d} |x-\a_l|^n |\partial^\nu G_l(t)| |x-\a_k|^m |G_k(t)| dx
   \leq Ce^{-\frac{\delta}{T-t}}, \ \ t\in [T^*,T).
\end{align}
Moreover, for any $h\in L^1$ or $L^2$,
$1\leq k\not = l\leq K$,
$m,n\in \mathbb{N}$,
multi-index $\nu$ with $|\nu|\leq 2$
and $T^*$ close to $T$,
\begin{align}  \label{Gj-hl-decoup}
   \int\limits_{\bbr^d} |x-\a_l|^n |\partial^\nu G_l(t)| |x-\a_k|^m |h| \Phi_kdx
   \leq Ce^{-\frac{\delta}{T-t}} \min\{\|h\|_{L^1}, \|h\|_{L^2}\}, \ \ t\in [T^*,T).
\end{align}
\end{lemma}

{\bf Coercivity of linearized operators.}
We recall the linearized operators from \cite{RS11,SZ19,SZ20,W85}.
Let $Q$ denote the ground state that solves  the elliptic equation \eqref{equa-Q}.
It is known that
$Q$ is smooth and decays at infinity exponentially fast,
i.e., there exist $C, \delta>0$ such that for any multi-index $|\nu|\leq 3$,
\be\label{Q-decay}
|\partial_x^\nu Q(x)|\leq C e^{-\delta |x|}, \ \ x\in \bbr^d.
\ee

Define the linearized operator $L=(L_+,L_-)$
around the ground state by
\begin{align} \label{L+-L-}
     L_{+}:= -\Delta + I -(1+{\frac{4}{d}})Q^{\frac{4}{d}}, \ \
    L_{-}:= -\Delta +I -Q^{\frac{4}{d}},
\end{align}
which has the generalized null space
spanned by $\{Q, xQ, |x|^2 Q, \na Q, \Lambda Q, \rho\}$.
Here, the operator
$\Lambda := \frac{d}{2}I_d + x\cdot \na$,
and $\rho$ is the unique radial solution to the equation
\begin{align} \label{rho-def}
L_{+}\rho= - |x|^2Q.
\end{align}
One has the exponential decay of $\rho$ (see, e.g., \cite{K-M-R, MP18})
\begin{align} \label{rho-decay}
|\rho(x)|+|\nabla \rho(x)|
\leq Ce^{-\delta|x|},
\end{align}
where $C,\delta>0$,
and the algebraic identities (see, e.g., \cite[(B.1), (B.10), (B.15)]{W85})
\be \ba \label{Q-kernel}
&L_+ \na Q =0,\ \ L_+ \Lambda Q = -2 Q,\ \ L_+ \rho = -|x|^2 Q, \\
&L_{-} Q =0,\ \ L_{-} xQ = -2 \na Q,\ \ L_{-} |x|^2 Q = - 4 \Lambda Q.
\ea\ee

For any complex valued $H^1$ function $f = f_1 + i f_2$
in terms of the real and imaginary parts,
set
\begin{align}
(Lf,f) :=\int f_1L_+f_1dx+\int f_2L_-f_2dx.
\end{align}
The scalar product  along the unstable directions
in the null space is defined by
\begin{align} \label{Scal-def}
{\rm Scal}(f)=\<f_1,Q\>^2+\<f_1,xQ\>^2+\<f_1,|x|^2Q\>^2+\<f_2,\nabla Q\>^2+\<f_2,\Lambda Q\>^2+\<f_2,\rho\>^2.
\end{align}

The localized coercivity of linearized operators is stated in Lemma \ref{Lem-coer-f-local} below.
\begin{lemma}(Localized coercivity \cite[Corollary 3.4]{SZ20})   \label{Lem-coer-f-local}
Let $\phi$ be a positive smooth radial function on $\R^d$,
such that
$\phi(x) = 1$ for $|x|\leq 1$,
$\phi(x) = e^{-|x|}$ for $|x|\geq 2$,
$0<\phi \leq 1$,
and $\lf|\frac{\nabla\phi}{\phi}\rt|\leq C$ for some $C>0$.
Set $\phi_A(x) :=\phi\lf(\frac{x}{A}\rt)$, $A>0$.
Then,
for $A$ large enough we have
\begin{align} \label{coer-f-local}
\int (|f|^2+|\nabla f|^2)\phi_A -(1+\frac 4d)Q^{\frac4d}f_1^2-Q^{\frac4d}f_2^2dx
\geq C_1\int(|f|^2+|\nabla f|^2)\phi_A dx-C_2 {\rm Scal}(f),
\end{align}
where $C_1, C_2>0$, and $f_1, f_2$ are the real and imaginary parts of $f$, respectively.
\end{lemma}

Below we prove Theorem \ref{Thm-Mod}.
For this purpose,
we set
$f(z):= |z|^{\frac 4d} z$, $d=1,2$,
and
\begin{align}
   & f'(U)\cdot R := \partial_z f(U) R + \partial_{\ol{z}} f(U) \ol{R}
      = (1+\frac 2d) |U|^\frac 4d R + \frac 2d |U|^{\frac 4d-2} U^2 \ol{R}, \label{f'UR} \\
   &f''(U)\cdot R^2
   := \frac 1d (1+\frac 2d) |U|^{\frac 4d -2} \ol U R^2
     + \frac 2d (1+\frac 2d) |U|^{\frac 4d -2}  U |R|^2
     + \frac 1d (\frac 2d-1) |U|^{\frac 4d -4}  U^3 \ol R^2, \label{f''z-R2} \\
   &f''(U,R)\cdot R^2:=R^2 \int_0^1 t \int _0^1 \partial_{zz} f(U+st R)ds dt
              + 2 |R|^2  \int_0^1 t \int _0^1 \partial_{z\ol{z}} f(U+st R)ds dt   \nonumber \\
             &\qquad \qquad \qquad \ \ + \ol{R}^2 \int_0^1 t \int _0^1 \partial_{\ol{z}\ol{z}} f(U+st R)ds dt.    \label{f''UR-R2}
\end{align}
Then,
we have the expansion
\begin{align} \label{f-expan}
   f(U+R) = f(U) + f'(U)\cdot R + f''(U,R)\cdot R^2.
\end{align}
We also get from equation \eqref{equa-NLS} and \eqref{v-dec} that
(see \cite[(4.11)]{SZ20})
\begin{align}\label{eq-U-R}
 &i\partial_tR
   +\sum_{k=1}^{K}
   (\Delta R_{k}+(1+\frac 2d)|U_{k}|^{\frac{4}{d}}R_{k}
   +\frac 2d|U_{k}|^{\frac{4}{d}-2}U_{k}^2 \ol{R_{k}}
   +  i\partial_tU_{k}+\Delta U_{k}+|U_{k}|^{\frac{4}{d}}U_{k})
   + \sum\limits_{k=1}^K f''(U_k) \cdot R^2  \nonumber \\
= &-H_1-H_2-H_3.
\end{align}
Here,
$H_1, H_2$ and $H_3$ contain the interactions between different blow-up profiles
\begin{align}
  H_1 : =& (1+\frac 2d)|U|^{\frac{4}{d}}R+\frac 2d|U|^{\frac{4}{d}-2}U^2\ol{R}
-\sum_{k=1}^{K}((1+\frac 2d)|U_{k}|^{\frac{4}{d}}R_{k}+\frac 2d|U_{k}|^{\frac{4}{d}-2}U_{k}^2\ol{{R}_{k}}), \label{G} \\
  H_2 : =&  |U|^{\frac{4}{d}}U-\sum_{k=1}^{K} |U_{k}|^{\frac{4}{d}}U_{k}, \label{H}\\
  H_3 : =&  f''(U,R)\cdot R^2-\sum_{k=1}^{K} f''(U_k)\cdot R^2. \label{H3}
\end{align}

Moreover,
we get from equation \eqref{equa-NLS} and \eqref{Uj-Qj-Q} that
(see \cite[(4.14)]{SZ20})
\begin{align}   \label{equa-Ut}
 i\partial_tU_{k}+\Delta U_{k}+|U_{k}|^{\frac{4}{d}}U_{k}
=&\frac{e^{i\theta_{k}}}{\lambda_{k}^{2+\frac d2}}
\bigg\{-(\lambda_{k}^2\dot{\theta}_{k}-1-|\beta_{k}|^2)Q_{k}
      -(\lambda_{k}^2\dot{\beta}_{k}+\gamma_{k}\beta_{k})\cdot yQ_{k}
     +\frac 14 (\lambda_{k}^2\dot{\gamma}_{k}+\gamma_{k}^2) |y|^2 Q_{k} \nonumber \\
  &\qquad \ \ -i(\lambda_{k}\dot{\alpha}_{k}-2\beta_{k})\cdot \nabla Q_{k}
     -i(\lambda_{k}\dot{\lambda_{k}}+\gamma_{k})\Lambda Q_{k} \bigg\}\(t,\frac{x-\alpha_{k}}{\lambda_{k}}\),
\end{align}
where $Q_k$ is given by \eqref{Qj-Q}, $1\leq k\leq K$.
The following identity also holds (see \cite[(4.28)]{SZ20})
\begin{align} \label{equa-Qk}
  \Delta Q_k - Q_k + |Q_k|^{\frac 4d} Q_k
  = |\beta_k - \frac{\g_k}{2}|^2 Q_k
     - i \g_k \Lambda Q_k + 2i \beta_k \cdot \na Q_k.
\end{align}

{\bf Proof of Theorem \ref{Thm-Mod}.}
The proof is similar to that of \cite[Proposition 4.3]{SZ20},
it mainly relies
on the  almost orthogonality in Lemma \ref{Lem-almost-orth}
and the decoupling Lemma \ref{Lem-inter-est}.

Taking the inner product of \eqref{eq-U-R} with
$\Lambda_k  {U_{k}}$ and then taking the real part
we get
\begin{align}\label{eq-Ul-Rl}
   &-{\rm Im}\langle\partial_tR,\Lambda_k U_{k}\rangle
   +{\rm Re}\langle\Delta R_{k}+(1+\frac 2d)|U_{k}|^{\frac{4}{d}}R_{k}+\frac 2d|U_{k}|^{\frac{4}{d}-2}U_{k}^2\ol{{R}_{k}},\Lambda_k U_{k}\rangle  \nonumber \\
   &+{\rm Re}\langle i\partial_tU_{k}+\Delta U_{k}+|U_{k}|^{\frac{4}{d}}U_{k},\Lambda_k U_{k}\rangle
     + {\rm Re} \<f''(U_k) \cdot R^2, \Lambda_k U_k\> \nonumber \\
 = &-{\rm Re} \bigg<\sum_{j\neq k} (\Delta R_{j}+(1+\frac 2d)|U_{j}|^{\frac{4}{d}}R_{j}+\frac 2d|U_{j}|^{\frac{4}{d}-2}U_{j}^2\ol{{R}_{j}}) + H_1,\Lambda_k U_{k} \bigg>   \nonumber \\
   & -{\rm Re} \bigg< \sum_{j\neq k}( i\partial_tU_{j}+\Delta U_{j}+|U_{j}|^{\frac{4}{d}}U_{j}) + H_2,\Lambda_k U_{k} \bigg>\nonumber \\
   &-{\rm Re} \bigg< \sum_{j\neq k}( f''(U_j)\cdot R^2 +H_3,\Lambda_k U_{k} \bigg>.
\end{align}

For the R.H.S. of equation \eqref{eq-Ul-Rl},
we claim that
\begin{align} \label{Mod-RHS-bdd}
  \lbb_k^2 \times({\rm R.H.S.\ of}\ \eqref{eq-Ul-Rl})
   = \calo\(e^{-\frac{\delta}{T-t}}(1+Mod+ \|R\|_{L^2}+ \|R\|^2_{L^2})
               + D^3(t) \).
\end{align}

To this end,
we have from \cite[(4.19), (4.20)]{SZ20} that the first two terms
on the R.H.S. of \eqref{eq-Ul-Rl} are bounded by
\begin{align} \label{Mod-RHS-bdd.1}
   C (T-t)^{-2} e^{-\frac{\delta}{T-t}} (1+ Mod+ \|R\|_{L^2}),
\end{align}
where $C,\delta>0$.
Moreover, By Lemma \ref{Lem-inter-est},
\begin{align} \label{Mod-RHS-bdd.2}
   \bigg|{\rm Re} \< \sum_{j\neq k}( f''(U_j)\cdot R^2 ,\Lambda_k U_{k} \> \bigg|
   \leq C (T-t)^{-2}  e^{-\frac{\delta}{T-t}} \|R\|^2_{L^2}.
\end{align}
Using \eqref{f''UR-R2} we also have
\begin{align} \label{Mod-RHS-bdd.3}
  |{\rm Re} \<H_3, \Lambda U_k\>|
    \leq&   C \sum\limits_{\tilde z, \tilde z^* \in \{z, \ol{z}\}}
     \bigg(\bigg<|R|^2 \int_0^1 t \int_0^1 |\partial_{\tilde z\tilde z^*} f(U) - \sum\limits_{k=1}^K \partial_{\tilde z\tilde z^*}f(U_k)|ds dt, |\Lambda_k U_k| \bigg> \nonumber \\
        & \qquad \qquad + \bigg< |R|^2 \int_0^1 t \int_0^1 |\partial_{\tilde z\tilde z^*} f(U+stR) - \partial_{\tilde z\tilde z^*}f(U)| ds dt, |\Lambda_k U_k|\bigg>\bigg).
\end{align}
Using Lemma \ref{Lem-inter-est} again
we see that the first inner product above only contributes
$e^{-\frac{\delta}{T-t}}\|R\|_{L^2}^2$,
while the second one can be bounded by
\begin{align} \label{Mod-RHS-bdd.4}
    \int |R|^3 (|U|^{\frac 4d -2} + |R|^{\frac 4d -2}) |\Lambda_k U_k| dx
    \leq& C (T-t)^{-2} (\|\ve_k\|_{L^3}^3 +\|\ve_k\|_{L^{1+\frac 4d}}^{1+\frac 4d})
    \leq C (T-t)^{-2} D^3,
\end{align}
where the renormalized remainder $\ve_k$ is given by \eqref{Rj-ej}
and we also used \eqref{fLp-fH1} in the last step.

Thus, combining  estimates \eqref{Mod-RHS-bdd.1}-\eqref{Mod-RHS-bdd.4}
and using \eqref{lbbk-t-approx}
we obtain \eqref{Mod-RHS-bdd}, as claimed.

Regarding the L.H.S. of \eqref{eq-Ul-Rl},
using the almost orthogonality \eqref{Orth-almost},
equation \eqref{equa-Ut} and \eqref{equa-Qk}
we have that
(see \cite[(4.26), (4.27), (4.31)]{SZ20})
\begin{align}
{\rm Im}\langle\partial_t R,\Lambda_k U_{k}\rangle
=& {\rm Im}\langle \Lambda_k R_{k},\partial_t U_{k}\rangle
      + \calo(e^{-\frac{\delta}{T-t}}(1+Mod) \|R\|_{L^2}),   \label{ptR-Ul} \\
 \lbb_k^2 {\rm Im}\langle \Lambda_k R_{k},\partial_t U_{k}\rangle
=& {\rm Re}\langle  \varepsilon_{k},\Lambda Q_{k}\rangle+ \g_k{\rm Im}\langle \Lambda \varepsilon_{k},\Lambda Q_{k}\rangle
   -2 \b_k{\rm Im}\langle \Lambda \varepsilon_{k},\nabla Q_{k}\rangle  \nonumber \\
&+ \calo((Mod+P^2)\|R\|_{L^2}),  \label{lef1}
\end{align}
and
\begin{align} \label{ptUj-DUj-fUj}
     \lbb_k^2 {\rm Re}\langle i\partial_tU_{k}+\Delta U_{k}+|U_{k}|^{\frac{4}{d}}U_{k},\Lambda_k U_{k}\rangle
   = - \frac{1}{4}\|yQ\|_2^2(\lbb_k^2\dot{\g_k}+\g_k^2)
     + \calo(|\beta_k|Mod).
\end{align}
Moreover, by \eqref{f''z-R2} and the change of variables,
\begin{align} \label{lef2}
  \lbb_k^2 {\rm Re}\langle f''(U_k)\cdot R^2,\Lambda_k U_{k}\rangle
=& {\rm Re}\int
   (1+\frac 2d)|Q_k|^{\frac {4}{d}}|\ve_k|^2
   + \frac 2d |Q_k|^{\frac 4d-2} Q_k^2 \ol{\ve_k}^2 dy  \nonumber \\
&+{\rm Re}\int (y \cdot\nabla \ol{Q_k}) (f''(Q_k)\cdot \ve^2_k) dy
  + \calo(e^{-\frac{\delta}{T-t}}\|R\|_{L^2}^2).
\end{align}

Thus,
combining \eqref{ptR-Ul}-\eqref{lef2} altogether
and rearranging the terms according to the orders of renormalized variable $\ve_k$
we get
\begin{align}
   & \lbb_k^2 \times ({\rm L.H.S.\ of}\ \eqref{eq-Ul-Rl})  \nonumber \\
=&- \frac{1}{4}\|yQ\|_2^2(\lbb_k^2\dot{\g_k}+\g_k^2)  \nonumber \\
 &  + {\rm Re}\langle\Delta \varepsilon_{k}-\varepsilon_k+(1+\frac 2d)|Q_{k}|^{\frac{4}{d}}\varepsilon_{k}+\frac 2d|Q_{k}|^{\frac{4}{d}-2}Q_{k}^2\ol{{\varepsilon}_{k}},\Lambda Q_{k}\rangle  \label{LHS-Modequa.1} \\
&- \g_k{\rm Im}\langle \Lambda \varepsilon_{k},\Lambda Q_{k}\rangle+2\b_k{\rm Im}\langle \Lambda \varepsilon_{k},\nabla Q_{k}\rangle
 \label{LHS-Modequa.2} \\
&+ {\rm Re}\int
   (1+\frac 2d)|Q_k|^{\frac {4}{d}}|\varepsilon_k|^2
   + \frac 2d  |Q_k|^{\frac 4d-2}Q_k^2 \ol{\varepsilon_k}^2 dy
    + {\rm Re}\int
   (y \cdot\nabla \ol{Q_k})(f''(Q_k)\cdot \ve_k^2) dy \label{LHS-Modequa.3} \\
&  + \calo\((Mod+P^2)\|R\|_{L^2} + P \ Mod + e^{-\frac{\delta}{T-t}}(1+Mod)\). \label{LHS-Modequa.4}
\end{align}

For the linear part on the R.H.S. above,
using the linearized operators $L_{\pm}$ in \eqref{L+-L-}
we have (see \cite[(4.37)]{SZ20})
\begin{align} \label{linear-Modequa}
    \eqref{LHS-Modequa.1}+ \eqref{LHS-Modequa.2}
    = {M_k } - \int |R |^2\Phi_kdx + \calo( (P^2+e^{-\frac{\delta}{T-t}}) \|R\|_{L^2}) .
\end{align}

Hence, we obtain that
\begin{align} \label{LHS-Modequa*}
  \lbb_k^2 \times ({\rm L.H.S.\ of}\ \eqref{eq-Ul-Rl})
  =& - \frac 14 \|yQ\|_{L^2}^2 ( \lbb_k^2 \dot{\g}_k + \g_k^2)
    + M_k - \int |R|^2 \Phi_k dx \nonumber \\
   &  + {\rm Re}\int
   (1+\frac 2d)|Q_k|^{\frac {4}{d}}|\varepsilon_k|^2
   + \frac 2d  |Q_k|^{\frac 4d-2}Q_k^2 \ol{\varepsilon_k}^2 dy    \nonumber \\
    &+ {\rm Re}\int
   (y \cdot\nabla \ol{Q_k}) (f''(Q_k)\cdot \ve_k^2) dy  \nonumber \\
   & + \calo\((P+\|R\|_{L^2}+e^{-\frac{\delta}{T-t}})Mod + P^2\|R\|_{L^2} + e^{-\frac{\delta}{T-t}}\),
\end{align}
which along with \eqref{R-t-0} and \eqref{Mod-RHS-bdd} yields \eqref{lbb2g+g2-f''Qve2}.

Moreover,
we can also bound the quadratic terms of $\ve_k$
in \eqref{LHS-Modequa.3}  by
\begin{align} \label{quadra-Modequa}
   \eqref{LHS-Modequa.3} \leq C \|\ve_k\|_{L^2}^2 \leq C\|R\|_{L^2}^2.
\end{align}
This along with \eqref{LHS-Modequa*} yields that
\begin{align} \label{Mod-LHS-bdd}
   \lbb_k^2 \times ({\rm L.H.S.\ of}\ \eqref{eq-Ul-Rl})
  =&- \frac{1}{4}\|yQ\|_2^2(\lbb_k^2\dot{\g_k}+\g_k^2)+ M_k  \nonumber \\
  &+ \calo\((P+\|R\|_{L^2}+e^{-\frac{\delta}{T-t}})Mod + P^2 \|R\|_{L^2} + \|R\|_{L^2}^2 + e^{-\frac{\delta}{T-t}}\).
\end{align}

Thus, combining \eqref{Mod-RHS-bdd} and \eqref{Mod-LHS-bdd} together
we obtain that for every $1\leq k\leq K$,
\begin{align}  \label{lbb2g+g2-Mod.0}
|\lbb_k^2\dot{\g_k}+\g_k^2|
   \leq& C \bigg( (P + \|R\|_{L^2}+e^{-\frac{\delta}{T-t}} )Mod  +|M_k|
         + P^2\|R\|_{L^2}+\|R\|_{L^2}^2
     + D^3  + e^{-\frac{\delta}{T-t}} \bigg),
\end{align}
which along with \eqref{D-def} and \eqref{D-0} yields that
\begin{align} \label{lbb2g+g2-Mod}
\sum_{k=1}^{K}|\lbb_k^2\dot{\g_k}+\g_k^2|
   \leq& C \bigg( (P+ \|R\|_{L^2}+e^{-\frac{\delta}{T-t}})Mod
           +\sum_{k=1}^{K} |M_k| +P^2 D +D^2
       + e^{-\frac{\delta}{T-t}} \bigg).
\end{align}

Similar arguments apply also to the remaining four modulation equations.
Actually, taking the inner products of equation (\ref{eq-U-R})
with $i(x-\a_k) U_k$, $i|x-\a_k|^2 {U_k}$, $\nabla {U_k}$, $\varrho_{k}$,
respectively,
then taking the real parts
and using analogous arguments as above, we obtain
\begin{align}\label{Mod-bdd-l}
Mod
\leq& C\bigg( (P+\|R\|_{L^2}+e^{-\frac{\delta}{T-t}})Mod
       +\sum_{k=1}^{K} |M_k| + P^2 D
+ D^2  + e^{-\frac{\delta}{T-t}} \bigg).
\end{align}

Therefore, using \eqref{lbb-a-t-0} and \eqref{P-t-0}
we may take $t$ close to $T$ such that
$P(t)+ \|R(t)\|_{L^2} + e^{-\frac{\delta}{T-t}}$ is sufficiently small
to obtain \eqref{Mod-bdd}.
The proof of Theorem \ref{Thm-Mod} is complete.
\hfill $\square$

\vspace{0.3cm}

\textbf{Acknowledgements:} D. Cao is supported by NNSF of China Grant 11831009,
Y. Su is supported by NSFC (No. 11601482) ,  D. Zhang  is supported by NSFC (No. 11871337) and Shanghai Rising-Star Program.

\end{document}